\documentclass[11pt]{article}
\usepackage{mathrsfs}
\usepackage{amssymb}
\usepackage{amsmath}
\usepackage{amsfonts}
\usepackage{amsthm}
\usepackage{mathrsfs}
\usepackage[top=0.9in,bottom=1.1in,left=0.8in,right=0.8in]{geometry}
\usepackage{lscape}
\usepackage{amsmath}
\usepackage{rotating}
\usepackage{slashbox}
\usepackage{setspace}
\usepackage{enumerate}
\usepackage{epsfig}
\usepackage{graphicx}
\usepackage{bm}
\usepackage{natbib}
\usepackage{subfigure}
\usepackage{amssymb}

\def\be{\begin{equation}}
\def\ee{\end{equation}}
\def\bea{\begin{eqnarray}}
\def\eea{\end{eqnarray}}
\def\bd{\begin{displaymath}}
\def\ed{\end{displaymath}}
\def\bda{\begin{eqnarray*}}
\def\eda{\end{eqnarray*}}

\def\ha1{\hat \beta_1}

\def\bsc{\begin{scriptsize}}
\def\esc{\end{scriptsize}}

\newcommand{\btheta} {\boldsymbol{\theta}}

\newtheorem{tm}{Theorem}
\newtheorem{la}{Lemma}
\newtheorem{cy}{Corollary}
\newtheorem{pn}{Proposition}

\theoremstyle{definition}

\begin{document}

\title{High Dimensional Generalized Empirical Likelihood for Moment
Restrictions with Dependent Data}
\author{~~~~~~~~~~{Jinyuan Chang\thanks{%
Department of Mathematics and Statistics, The University of Melbourne, Parkville, VIC, 3010, Australia }~~~~~~~~~~~~~~~~~~~~ Song Xi Chen\thanks{ Corresponding author. Guanghua School of Management and Center for Statistical Science, Peking University, Beijing, 100871, China. Tel.: +86 10 62760736; fax: +86 10 62760736. E-mail address: csx@gsm.pku.edu.cn}~~~~~~~~~~~~~~~
Xiaohong Chen\thanks{%
Department of Economics, Yale University, New Haven, CT, 06520, U.S.A.}~~~~~~~~~~} \\
\\
{The University of Melbourne~~~~~~~~ Peking University ~~~~~~~~~~~ Yale University~~~~~~~~}}
\date{First version: May 2013~~~~~ This version: October 2014}
\maketitle

\begin{abstract}
This paper considers the maximum generalized empirical likelihood (GEL)
estimation and inference on parameters identified by high dimensional moment
restrictions with weakly dependent data when the dimensions of the moment
restrictions and the parameters diverge along with the sample size. The
consistency with rates and the asymptotic normality of the GEL estimator are
obtained by properly restricting the growth rates of the dimensions of the
parameters and the moment restrictions, as well as the degree of data
dependence. It is shown that even in the high dimensional time series
setting, the GEL ratio can still behave like a chi-square random variable
asymptotically. A consistent test for the over-identification is proposed. A
penalized GEL method is also provided for estimation under sparsity setting.
%\footnote{Quite long the abstract. {\bf JY: I
%modified some thing, it seems be shorter} }
\end{abstract}



\section{Introduction}

In economic, financial and statistical applications, econometric models
defined with a growing number of parameters and moment restrictions are
increasingly employed. Vector autoregressive models, dynamic asset pricing
models, dynamic panel data models and high dimensional dynamic factor models
are specific examples; see, e.g., \cite{BaiNg_2002_Ecma}, \cite%
{StockWatson_2010} and \cite{FanLiao_2014}. Due to the desire to better
capture large scale dynamic fundamental relations, these models with large
number of unknown parameters of interest are typically used to for time
series data of high dimension due to a large number of variables (relative
to the sample size).

The unconditional moment restriction models are the inferential settings of
the Generalized Method of Moment (GMM) of \cite{Hansen_1982_Ecma}, which is
perhaps the most popular econometric method for semiparametric statistical
inference. There are two dimensions that play essential roles in this
method: the dimension of the moment restrictions and the dimension of the
unknown parameters of interest. When both dimensions are fixed and finite,
there is a huge established literature on inferential procedures, which
include but not restrict to \cite{Rothenberg_1973} for the minimum distance,
\cite{Hansen_1982_Ecma} and \cite{HansenSingleton_1982_Ecma} for the GMM,
\cite{Owen_1988_Biometrika}, \cite{QinLawless_1994_AOS} and \cite%
{Kitamura_1997_AOS} for the empirical likelihood (EL), \cite{Smith_1997},
\cite{NeweySmith_2004_Ecma} and {\cite{Anatolyev_2005} for the} generalized
empirical likelihood (GEL). Among these methods, some members of the GEL
(especially the EL) have the attractive properties of the Wilks' theorem %
\citep{Owen_1988_Biometrika,Owen_1990_AOS, QinLawless_1994_AOS}, Bartlett
correction \citep{ChenCui_2006_Biometrika, ChenCui_2007_JOE}, and a smaller
second order bias \citep{NeweySmith_2004_Ecma, Anatolyev_2005}. See \cite%
{Owen_2001}, \cite{Kitamura_2007} and \cite{ChenVanKeilegom_2009_Test} for
reviews.

This paper investigates high dimensional GEL estimation and testing for
weakly dependent observations when the dimensions of both the moment
restrictions and the unknown parameters of interest may grow with the sample
size $n$. Let $p$ and $r$ denote the dimension of the unknown parameters and
the number of moment restrictions, respectively. When $r\geq p$, we
investigate the impacts of $p$ and $r$ on the consistency, the rate of
convergence and the asymptotic normality of the GEL estimator, the limiting
behavior of the GEL ratio statistics as well as the overidentification test.
To accommodate the potential serial dependence in the estimating functions
induced by the original time series data, the blocking technique is
employed. This paper establishes the consistency (with rate) and the
asymptotic normality of the GEL estimator under either (fixed) finite or
diverging block size $M$ under some suitable restrictions on $r$, $p$, $M$
and $n$. It is demonstrated that in general the blocking technique with a
diverging block size delivers the estimation efficiency. We also discuss the
impact of the smallest eigenvalue of the covariance matrix of the averaged
estimating function on the consistency and the asymptotic normality of the
GEL estimator. We show that, even in high dimensional nonlinear time series
setting (with diverging $M$), the GEL ratio still behaves like a chi-square
random variable asymptotically, which echoes a similar result by \cite%
{FanZhangZhang_2001} for nonparametric regression with iid data. A GEL based
over-identification specification test is also presented for high
dimensional time series models, which extends that of \cite%
{DonaldImbensNewey_2003_JOE} for iid data from increasing dimension of
moments ($r$) but fixed finite dimension of parameters ($p)$ to both
dimensions are allowed to diverge (as long as $r-p>0$). Finally, when the
parameter space is sparse, a penalized GEL method is proposed to allow for $%
p>r$, and is shown to attain the oracle property in the selection
consistency as well as the asymptotic normality of the estimated non-zero
parameters.

There are some studies on the EL and its related methods under high
dimensionality of both the moment restrictions and the parameters of
interest. \cite{ChenPengQin_2009_Biometrika} and \cite%
{HjortMcKeagueVanKeilegom_2009_AOS} evaluated the EL ratio statistic for the
mean under high dimensional setting. \cite{TangLeng_2010_Biometrika} and
\cite{LengTang_2011_Biometrika} evaluated a penalized EL when the underlying
parameter is sparse in the context of the mean parameters and the general
estimating equations, respectively. \cite{FanLiao_2014} considered penalized
GMM estimation under high dimensionality and sparsity assumption. These
papers assume independent data. Recently, by allowing for dependent data but
losing the self-standardization property of the EL, \cite{Lahiri_2012}
proposed a modified EL method by adding a penalty term to the original EL
criterion for estimating the high-dimensional mean parameters with $r=p>n$.
\cite{Lahiri_2012} did not implement data blocking in their modified
(penalized) EL for the means despite the moment equations are serially
dependent. The EL ratio statistic based on their modified EL method is no
longer asymptotically pivotal. As a result, any inference based on this
modified EL has to use data blocking or other HAC long-run variance
estimation. The rationale in our paper is to preserve the attractive
self-standardization property of the GEL in high dimensional time series
setting; doing so makes our allowed dimensionality smaller than that in \cite%
{Lahiri_2012} but maintains simple GEL inference.

The rest of the paper is organized as follows. Section 2 introduces the high
dimensional model framework and the basic regularity conditions. Sections 3
and 4 establish the consistency, the rate of convergence and the asymptotic
normality of the GEL estimator. Sections 5 and 6 derive the asymptotic
properties of the GEL ratio statistic and the overidentification
specification test respectively. {Section 7 presents a penalized GEL
approach for parameter estimation and variable selection when the unknown
parameter is sparse}. Section 8 reports some simulation results and Section
9 briefly concludes. Technical lemmas and all the proofs are given in
Appendix.

\section{Preliminaries}

\subsection{Empirical Likelihood and Generalization}

Let $\{X_{t}\}_{t=1}^{n}$ be a sample of size $n$ from an $\mathbb{R}^{d}-$%
valued strictly stationary stochastic process, where $d$ denotes the
dimension of $X_{t}$, and $\boldsymbol{\theta }=(\theta _{1},\ldots ,\theta
_{p})^{\prime }$ be a $p$-dimensional parameter taking values in a parameter
space $\Theta $. Consider a sequence of $r$-dimensional estimating equation
\begin{equation*}
g(X_{t},\boldsymbol{\theta })=(g_{1}(X_{t},\boldsymbol{\theta }),\ldots
,g_{r}(X_{t},\boldsymbol{\theta }))^{\prime }
\end{equation*}%
for $r\geq p$. The model information regarding the data and the parameter is
summarized by moment restrictions
\begin{equation}
{E}\{g(X_{t},\boldsymbol{\theta }_{0})\}=\boldsymbol{\mathbf{0}}
\label{eq:1}
\end{equation}%
where $\boldsymbol{\theta }_{0}\in \Theta $ is the true parameter. As argued
in \cite{HjortMcKeagueVanKeilegom_2009_AOS}, the moment restrictions (\ref%
{eq:1}) can be viewed as a triangular array where $r,d,X_{t},\boldsymbol{%
\theta }$ and $g(x,\boldsymbol{\theta })$ may all depend on the sample size $%
n$. We will explicitly allow $r$ and/or $p$ grow with $n$ while considering
inference for $\boldsymbol{\theta }_{0}$ identified by (\ref{eq:1}).
Although there is often a connection between $d$ and $r$ which is dictated
by the context of an econometrical or statistical analysis, the theoretical
results established in this paper are written directly on the growth rates
of $r$ and $p$ relative to $n$. Hence, we will not impose explicit
conditions on $d$ which can be either growing or fixed. Certainly, when $d$
diverges, it would indirectly affect the underlying assumptions made in
Section 2.3, for instance the moment condition and the rate of the mixing
coefficients.

We assume the dependence in the time series $\{X_{t}\}$ satisfies the $%
\alpha $-mixing condition \citep{Doukhan_1994}. Specifically, let $%
\mathscr{F}_{u}^{v}=\sigma (X_{t}:u\leq t\leq v)$ be the $\sigma $-field
generated by the data from a time $u$ to a time $v$ for $v\geq u$. Then, the
$\alpha $-mixing coefficients are defined as
\begin{equation*}
\alpha _{X}(k)=\sup_{d}\sup_{A\in \mathscr{F}_{-\infty }^{0},B\in \mathscr{F}%
_{k}^{\infty }}|{P}(A\cap B)-{P}(A){P}(B)|~~\text{for each}~k\geq 1.
\end{equation*}%
The $\alpha $-mixing condition means that $\alpha _{X}(k)\rightarrow 0$ as $%
k\rightarrow \infty $. When $\{X_{t}\}$ are independent, $\alpha _{X}(k)=0$
for all $k\geq 1$.

We employ the blocking technique
\citep{Hall_1985,Carlstein_1986_AOS,
Kunsch_1989_AOS} to preserve the dependence among the underlying data. Let $%
M $ and $L$ be two integers denoting the block length and separation between
adjacent blocks, respectively. Then, the total number of blocks is $%
Q=\lfloor (n-M)/L\rfloor +1$, where $\lfloor \cdot \rfloor $ is the integer
truncation operator. For each $q=1,\ldots ,Q$, the $q$-th data block $%
B_{q}=(X_{(q-1)L+1},\ldots ,X_{(q-1)L+M})$. The average of the estimating
equation over the $q$-th block is
\begin{equation}
\phi _{M}\left( B_{q},\boldsymbol{\theta }\right) =\frac{1}{M}%
\sum\limits_{m=1}^{M}g(X_{(q-1)L+m},\boldsymbol{\theta }).  \label{eq:phi}
\end{equation}%
Clearly, ${E}\{\phi _{M}(B_{q},\boldsymbol{\theta }_{0})\}=\boldsymbol{%
\mathbf{0}}$. For any $n$ and $\boldsymbol{\theta }\in \Theta $, $\{\phi
_{M}(B_{q},\boldsymbol{\theta })\}_{q=1}^{Q}$ is a new stationary sequence.
The blockwise EL \citep{Kitamura_1997_AOS} is defined as
\begin{equation}
\mathcal{L}(\boldsymbol{\theta })=\sup \bigg\{\prod\limits_{q=1}^{Q}\pi _{q}%
\bigg{|}\pi _{q}>0,\sum\limits_{q=1}^{Q}\pi _{q}=1,\sum\limits_{q=1}^{Q}\pi
_{q}\phi _{M}(B_{q},\boldsymbol{\theta })=\boldsymbol{\mathbf{0}}\bigg\}.
\label{eq:l}
\end{equation}%
Employing the routine optimization procedure for the blockwise EL leads to
\begin{equation}
\mathcal{L}(\boldsymbol{\theta })=\prod\limits_{q=1}^{Q}\bigg\{\frac{1}{Q}%
\frac{1}{1+\widehat{\lambda }(\boldsymbol{\theta })^{\prime }\phi _{M}(B_{q},%
\boldsymbol{\theta })}\bigg\},  \label{eq:L}
\end{equation}%
where $\widehat{\lambda }(\boldsymbol{\theta })$ is a stationary point of
the function $q(\lambda )=-\sum_{q=1}^{Q}\log \{1+\lambda ^{\prime }\phi
_{M}(B_{q},\boldsymbol{\theta })\}$.

The EL estimator for $\boldsymbol{\theta}_0$ is $\widehat{\boldsymbol{\theta}%
}_{EL}=\arg\max_{\boldsymbol{\theta}\in\Theta}\log\mathcal{\ L}(\boldsymbol{%
\theta}).$
% See \cite{QinLawless_1994_AOS} and \cite{Kitamura_1997_AOS} for details.
The maximization in (\ref{eq:l}) can be carried out more efficiently by
solving the corresponding dual problem, which implies that $\widehat{%
\boldsymbol{\theta}}_{EL}$ can be obtained as
\begin{equation}
\widehat{\boldsymbol{\theta}}_{EL}=\arg\min_{\boldsymbol{\theta}%
\in\Theta}\max_{\lambda\in\widehat{\Lambda}_n(\boldsymbol{\theta}%
)}\sum_{q=1}^Q\log\{1+\lambda^{\prime }\phi_M(B_q,\boldsymbol{\theta})\},
\label{eq:el}
\end{equation}
where $\widehat{\Lambda}_{n}(\boldsymbol{\theta} )=\left\{ \lambda \in
\mathbb{R}^{r}:\lambda ^{\prime }\phi _{M}(B_q,\boldsymbol{\theta} )\in
\mathcal{V},q=1,\ldots ,Q\right\} $ for any $\boldsymbol{\theta} \in \Theta $
and $\mathcal{V}$ is an open interval containing zero.

The link function $\log (1+v)$ in (\ref{eq:el}) can be replaced by a general concave function $\rho (v)$ \citep{Smith_1997}. The domain of $\rho (\cdot )$ contains $0$ as an interior point, and $\rho (\cdot )$ satisfies $\rho _{v}(0)\neq 0$ and $\rho _{vv}(0)<0$ where $\rho _{v}(v)=\partial \rho (v)/\partial v$ and $\rho _{vv}(v)=\partial ^{2}\rho(v) /\partial v^{2}$. The GEL estimator \citep{Smith_1997, NeweySmith_2004_Ecma} is
\begin{equation}
\widehat{\boldsymbol{\theta }}_{n}=\arg \min_{\boldsymbol{\theta }\in \Theta
}\max_{\lambda \in \widehat{\Lambda }_{n}(\boldsymbol{\theta }%
)}\sum_{q=1}^{Q}\rho (\lambda ^{\prime }\phi _{M}(B_{q},\boldsymbol{\theta }%
)),  \label{eq:gel}
\end{equation}%
which includes the EL estimator $\widehat{\boldsymbol{\theta }}_{EL}$ of
\cite{Owen_1988_Biometrika}, the exponential tilting (ET) estimator of \cite%
{KitamuraStutzer_1997} and \cite{ImbensSpadyJonson_1998} (with $\rho
(v)=-\exp (v)$), the continuous updating (CU) GMM estimator of \cite%
{HansenHeatonYaron_1996} (with a quadratic $\rho (v)$), and many others as
special cases. Define
\begin{equation}
\widehat{S}_{n}(\boldsymbol{\theta },\lambda )=\frac{1}{Q}%
\sum\limits_{q=1}^{Q}\rho (\lambda ^{\prime }\phi _{q}(\boldsymbol{\theta }%
)).  \label{eq:sn}
\end{equation}%
Then $\widehat{\boldsymbol{\theta }}_{n}$ and its Lagrange multiplier $%
\widehat{\lambda }$ satisfy the score equation
\begin{equation*}
\nabla _{\lambda }\widehat{S}_{n}(\widehat{\boldsymbol{\theta }}_{n},%
\widehat{\lambda })=\boldsymbol{\mathbf{0}}.
\end{equation*}%
By the implicit function theorem [Theorem 9.28 of \cite{Rudin_1976}], for
all $\boldsymbol{\theta }$ in a $\Vert \cdot \Vert _{2}$-neighborhood of $%
\widehat{\boldsymbol{\theta }}_{n}$, there is a $\widehat{\lambda }(%
\boldsymbol{\theta })$ such that $\nabla _{\lambda }\widehat{S}_{n}(%
\boldsymbol{\theta },\widehat{\lambda }(\boldsymbol{\theta }))=\boldsymbol{%
\mathbf{0}}$ and $\widehat{\lambda }(\boldsymbol{\theta })$ is continuously
differentiable in $\boldsymbol{\theta }$. By the concavity of $\widehat{S}%
_{n}(\boldsymbol{\theta },\lambda )$ with respect to $\lambda $, $\widehat{S}%
_{n}(\boldsymbol{\theta },\widehat{\lambda }(\boldsymbol{\theta }%
))=\max_{\lambda \in \widehat{\Lambda }_{n}(\boldsymbol{\theta })}\widehat{S}%
_{n}(\boldsymbol{\theta },\lambda )$. From the envelope theorem,
\begin{equation}
\boldsymbol{\mathbf{0}}=\nabla _{\boldsymbol{\theta }}\widehat{S}_{n}(%
\boldsymbol{\theta },\widehat{\lambda }(\boldsymbol{\theta }))\big|_{%
\boldsymbol{\theta }=\widehat{\boldsymbol{\theta }}_{n}}=\frac{1}{Q}%
\sum_{q=1}^{Q}\rho _{v}(\widehat{\lambda }(\widehat{\boldsymbol{\theta }}%
_{n})^{\prime }\phi _{q}(\widehat{\boldsymbol{\theta }}_{n}))\{\nabla _{%
\boldsymbol{\theta }}\phi _{q}(\widehat{\boldsymbol{\theta }}_{n})\}^{\prime
}\widehat{\lambda }(\widehat{\boldsymbol{\theta }}_{n}).  \label{eq:0.1}
\end{equation}%
The role of the block size $M$ played in the consistency and the asymptotic
normality of the GEL estimator $\widehat{\boldsymbol{\theta }}_{n}$ will be
discussed in Sections 3 and 4, respectively.

\subsection{Examples}

We illustrate the model setting of high dimensional moment restrictions
framework through three examples.

%Examples 1 and 2 are two specific models and Example
%3 is a general model form which includes a wide range of specific
%models in econometrics, finance and statistics. The estimation of
%the parameter identified by the model proposed in Example 3 can be
%transferred to estimate the parameter defined by a new generalized
%estimating equation follows our setting. This shows our model
%setting is valuable.
\bigskip

\textbf{Example 1} (High dimensional means): Suppose $\{X_{t}\}_{t=1}^{n}$
is a stationary sequence of observations, where $X_{t}\in \mathbb{R}^{d}$
and $\boldsymbol{\theta }_{0}={E}(X_{t})$. For high dimensional data, $d$
diverges and $g(X_{t},\boldsymbol{\theta })=X_{t}-\boldsymbol{\theta }$
constitutes the simplest high dimensional moment equation, which implies the
dimension of observation $d$, the number of moment restrictions $r$ and the
number of parameters $p$ all are the same. Under this setting and for
independent data, \cite{ChenPengQin_2009_Biometrika} and \cite%
{HjortMcKeagueVanKeilegom_2009_AOS} considered the asymptotic normality of
the EL ratio, that mirrors the Wilks' theorem for finite dimensional case.

This framework can be used in other inference problems. For instance
checking if two univariate stationary time series $\{Y_{t}\}$ and $\{Z_{t}\}$
have identical marginal distribution. Let
\begin{equation*}
f_{Y}(s)=E(e^{\mathbf{i}sY_{t}})~~\text{and}~~f_{Z}(s)=E(e^{\mathbf{i}%
sZ_{t}})
\end{equation*}%
denote the characteristic functions of the two series, respectively. Suppose
all the moments of $Y_{t}$ and $Z_{t}$ exist, then the characteristic
functions can be expressed as
\begin{equation*}
f_{Y}(s)=1+\sum_{k=1}^{\infty }\frac{(\mathbf{i}s)^{k}}{k!}E(Y_{t}^{k})~~%
\text{and}~~f_{Z}(s)=1+\sum_{k=1}^{\infty }\frac{(\mathbf{i}s)^{k}}{k!}%
E(Z_{t}^{k}).
\end{equation*}%
Let $X_{t}=(Y_{t},Z_{t})$ and $g(X_{t},\boldsymbol{\theta }%
)=(a_{1}(Y_{t}-Z_{t}-\theta _{1}),\ldots ,a_{r}(Y_{t}^{r}-Z_{t}^{r}-\theta
_{r}))^{\prime }$ for some nonzero constants $a_{1},\ldots ,a_{r}$. Here $%
\theta _{l}$ measures $E(Y_{t}^{l})-E(Z_{t}^{l})$ for $l=1,\ldots ,r$, and
the $a_{i}$'s are used to account for the potential diverging moments case,
i.e, either $E(Y_{t}^{l})$ or $E(Z_{t}^{l})$ may diverge as $l\rightarrow
\infty $. Then, the test for whether $Y_{t}$ and $Z_{t}$ having the same
marginal distribution can be conducted by testing if $\boldsymbol{\theta }%
_{0}=\mathbf{0}$ via the growing dimensional moment restrictions $E\{g(X_{t},%
\boldsymbol{\theta }_{0})\}=\boldsymbol{\mathbf{0}}$ by letting $%
r\rightarrow \infty $.

\bigskip

\textbf{Example 2} (Time series regression): We assume a structural model
for $s$-dimensional time series $Y_t$ which involve unknown parameter $%
\boldsymbol{\theta}\in\mathbb{R}^p$ of interest as well as time innovations
with unknown distributional form. Specifically, assume
\begin{equation}
h(Y_t,\ldots,Y_{t-m};\boldsymbol{\theta}_0)=\boldsymbol{\varepsilon}_t\in%
\mathbb{R}^r  \label{eq:hd}
\end{equation}
where $m\geq1$ is some constant. In this model, we can view $%
X_t=(Y_t^{\prime },\ldots,Y_{t-m}^{\prime })^{\prime }\in\mathbb{R}^d$ with $%
d=sm$ and $g(X_t,\boldsymbol{\theta})=h(Y_t,\ldots,Y_{t-m};\boldsymbol{\theta%
})$. If $E(\boldsymbol{\varepsilon}_t)=\boldsymbol{\mathbf{0}}$, it implies
\begin{equation*}
E\{g(X_t,\boldsymbol{\theta}_0)\}=\boldsymbol{\mathbf{0}}.
\end{equation*}
For conventional vector autoregressive models
\begin{equation}
Y_t=\mathbf{A}_1Y_{t-1}+\cdots+\mathbf{A}_mY_{t-m}+\eta_t  \label{eq:var}
\end{equation}
where $\mathbf{A}_1,\ldots,\mathbf{A}_m$ are some coefficient matrices
needed to be estimated and $\eta_t$ is the white noise series. This model is
the special case of (\ref{eq:hd}) with
\begin{equation*}
h(Y_t,\ldots,Y_{t-m};\boldsymbol{\theta}_0)=(Y_t-\mathbf{A}_1Y_{t-1}-\cdots-%
\mathbf{A}_mY_{t-m})\otimes(Y_t^{\prime },\ldots,Y_{t-m}^{\prime })^{\prime
}.  \label{eq:h}
\end{equation*}
In modern high dimensional time series analysis, we always assume the
dimensionality of $Y_t$ is large in relation to sample size, i.e., $%
s\rightarrow\infty$ as $n\rightarrow\infty$. Under such background, the
numbers of estimating equation and unknown parameters are both $s^2m$. If we
replace $(Y_t^{\prime },\ldots,Y_{t-m}^{\prime })^{\prime }$ by $%
(Y_t^{\prime },\ldots,Y_{t-m-l}^{\prime })^{\prime }$ for some fixed $l\geq1$%
, the model will be over-identified. The phenomenon of over-parametrization
in such model is well known \citep{Lutkepohl_(2006).}. \cite{Davis2012}
considered the estimation of (\ref{eq:var}) under the sparsity assumption on
$\mathbf{A}_i$'s. Under the sparsity, the penalized method proposed in
Section 7 can be applied. Some other models share the form (\ref{eq:hd}) can
be found in Section 3.1 of \cite{NordmanLahiri_2013}.

\bigskip

\textbf{Example 3} (Conditional moment restrictions): Let $%
\{X_{t}=(Y_{t}^{\prime },Z_{t}^{\prime })^{\prime }\}_{t=1}^{n}$ be a set of
observations, and $\rho (y,z,\boldsymbol{\theta })$ be a known $J$%
-dimensional vector of generalized residual function. The parameter $%
\boldsymbol{\theta }_{0}$ is uniquely defined via the following conditional
moment restrictions
\begin{equation}
{E}\{\rho (Y_{t},Z_{t},\boldsymbol{\theta }_{0})|Y_{t}\}=\boldsymbol{\mathbf{%
0}}~~\hbox{almost
surely}.  \label{eq:condm}
\end{equation}%
By different choices of the functional forms of the generalized residual
function $\rho (y,z,\boldsymbol{\theta })$, the conditional moment
restrictions (\ref{eq:condm}) include many existing models in statistics and
econometrics as special cases. The popular generalized linear models are
special cases of (\ref{eq:condm}). To appreciate this point, let $\mu (y)={E}%
(Z|Y=y)$ and $h(\mu (y))=y^{\prime }\boldsymbol{\theta }_{0}$ for an
increasing link function $h(\cdot )$. Then the generalized linear models are
special cases of (\ref{eq:condm}) with $\rho (y,z,\boldsymbol{\theta }%
_{0})=z-h^{-1}(y^{\prime }\boldsymbol{\theta }_{0})$.

Let $q^{K}(y)=(q_{1K}(y),\ldots ,q_{KK}(y))^{\prime }$ denote a $K\times 1$
vector of known basis functions that can approximate any square integrable
functions of $Y$ well as $K\rightarrow \infty $, such as polynomial splines,
B-splines, power series, Fourier series, wavelets, Hermite polynomials and
others; see, e.g., \cite{AiChen_2003_Ecma} and \cite%
{DonaldImbensNewey_2003_JOE}. Then, (\ref{eq:condm}) implies
\begin{equation}
{E}\{\rho (Y_{t},Z_{t},\boldsymbol{\theta }_{0})\otimes q^{K}(Y_{t})\}=%
\boldsymbol{\mathbf{0}}.  \label{eq:uncond}
\end{equation}%
Moreover, the unknown parameter $\boldsymbol{\theta }_{0}$ is a solution to
this set of increasing dimensional ($r=JK$) unconditional moment
restrictions (\ref{eq:uncond}). The dimension $K$ will increase with $n$ to
guarantee the consistency of the estimator for $\boldsymbol{\theta }_{0}$
and its asymptotic efficiency. Define $g(X_{t},\boldsymbol{\theta })=\rho
(Y_{t},Z_{t},\boldsymbol{\theta })\otimes q^{K}(Y_{t})$, then (\ref%
{eq:uncond}) is a special case of (\ref{eq:1}). The number of moment
restrictions $r=JK$ increases as $K$ does. For this model with iid data,
\cite{DonaldImbensNewey_2003_JOE} apply the GEL method to the increasing
number of the unconditional moment restrictions (\ref{eq:uncond}) to obtain
efficient estimation for finite fixed dimensional $\boldsymbol{\theta }_{0}$%
. They find that the diverging rate of the moment restrictions $r=JK$
depends on the choice of the basis functions $q^{K}(y)$. For example, if $%
q^{K}(y)$ is a spline basis then $r=JK$ could grow at the rate of $%
K=o(n^{1/3})$.

\subsection{Notations and Technical Conditions}

Throughout the paper, we use $C$s, with different subscripts, to denote
positive finite constants which does not depend on the sample size $n$. For
a matrix $A$, we use $\Vert A\Vert _{F}$ and $\Vert A\Vert _{2}$ to denote
its Frobenius-norm and operator-norm respectively, i.e., $\Vert A\Vert
_{F}=\left\{ \text{tr}\left( A^{\prime }A\right) \right\} ^{1/2}$ and $\Vert
A\Vert _{2}=\left\{ \lambda _{\text{max}}\left( A^{\prime }A\right) \right\}
^{1/2}$. If $a$ is a vector, $\Vert a\Vert _{2}$ denotes
%its Frobenius-norm and operator-norm are also
its $L_{2}$-norm. %-norm of vectors denoted by .
Without causing much confusion, we denote {the $i$-th component of $g(x,%
\boldsymbol{\theta})$ by $g_i(x,\boldsymbol{\theta})$; and simplify $g(X_t,%
\boldsymbol{\theta})$ and $\phi_M(B_q,\boldsymbol{\theta})$ by $g_t(%
\boldsymbol{\theta})$ %\footnote{do u use $g_{t
%i}(\theta)$ notation ? {\bf JY: Yes, see Assumption (A.2)(iii).}}}
and $\phi_q(\boldsymbol{\theta})$, % for any $\theta\in\Theta$,
respectively, where $\phi_M(B_q,\boldsymbol{\theta})$ is defined in (\ref%
{eq:phi}). Furthermore, we use $g_{t,j}(\boldsymbol{\theta})$ and $%
\phi_{q,j}(\boldsymbol{\theta})$ to denote the $j$-th component of $g_t(%
\boldsymbol{\theta})$ and $\phi_q(\boldsymbol{\theta})$ respectively.
% and use $\hat{\lambda}$ for the Lagrange multiplier corresponding to $%
%\hat{\theta}_{n,p}$.
Let %We also define %\footnote{So, $V_M = V_n$ if $M=n$ ?  {\bf Yes.}}
$\bar{g}(\boldsymbol{\theta})={n}^{-1}\sum_{t=1}^ng_t(\boldsymbol{\theta})$
and $\bar{\phi}(\boldsymbol{\theta})={Q}^{-1}\sum_{q=1}^Q\phi_q(\boldsymbol{%
\theta})$. Additionally, define
\begin{equation*}
V_M=\text{Var}\{M^{1/2}{\phi}_q(\boldsymbol{\theta}_0)\}~~\hbox{and}~~V_n=%
\text{Var}\{n^{1/2}\bar{g}(\boldsymbol{\theta}_0)\}
\end{equation*}
which are the covariance of the averaged estimating functions over a block
and the entire sample respectively. Clearly $V_M = V_n$ if $M=n$. The
following regularity conditions are needed in our analysis. }

\bigskip

(A.1) (i) $\{X_t\}$ is strictly stationary and there exists $\gamma >2$ such
that $\sum_{k=1}^{\infty }k\alpha _{X}(k)^{1-2/\gamma }<\infty$; (ii) $M\geq
L$ and $M/L\rightarrow c\geq 1$; (iii) ${E}\{g_t(\boldsymbol{\theta} _{0})\}=%
\boldsymbol{\mathbf{0}}$ and there are positive functions $\Delta _{1}(r,p)$
and $\Delta _{2}(\varepsilon )$ such that for any $\varepsilon >0$,
\begin{equation*}
\inf_{\{\boldsymbol{\theta} \in \Theta:\Vert \boldsymbol{\theta} -%
\boldsymbol{\theta} _{0}\Vert _{2}\geq \varepsilon \}}\Vert {E}\{g_t(%
\boldsymbol{\theta})\}\Vert _{2}\geq \Delta _{1}(r,p)\Delta _{2}(\varepsilon
)>0,
\end{equation*}%
where $\liminf_{r,p\rightarrow \infty }\Delta _{1}(r,p)>0$; (iv) $\sup_{%
\boldsymbol{\theta} \in \Theta}\Vert \bar{g}(\boldsymbol{\theta} )-{E}\{g_t(%
\boldsymbol{\theta} )\}\Vert _{2}=o_{p}\{\Delta _{1}(r,p)\}$.

\bigskip

(A.2) (i) $\boldsymbol{\theta} _{0}\in \text{int}(\Theta)$ and $\Theta$
contains a small $\Vert \cdot \Vert _{2}$-neighborhood of $\boldsymbol{\theta%
} _{0}$ in which $g(x,\boldsymbol{\theta} )$ is continuously differentiable
with respect to $\boldsymbol{\theta} $ for any $x\in \mathcal{X}$, the
domain of $X_t$, and
\begin{equation*}
\bigg|\frac{\partial g_{i}(x,\boldsymbol{\theta} )}{\partial \theta _{j}}%
\bigg|\leq T_{n,ij}(x)~~(i=1,\ldots,r;j=1,\ldots,p)
\end{equation*}
for some functions $T_{n,ij}(x)$ with $E\{T_{n,ij}^2(X_t)\}\leq C$ for any $%
i,j$; (ii) $\sup_{\boldsymbol{\theta} \in \Theta}\Vert g(x,\boldsymbol{\theta%
} )\Vert _{2}\leq r^{1/2}B_{n}(x)$, where ${E}\{B_{n}^{\gamma }(X_{t})\}\leq
C$ for $\gamma $ given in (A.1)(i); (iii) ${E}\{|g_{t,j}(\boldsymbol{\theta}
_{0})|^{2\gamma }\}\leq C$ for all $j=1,\ldots ,r$; (iv) the eigenvalues of $%
[E\{\nabla _{\boldsymbol{\theta} }{g}_t(\boldsymbol{\theta} )\}]^{\prime
}[E\{\nabla _{\boldsymbol{\theta} }{g}_t(\boldsymbol{\theta} )\}]$ in a $%
\|\cdot\|_2$-neighborhood of $\boldsymbol{\theta}_0$ are uniformly bounded
away from zero and infinity;
%\footnote{this is vague. which one ? and what is the meaning of uniform for the first part ? {\bf See above. Do you mean $\sup_{\theta \in \Theta}\lambda _{\text{%
%max}}\{n^{-1}\sum_{t=1}^{n}g_{t}(\theta )g_{t}^{\prime }(\theta
%)\}$?}}
$\sup_{\boldsymbol{\theta} \in \Theta}\lambda _{\text{max}%
}\{n^{-1}\sum_{t=1}^{n}g_{t}(\boldsymbol{\theta} )g_{t}(\boldsymbol{\theta}
)^{\prime }\}\leq C$ with probability approaching to 1.

\bigskip

(A.3) In a $\|\cdot\|_2$-neighborhood of $\boldsymbol{\theta}_{0}$, $g(x,%
\boldsymbol{\theta})$ is twice continuously differentiable with respect to $%
\boldsymbol{\theta}$ for any $x\in\mathcal{X}$, and for some functions $%
K_{n,ijk}(x)$ with $E\{K_{n,ijk}^2(X_t)\}\leq C$ for any $i,j,k$,
\begin{equation*}
\bigg{|}\frac{\partial ^{2}g_{i}(x,\boldsymbol{\theta} )}{\partial \theta
_{j}\partial \theta _{k}}\bigg{|}\leq
K_{n,ijk}(x)~~(i=1,\ldots,r;j,k=1,\ldots,p).
\end{equation*}

\bigskip

Condition (A.1)(i) specifies the rate of decay for the mixing coefficients
via a tuning parameter $\gamma $ as commonly assumed in the analysis of
weakly dependent data. When the data are independent, $\alpha _{X}(k)=0$ for
all $k\geq 1$ and this condition is automatically satisfied for any $\gamma
>2$. \cite{Kitamura_1997_AOS} assumed $\sum_{k=1}^{\infty }\alpha
_{X}(k)^{1-2/\gamma }<\infty $ for fixed finite dimensional EL, which
implies $M^{-1}\sum_{k=1}^{M}k\alpha _{X}(k)^{1-2/\gamma }\rightarrow 0$ by
Kronecker's lemma. In the current high dimensional setting, we need stronger
condition on the mixing coefficients in order to control remainder terms
when analyzing the asymptotic properties of the GEL estimator and the GEL
ratio. If $\{X_{t}\}$ is exponentially strong mixing \citep{FanYao_2003} so
that $\alpha _{X}(k)\sim \varrho ^{k}$ for some $\varrho \in (0,1)$, then
(A.1)(i) is automatically valid for any $\gamma >2$. (A.1)(ii) imposes a
condition regarding the two blocking quantities $M$ and $L$, which is
commonly assumed in the works of block bootstrap and blockwise EL.
(A.1)(iii) is the population identification condition for the case of
diverging parameter space. A similar assumption can be found in \cite%
{Chen_2007} and \cite{ChenPouzo_2012_Ecma}. The last part of (A.1) is an
extension of the uniform convergence. If $p$ is fixed, under the assumption
of the compactness of $\Theta $ and some other regularity conditions,
following \cite{Newey_1991}, $\sup_{\boldsymbol{\theta }\in \Theta }\Vert
\bar{g}(\boldsymbol{\theta })-{E}\{g_{t}(\boldsymbol{\theta })\}\Vert
_{2}=o_{p}(1)$ which is a special case of (A.1)(iv) with $\Delta _{1}(r,p)$
being a constant.

As conditions (A.1)(iii) and (iv) are rather abstractive, we illustrate them
via the examples given in Section 2.2. For Example 1, we can choose $%
\Delta_1(r,p)=1$ and $\Delta_2(\varepsilon)=\varepsilon$. For the
conditional moment restrictions model (Example 3), a common assumption in
the literature is that for any $a(Y_t)$ with $E\{a^2(Y_t)\}<\infty$ there
exists a $K\times1$ vector $\gamma_K$ such that $E[\{a(Y_t)-\gamma_K^{\prime
K}(Y_t)\}^2]\rightarrow0$ as $K\rightarrow\infty$. For any $\boldsymbol{%
\theta}\in\{\boldsymbol{\theta}:\|\boldsymbol{\theta}-\boldsymbol{\theta}%
_0\|_2\geq\varepsilon\}$, let $\Gamma_K(\boldsymbol{\theta})$ satisfy $%
E[\|E\{\rho(Y_t,Z_t,\boldsymbol{\theta})|Y_t\}-\Gamma_K(\boldsymbol{\theta}%
)q^K(Y_t)\|_2^2]\rightarrow0$ as $K\rightarrow\infty$. If $%
\sup_{y}\|E\{\rho(Y_t,Z_t,\boldsymbol{\theta})|Y_t=y\}-\Gamma_K(\boldsymbol{%
\theta})q^K(y)\|_2=O(K^{-\lambda})$ for some $\lambda>1/2$, then
\begin{equation*}
\begin{split}
\|E\{g_t(\boldsymbol{\theta})\}\|_2\geq&~ \big\|E[\{\Gamma_K(\boldsymbol{%
\theta}) q^K(Y_t)\}\otimes q^K(Y_t)]\big\|_2-\big\|E[\{\rho(Y_t,Z_t,%
\boldsymbol{\theta})-\Gamma_K(\boldsymbol{\theta})q^K(Y_t)\}\otimes q^K(Y_t)]%
\big\|_2 \\
\geq&~\lambda_{\min}\big(E\{q^K(Y_t)q^K(Y_t)^{\prime }\}\big)\|\Gamma_K(%
\boldsymbol{\theta})\|_F-O(K^{-\lambda})\big(\text{tr}[E\{q^K(Y_t)q^K(Y_t)^{%
\prime }\}]\big)^{1/2}.
\end{split}%
\end{equation*}
Under the assumption that the eigenvalues of $E\{q^K(Y_t)q^K(Y_t)^{\prime
}\} $ are uniformly bounded away from zero and infinity, we have
\begin{equation*}
\begin{split}
\inf_{\{\boldsymbol{\theta} \in \Theta:\Vert \boldsymbol{\theta} -%
\boldsymbol{\theta} _{0}\Vert _{2}\geq \varepsilon \}}\|E\{g_t(\boldsymbol{%
\theta})\}\|_2\geq&~ C\bigg[\inf_{\{\boldsymbol{\theta} \in \Theta:\Vert
\boldsymbol{\theta} -\boldsymbol{\theta} _{0}\Vert _{2}\geq \varepsilon
\}}\|\Gamma_K(\boldsymbol{\theta})\|_F-K^{1/2-\lambda}\bigg] \\
\geq&~C\bigg[\inf_{\{\boldsymbol{\theta} \in \Theta:\Vert \boldsymbol{\theta}
-\boldsymbol{\theta} _{0}\Vert _{2}\geq \varepsilon \}}E\big[\big\|%
E\{\rho(Y_t,Z_t,\boldsymbol{\theta})|Y_t\}\big\|_2\big]-K^{1/2-\lambda}\bigg]%
. \\
\end{split}%
\end{equation*}
Hence, as $\boldsymbol{\theta}_0$ is the unique root of $E\{\rho(Y_t,Z_t,%
\boldsymbol{\theta})|Y_t\}=\mathbf{0}$, (A.1)(iii) holds provided that the
lower bound in the above inequality is greater than or equal to $\Delta
_{1}(r,p)\Delta _{2}(\varepsilon )$. In addition, if the generalized
residual function $\rho(y,z,\boldsymbol{\theta})$ is continuously
differentiable with respect to $\boldsymbol{\theta}$. Then,
\begin{equation*}
\big\|E\{\rho(Y_t,Z_t,\boldsymbol{\theta})|Y_t\}\big\|_2\geq\|\boldsymbol{%
\theta}-\boldsymbol{\theta}_0\|_2\lambda_{\min}^{1/2}\big(E\big[\{\nabla_{%
\boldsymbol{\theta}} \rho(Y_t,Z_t,\boldsymbol{\theta}^*)\}^{\prime
}\{\nabla_{\boldsymbol{\theta}} \rho(Y_t,Z_t,\boldsymbol{\theta}^*)\}|Y_t%
\big]\big)
\end{equation*}
where $\boldsymbol{\theta}^*$ is on the line joining $\boldsymbol{\theta}_0$
and $\boldsymbol{\theta}$. If the eigenvalues of $E[\{\nabla_{\boldsymbol{%
\theta}} \rho(Y_t,Z_t,\boldsymbol{\theta})\}^{\prime }\{\nabla_{\boldsymbol{%
\theta}} \rho(Y_t,Z_t,\boldsymbol{\theta})\}|Y_t]$ are uniformly bounded
away from zero, $\Delta_1(r,p)$ and $\Delta_2(\varepsilon)$ can be chosen as
some constant $C$ and $\varepsilon$, respectively.

Condition (A.2)(i) assumes that the first derivatives of $g_{i}(x,%
\boldsymbol{\theta })$ near $\boldsymbol{\theta }_{0}$ are uniformly bounded
by some functions which have bounded second moments. (A.2)(ii) generalizes
the moment conditions on $g(x,\boldsymbol{\theta })$ for fixed dimensional
case \citep{QinLawless_1994_AOS,Kitamura_1997_AOS,NeweySmith_2004_Ecma}.
More generally, we can replace the factor $r^{1/2}$ by some function $\zeta
(r)>0$. We let $\zeta (r)=r^{1/2}$ to simplify the presentation. (A.2)(iii)
is the moment assumption on each $g_{t,j}(\boldsymbol{\theta }_{0})$. The
first part of (A.2)(iv) is an extension of that assumed in the EL or the GEL
in the fixed dimensional case %
\citep{QinLawless_1994_AOS,Kitamura_1997_AOS,NeweySmith_2004_Ecma}. The
second one of (A.2)(iv) is to bound $\sup_{\boldsymbol{\theta }\in \Theta
}\lambda _{\max }\{Q^{-1}\sum_{q=1}^{Q}\phi _{q}(\boldsymbol{\theta })\phi
_{q}(\boldsymbol{\theta })\}$. Our proofs in Appendix can be easily extended
to allow for $\sup_{\boldsymbol{\theta }\in \Theta }\lambda _{\max
}\{Q^{-1}\sum_{q=1}^{Q}\phi _{q}(\boldsymbol{\theta })\phi _{q}(\boldsymbol{%
\theta })^{\prime }\}$ diverging in probability. Note that we do not assume
the eigenvalues of $V_{M}$ or $V_{n}$ being bounded away from zero and
infinity, but rather leave it open for specific treatments in Sections 3 and
4 for the consistency and the asymptotic normality of the GEL estimator. Our
subsequent analysis shows that, to obtain the main results of the paper, $%
\lambda _{\min }(V_{M})$ is allowed to decay to zero at certain rates by
properly restricting the diverging rates of $r$ and $p$. Condition (A.3)
ensures the second derivatives of $g_{i}(x,\boldsymbol{\theta })$ near $%
\boldsymbol{\theta }_{0}$ are uniformly bounded by functions which have
bounded second moments.

\section{Consistency and Convergence Rates}

To study the consistency of the GEL estimator $\widehat{\boldsymbol{\theta }}%
_{n}$ defined by (\ref{eq:gel}), we need the following conditions regarding
the dimensionality $r$, the block size $M$ and the sample size $n$:%
\begin{equation}
r^{2}M^{2-2/\gamma }n^{2/\gamma -1}=o(1)~~\text{and}~~r^{2}M^{3}n^{-1}=o(1).
\label{eq:cond1}
\end{equation}

\begin{tm}
Assume conditions \textrm{(A.1)}, \textrm{(A.2)} and that the eigenvalues of
$V_M$ are uniformly bounded away from zero and infinity. Then, if \textrm{(%
\ref{eq:cond1})} holds, $\|\widehat{\boldsymbol{\theta}}_{n}-\boldsymbol{%
\theta}_{0}\|_2\xrightarrow{p}0. $ If in addition, $r^2pM^2n^{-1}=o(1)$,
then $\|\widehat{\boldsymbol{\theta}}_{n}-\boldsymbol{\theta}%
_{0}\|_2=O_p(r^{1/2}n^{-1/2})$ and $\|\widehat{\lambda}(\widehat{\boldsymbol{%
\theta}}_{n})\|_2=O_p(r^{1/2}Mn^{-1/2})$.
\end{tm}

%We note that for independent observations $\alpha _{X}(k) \equiv 0$
%for all $k\geq1$
This theorem provides the consistency of the GEL estimator $\widehat{%
\boldsymbol{\theta }}_{n}$ for both independent and dependent data when the
blocking size $M$ is either finite or diverging.
%\footnote{It also contains fixed and diverging $M$ case, right ? or do you need conditions on the EVs of $V_n$ ? {\bf Yes, it includes the results for both fixed and diverging $M$. We impose the conditions on $V_M$, so we do not need the conditions on $V_n$.}}
%The difference between independent data and dependent data can be reflected by the condition imposed on the eigenvalues of $V_M$.
For independent data, $V_{M}=E\{g_{1}(\boldsymbol{\theta }_{0})g_{1}(%
\boldsymbol{\theta }_{0})^{\prime }\}$ for any $M\geq 1$. Thus, to make $r$
have a faster diverging rate, we select the block size $M=1$. For dependent
data,
\begin{equation*}
V_{M}=E\{g_{1}(\boldsymbol{\theta }_{0})g_{1}(\boldsymbol{\theta }%
_{0})^{\prime }\}+\sum_{k=1}^{M-1}\bigg(1-\frac{k}{M}\bigg)\big[E\{g_{1}(%
\boldsymbol{\theta }_{0})g_{1+k}(\boldsymbol{\theta }_{0})^{\prime
}\}+E\{g_{1+k}(\boldsymbol{\theta }_{0})g_{1}(\boldsymbol{\theta }%
_{0})^{\prime }\}\big].
\end{equation*}%
However, if $\{g_{t}(\boldsymbol{\theta }_{0})\}_{t=1}^{n}$ is a martingale
difference sequence, then $V_{M}\equiv E\{g_{1}(\boldsymbol{\theta }%
_{0})g_{1}(\boldsymbol{\theta }_{0})^{\prime }\}$ for any $M\geq 1$ and $M=1$
should be used to make $r$ have a faster diverging rate. Furthermore, if the
eigenvalues of $V_{M}$ are uniformly bounded away from zero and infinity for
some fixed $M$, (\ref{eq:cond1}) is simplified to
\begin{equation*}
r^{2}n^{2/\gamma -1}=o(1).
\end{equation*}%
Here $\gamma $ determines the number of moments of the estimating equation
as specified in (A.2)(ii) and (A.2)(iii). Then, $r=o(n^{1/2-1/\gamma })$
ensures the consistency of the GEL estimator $\widehat{\boldsymbol{\theta }}%
_{n}$. For large enough $\gamma $, $r$ will be made close to $o(n^{1/2})$,
which is the best rate we can established.

%For dependent data, Condition (A.1)(i) implies that
%$\sum_{k=1}^{\infty} k\alpha_X(k)^{1-2/\gamma}$ is bounded, which
%together with the third part of (\ref{eq:cond1}) implies that
%$r=o(M)$. % which means $M$ should diverge faster than $r$.
% Then,
%$\|\widehat{\btheta}_n-\btheta_0\|_2=O_p(r^{1/2}n^{-1/2})$ provided
%that $r^2M^{2-2/\gamma}n^{2/\gamma-1}=o(1)$, $r^2M^3n^{-1}=o(1)$ and
%$rM^{-1}=o(1)$. More explicitly,
%\begin{equation*}
%M=O(n^{\{(\gamma-2)/(4\gamma-2)\}\wedge1/5 }) \quad \hbox{and} \quad
%r=o(n^{\{ (\gamma-2)/(4\gamma-2)\} \wedge {1/5}}).
%\end{equation*}
%%These mean that the block size $M$ and the number of estimating
%%equations $r$ has to be smaller order than $n^{1/5}$.
%%which are both determined by $\gamma$.
%If $\gamma\geq8$%
%, $M=O(n^{1/5})$ and $r=o(n^{1/5})$. We see here a substantial
%slowdown in the growth rate of $r$ under the dependence as compared
%to $o(n^{1/2 - 1/\gamma}) $ for the independence case.
%% as elaborated above.

Theorem 1 encompasses the existing consistency results for the GEL estimator
in the literature. Indeed, if $r$ is fixed and the data are independent,
Theorem 1 implies that both $\Vert \widehat{\boldsymbol{\theta }}_{n}-%
\boldsymbol{\theta }_{0}\Vert _{2}$ and $\Vert \widehat{\lambda }(\widehat{%
\boldsymbol{\theta }}_{n})\Vert _{2}$ are $O_{p}(n^{-1/2})$, which are the
same as the rates obtained in \cite{QinLawless_1994_AOS} for the EL and \cite%
{NeweySmith_2004_Ecma} for the GEL. If $r$ is fixed but the data are
dependent, Theorem 1 means that $\Vert \widehat{\boldsymbol{\theta }}_{n}-%
\boldsymbol{\theta }_{0}\Vert _{2}=O_{p}(n^{-1/2})$ and $\Vert \widehat{%
\lambda }(\widehat{\boldsymbol{\theta }}_{n})\Vert _{2}=O_{p}(Mn^{-1/2})$,
which coincides with the result of \cite{Kitamura_1997_AOS} for the EL
estimator. If $r$ is diverging and the data are independent, both $\Vert
\widehat{\boldsymbol{\theta }}_{n}-\boldsymbol{\theta }_{0}\Vert _{2}$ and $%
\Vert \widehat{\lambda }(\widehat{\boldsymbol{\theta }}_{n})\Vert _{2}$ are $%
O_{p}(r^{1/2}n^{-1/2})$, which retain the results in \cite%
{DonaldImbensNewey_2003_JOE} and \cite{LengTang_2011_Biometrika}.

%As indicated at the start of this section, the blocking is not
%necessary to establish the consistency of estimation despite the
%underlying dependence in the data. Indeed,
%%if $
%%\lambda_{\min}(E\{g_t(\btheta_0)g_t(\btheta_0)'\})$ is uniformly
%%bounded away from zero,
%  the following corollary gives the
%consistency of the GEL estimator  without the blocking ($M=1$).
%
%
%%\begin{cy}
%%Under conditions (A.1) (i) and (iii) and (A.2), if $
%%\lambda_{\min}(E\{g_t(\btheta_0)g_t(\btheta_0)'\})$ is uniformly
%%bounded away from zero and $r^2n^{2/\gamma-1}=o(1)$, then the GEL
%%estimator with $M=1$ satisfies
%%$\|\widehat{\btheta}_{n}-\btheta_{0}\|_2\xrightarrow{p}0.$ In
%%addition, if
%%$r^2pn^{-1}=o(1)$,  then $\|\widehat{\btheta}_{n}-\btheta_{0}%
%%\|_2=O_p(r^{1/2}n^{-1/2})$.
%%\end{cy}
%%
%%The proof of Corollary 1 can be made by repeating the proof of
%%Theorem 1 given in the Appendix. We note that under the assumption
%%that $\lambda_{\min}(E\{g_t(\btheta_0)g_t(\btheta_0)'\})$ is
%%uniformly bounded away from zero, Lemmas 7 and 8 in Appendix still
%%hold without the last restriction of (\ref{eq:cond1}).
%%%$rM^{-1}\sum_{k=1}^M k\alpha_X(k)^{1-2/\gamma}= o(1)$.
%
%
%The purpose of the blocking with diverging $M$ is to capture the
%underlying dependence so as to make the GEL self-standardize, a
%virtue of the GEL method that we are interested in despite the high
%dimensionality.  The self-standardizing is the reason behind the GEL
%having the chi-square limit in the fixed dimensional case. Hence,
%from now on we assume $M \to \infty$ at a proper speed.

The following is an extension of Theorem 1 by allowing $V_{M}$ to be
asymptotically singular, namely $\lambda _{\min }(V_{M})\rightarrow 0$ as $%
r\rightarrow \infty $, with $M$ being either fixed or diverging.

\begin{cy}
\label{cy1} Assume conditions \textrm{(A.1)}, \textrm{(A.2)}, and that $%
\lambda_{\min}(V_M) \asymp r^{-\iota_1}$ for some $\iota_1>0$ and $%
\lambda_{\max}(V_M)$ is uniformly bounded away from zero and infinity. Then $%
\|\widehat{\boldsymbol{\theta}}_n-\boldsymbol{\theta}_0\|_2=O_p(r^{(1+%
\iota_1)/2}n^{-1/2})$ and $\|\widehat{\lambda}(\widehat{\boldsymbol{\theta}}%
_n)\|_2=O_p(r^{(1+3\iota_1)/2}Mn^{-1/2})$ provided that $r^{2+3%
\iota_1}M^{2-2/\gamma}n^{2/\gamma-1}=o(1)$, $r^{2+2\iota_1}M^3n^{-1}=o(1)$, $%
r^{2+3\iota_1}pMn^{-1}=o(1)$ and $r^{2+\iota_1}pM^2n^{-1}=o(1)$.
\end{cy}

This corollary shows that when the smallest eigenvalue of $V_{M}$ is not
bounded away from zero, the convergence rates for $\widehat{\boldsymbol{%
\theta }}_{n}$ and the Lagrange multiplier $\widehat{\lambda }(\widehat{%
\boldsymbol{\theta }}_{n})$ become slower. Theorem 1 can be viewed as a
special case of Corollary 1 with $\iota _{1}=0$.

The convergence rate of $\Vert \widehat{\boldsymbol{\theta }}_{n}-%
\boldsymbol{\theta }_{0}\Vert _{2}$ attained in Theorem 1 is dictated by $r$%
, the number of the moment restrictions, rather than by $p$, the dimension
of $\boldsymbol{\theta }$. Under slightly stronger conditions the next
proposition improves the convergence rate to $O_{p}(p^{1/2}n^{-1/2})$.

\begin{pn}
Under conditions \textrm{(A.1)-(A.3)}, assume that the eigenvalues of $V_{M}$
and $V_{n}$ are uniformly bounded away from zero and infinity. Then $\Vert
\widehat{\boldsymbol{\theta }}_{n}-\boldsymbol{\theta }_{0}\Vert
_{2}=O_{p}(p^{1/2}n^{-1/2})$ provided that $r^{3}M^{2-2/\gamma }n^{2/\gamma
-1}=o(1)$, $r^{3}M^{3}n^{-1}=o(1)$, $r^{3}pMn^{-1}=o(1)$ and $%
r^{3}p^{2}n^{-1}=o(1)$.
\end{pn}

%Assuming both $V_M$ and $V_n$'s eigenvalues to be bounded in the above proposition is to contain the cases of  finite and diverging $M$. For $M \to \infty$,  it is only required that  the eigenvalues of $V_n$ being uniformly bounded away from zero and infinity.
% ? {\bf JY: If we let $M\rightarrow\infty$, we can do this. However, the AE fights against letting $M\rightarrow\infty$ directly. I guess this may be better and more general. I have added some comments below this proposition.}}.
%As shown in Appendix, the proof of Proposition 1 is based on the asymptotic expansion given in Proposition 2. Hence, the restrictions among $r, p, M$ and $n$ are applied to control the remainder terms specified in Proposition 2.

%By comparing the restrictions among $r, p, M$ and $n$ employed in
%Theorem 1, we see an extra factor $rp^{-1}$ or $r^{1/2}p^{-1/2}$ used in Proposition 1 for the dependent case
% to improve the convergence rate of
%$\|\widehat{\btheta}_n-\btheta_0\|_2$ from $r^{1/2}n^{-1/2}$ to
%$p^{1/2} n^{-1/2}$. A routine algebra reveals specific rates for $M$
%and $r$: $M=O(n^{\{(\gamma -2)/(4\gamma -2)\}\wedge 1/5})$ and
%$r^3p^{-1}=o(n^{\{(\gamma -2)/(2\gamma -1)\}\wedge 2/5})$ to ensure
%$\|\widehat{\btheta}_{n}-\btheta_{0}\|_2=O_p(p^{1/2}n^{-1/2})$.
%

%\end{spacing}
%\begin{spacing}{1.3}

\section{Asymptotic Normality}

We now turn to the asymptotic normality of the GEL estimator $\widehat{%
\boldsymbol{\theta}}_{n}$. We are in particularly interested in the effect
of the block size $M$ on the estimation efficiency. Based on the consistency
of $\widehat{\boldsymbol{\theta}}_{n}$ and $\widehat{\lambda}(\widehat{%
\boldsymbol{\theta}}_{n})$ given in Theorem 1, expanding $\nabla _{\lambda }%
\widehat{S}_{n}(\widehat{\boldsymbol{\theta}}_{n},\widehat{\lambda}(\widehat{%
\boldsymbol{\theta}}_{n}))=\boldsymbol{\mathbf{0}}$ for $\widehat{\lambda}(%
\widehat{\boldsymbol{\theta}}_{n})$ around $\lambda =\boldsymbol{\mathbf{0}}$
gives
\begin{equation}
\boldsymbol{\mathbf{0}}=\frac{1}{Q}\sum_{q=1}^{Q}\rho_v(0)\phi _{q}(\widehat{%
\boldsymbol{\theta}}_{n})+\frac{1}{Q}\sum_{q=1}^{Q}\rho_{vv}(\widetilde{%
\lambda}^{\prime }\phi_q(\widehat{\boldsymbol{\theta}}_n))\phi _{q}(\widehat{%
\boldsymbol{\theta}}_{n})\phi _{q}(\widehat{\boldsymbol{\theta}}%
_{n})^{\prime }\widehat{\lambda}(\widehat{\boldsymbol{\theta}}_{n}),
\label{eq:0.2}
\end{equation}%
where $\widetilde{\lambda}$ is on the line joining $\boldsymbol{\mathbf{0}}$
and $\widehat{\lambda}(\widehat{\boldsymbol{\theta}}_{n})$. From (\ref%
{eq:0.1}) and (\ref{eq:0.2}), it yields
\begin{equation}
\bigg[\frac{1}{Q}\sum_{q=1}^{Q}\rho_v(\widehat{\lambda}(\widehat{\boldsymbol{%
\theta}}_n)^{\prime }\phi_q(\widehat{\boldsymbol{\theta}}_n))\{\nabla _{%
\boldsymbol{\theta} }\phi _{q}(\widehat{\boldsymbol{\theta}}_{n})\}^{\prime }%
\bigg]\bigg\{\frac{1}{Q}\sum_{q=1}^{Q}\rho_{vv}(\widetilde{\lambda}^{\prime
}\phi_q(\widehat{\boldsymbol{\theta}}_n))\phi _{q}(\widehat{\boldsymbol{%
\theta}}_{n})\phi _{q}(\widehat{\boldsymbol{\theta}}_{n})^{\prime }\bigg\}%
^{-1}\bar{\phi} (\widehat{\boldsymbol{\theta}}_{n})=\boldsymbol{\mathbf{0}}.
\label{eq:key}
\end{equation}%
Based on (\ref{eq:key}), we can establish the following proposition which is
the starting point in our study of the asymptotic normality of $\widehat{%
\boldsymbol{\theta}}_{n}$.

\begin{pn}
Under conditions \textrm{(A.1)-(A.3)}, assume that the eigenvalues of $V_M$
and $V_n$ are uniformly bounded away from zero and infinity. If $%
r^2pM^2n^{-1}=o(1)$ and \textrm{(\ref{eq:cond1})} holds, then for any vector
$\boldsymbol{\alpha}_{n}\in\mathbb{R}^p$ with unit $L_2$-norm,
\begin{equation*}
\begin{split}
&\sqrt{n}\boldsymbol{\alpha}_{n}^{\prime }([{E}\{\nabla_{\boldsymbol{\theta}%
} g_t(\boldsymbol{\theta}_{0})\}]^{\prime }V_M^{-1}V_nV_M^{-1}[{E}\{\nabla_{%
\boldsymbol{\theta}} g_t(\boldsymbol{\theta}_{0})\}])^{-1/2}[E\{\nabla_{%
\boldsymbol{\theta}} g_t(\boldsymbol{\theta}_{0})\}]^{\prime }V_M^{-1}[{E}%
\{\nabla_{\boldsymbol{\theta}} g_t(\boldsymbol{\theta}_{0})\}](\widehat{%
\boldsymbol{\theta}}_{n}-\boldsymbol{\theta}_{0}) \\
=&-\sqrt{n}\boldsymbol{\alpha}_{n}^{\prime }([{E}\{\nabla_{\boldsymbol{\theta%
}} g_t(\boldsymbol{\theta}_{0})\}]^{\prime }V_M^{-1}V_nV_M^{-1}[{E}\{\nabla_{%
\boldsymbol{\theta}} g_t(\boldsymbol{\theta}_{0})\}])^{-1/2}[{E}\{\nabla_{%
\boldsymbol{\theta}} g_t(\boldsymbol{\theta}_{0})\}]^{\prime }V_M^{-1}\bar{g}%
(\boldsymbol{\theta}_{0}) \\
&+O_p(r^{3/2}p^{1/2}M^{1/2}n^{-1/2})+O_p(r^{3/2}pn^{-1/2})+O_p(r^{3/2}M^{1-1/\gamma}n^{1/\gamma-1/2})+O_p(r^{3/2}M^{3/2}n^{-1/2}).
\end{split}%
\end{equation*}
\end{pn}

Proposition 2 covers both the finite and diverging $M$ cases. When $M$ is
diverging, $\Vert V_{M}-V_{n}\Vert _{2}\rightarrow 0$ provided that $r=o(M)$%
. We can replace $V_{M}^{-1}V_{n}V_{M}^{-1}$ on the left-hand side of above
asymptotic expansion by $V_{n}^{-1}$ via adding an extra high order term on
the right-hand side of the asymptotic expansion. Let
\begin{equation}
\boldsymbol{\beta }_{n}=-V_{M}^{-1}[{E}\{\nabla _{\boldsymbol{\theta }}g_{t}(%
\boldsymbol{\theta }_{0})\}]([{E}\{\nabla _{\boldsymbol{\theta }}g_{t}(%
\boldsymbol{\theta }_{0})\}]^{\prime }V_{M}^{-1}V_{n}V_{M}^{-1}[{E}\{\nabla
_{\boldsymbol{\theta }}g_{t}(\boldsymbol{\theta }_{0})\}])^{-1/2}\boldsymbol{%
\alpha }_{n}  \label{eq:betan}
\end{equation}%
and $U_{n,t}=n^{-1/2}\boldsymbol{\beta }_{n}^{\prime }g_{t}(\boldsymbol{%
\theta }_{0})$ for $t=1,\ldots ,n$. From Proposition 2, a major point of
interest is under what conditions $\sum_{t=1}^{n}U_{n,t}$ is asymptotically
normal.

Let us first consider the easier case where the observations $%
\{X_t\}_{t=1}^n $ are independent. From Lindeberg-Feller theorem %
\citep{Durrett_2010}, to attain the asymptotic normality of $%
\sum_{t=1}^nU_{n,t}$, it suffices to verify the following two conditions,%
%the following two conditions:
\begin{equation*}
\text{(i)}~\sum_{t=1}^n{E}(U_{n,t}^2)\rightarrow1~~\text{and}~~\text{(ii)}%
~\sum_{t=1}^n{E}\{U_{n,t}^21_{(|U_{n,t}|>\varepsilon)}\}\rightarrow0 ~~\text{%
as}~~ n\rightarrow\infty~~\text{for any}~~\varepsilon>0.
\end{equation*}
Actually, $V_n=V_M={E}\{ g_t(\boldsymbol{\theta}_{0})g_t(\boldsymbol{\theta}%
_{0})^{\prime }\}$ in this case. Hence, $\sum_{t=1}^n{E}(U_{n,t}^2)=1$. Note
that $\|\boldsymbol{\beta}_n\|_2\leq \lambda_{\text{min}}^{-1/2}(V_n)$ which
is uniformly bounded away from infinity if $\lambda_{\min}(V_n)$ is
uniformly bounded away from zero. Hence, by (A.2)(ii),
\begin{equation*}
(n^{1/2}\varepsilon)^{\gamma-2}{E}[|\boldsymbol{\beta}_n^{\prime }g_t(%
\boldsymbol{\theta}_{0})|^21_{\{|\boldsymbol{\beta}_n^{\prime }g_t(%
\boldsymbol{\theta}_{0})|>n^{1/2}\varepsilon\}}]\leq{E}\{|\boldsymbol{\beta}%
_n^{\prime }g_t(\boldsymbol{\theta}_{0})|^{\gamma}\}\leq Cr^{\gamma/2},
\end{equation*}
which implies that part (ii) holds if $rn^{2/\gamma-1}=o(1)$.
%\footnote{Should it be $\zeta_2^2(r)n^{2/\gamma-1}=o(1)$ ?}
Therefore, %if we require
%\begin{equation}
%r^3p^2n^{-1}=o(1) \quad \hbox{and} \quad
%\zeta_2^2(r)r^2n^{2/\gamma-1}=o(1), \label{eq:cond2}
%\end{equation}
$\sum_{t=1}^nU_{n,t}\xrightarrow{d}N(0,1)$ provided that $%
rn^{2/\gamma-1}=o(1)$ for any selection of $\boldsymbol{\alpha}_{n} \in
\mathbb{R}^p$ with unit $L_2$-norm.

For dependent data, we need to assume $\sup_{n}{E}\{|\boldsymbol{\beta }%
_{n}^{\prime }g_{t}(\boldsymbol{\theta }_{0})|^{2+v}\}<\infty $ for some $%
v>0 $, namely $|\boldsymbol{\beta }_{n}^{\prime }g_{t}(\boldsymbol{\theta }%
_{0})| $ has a higher than two uniformly bounded moment. This is required in
the central limit theorem for dependent processes as carried out in \cite%
{PeligradUtev_1997_AOP} and \cite{FrancqZakoian_2005_ET}.
% The higher than the second moment on $\beta_n^{\prime
%}g_t(\theta_{0})$
It is used to guarantee the limit of $\text{Var}\{n^{-1/2}\sum_{t=1}^{n}%
\boldsymbol{\beta }_{n}^{\prime }g_{t}(\boldsymbol{\theta }_{0})\}$ can be
well defined as $n\rightarrow \infty $. More specifically, notice that
\begin{equation*}
\text{Var}\bigg\{\frac{1}{n^{1/2}}\sum_{t=1}^{n}\boldsymbol{\beta }%
_{n}^{\prime }g_{t}(\boldsymbol{\theta }_{0})\bigg\}=E\{|\boldsymbol{\beta }%
_{n}^{\prime }g_{1}(\boldsymbol{\theta }_{0})|^{2}\}+2\sum_{k=1}^{n-1}\bigg(%
1-\frac{k}{n}\bigg)E\{\boldsymbol{\beta }_{n}^{\prime }g_{1}(\boldsymbol{%
\theta }_{0})g_{1+k}(\boldsymbol{\theta }_{0})^{\prime }\boldsymbol{\beta }%
_{n}\},
\end{equation*}%
to define the limit of above sum of series, we need that $\text{Var}%
\{n^{-1/2}\sum_{t=1}^{n}\boldsymbol{\beta }_{n}^{\prime }g_{t}(\boldsymbol{%
\theta }_{0})\}$ is absolutely convergent, i.e.,
\begin{equation*}
\lim_{n\rightarrow \infty }\bigg[E\{|\boldsymbol{\beta }_{n}^{\prime }g_{1}(%
\boldsymbol{\theta }_{0})|^{2}\}+2\sum_{k=1}^{n-1}\bigg(1-\frac{k}{n}\bigg)%
|E\{\boldsymbol{\beta }_{n}^{\prime }g_{1}(\boldsymbol{\theta }_{0})g_{1+k}(%
\boldsymbol{\theta }_{0})^{\prime }\boldsymbol{\beta }_{n}\}|\bigg]<\infty .
\end{equation*}%
By Davydov inequality \citep{Davydov_1968_TPIA, Rio_1993_AIHP}, the absolute
convergence of $\text{Var}\{n^{-1/2}\sum_{t=1}^{n}\boldsymbol{\beta }%
_{n}^{\prime }g_{t}(\boldsymbol{\theta }_{0})\}$ will hold by requiring $%
\sup_{n}E\{|\boldsymbol{\beta }_{n}^{\prime }g_{t}(\boldsymbol{\theta }%
_{0})|^{2+v}\}<\infty $ for some suitable $v$.

%the limit of $\textrm{Var}\{n^{-1/2}\sum_{t=1}^n\bbetaa_n'g_t(\btheta_0)\}$ can be well defined as $n\rightarrow\infty$. More specifically, notice that
%\[
%\textrm{Var}\bigg\{\frac{1}{n^{1/2}}\sum_{t=1}^n\bbetaa_n'g_t(\btheta_0)\bigg\}=E\{|\bbetaa_n'g_1(\btheta_0)|^2\}+2\sum_{k=1}^{n-1}\bigg(1-\frac{k}{n}\bigg)E\{\bbetaa_n'g_1(\btheta_0)g_{1+k}'(\btheta_0)\bbetaa_n\},
%\]
%to define the limit of above sum of series, we need

For high dimensional moment equation $g(x,\theta)$ with diverging $r$, we
need
\begin{equation}
\sup_n{E}\{|\boldsymbol{\beta}_n^{\prime }g_t(\boldsymbol{\theta}%
_{0})|^{\gamma}\}<\infty  \label{eq:moment-cond}
\end{equation}
for $\boldsymbol{\beta}_n$ defined via (\ref{eq:betan}) and $\gamma > 2$
specified in (A.1)(i).
% is defined in (A.1)(i) and used to describe the speed of
%decay in the $\alpha$-mixing coefficients.
A sufficient condition for (\ref{eq:moment-cond}) is to restrict%
% $\beta_n \in \mathcal{D}(K)$ where %(and hence indirectly on
%$\alpha_n$) to ensure $\beta_n'g_i(\theta_0)$ has higher than the second moments uniformly . %It seems like that lower finite moment of a random
%variable can not guarantee it has higher finite moment.
% Define
\begin{equation*}
\boldsymbol{\beta}_n\in\mathcal{D}(K):=\bigg\{(v_1,v_2,\ldots)\in\mathbb{R}%
^\infty:\sum_{k=1}^{\infty}|v_k|\leq K\bigg\},
\end{equation*}
where $K$ is a given finite constant. To appreciate this, write $\boldsymbol{%
\beta}_n=(\beta_{n,1},\ldots,\beta_{n,r})^{\prime }$ and let $%
\kappa_r=\sum_{j=1}^r|\beta_{n,j}|$. Then,
\begin{equation*}
{E}\{|\boldsymbol{\beta}_n^{\prime }g_t(\boldsymbol{\theta}%
_{0})|^{\gamma}\}=\kappa_r^{\gamma}{E}\bigg\{\bigg|\sum_{j=1}^r\frac{%
|\beta_{n,j}|}{\kappa_r}\text{sign}(\beta_{n,j})g_{t,j}(\boldsymbol{\theta}%
_{0})\bigg|^{\gamma}\bigg\}\leq K^{\gamma} C,
\end{equation*}
where the last step is based on the Jensen's inequality and (A.2)(iii). If $%
\sum_{j=1}^r|\beta_{n,j}|\rightarrow\infty$ as $n\rightarrow\infty$, we can
construct a counter-example such that $\sup_n{E}\{|\boldsymbol{\beta}%
_n^{\prime }g_t(\boldsymbol{\theta}_{0})|^{2+v}\}\rightarrow\infty$ for any $%
v>0$. The following theorem establishes the asymptotic normality of $%
\widehat{\boldsymbol{\theta}}_{n}$.
% for high dimensional estimating equations under dependence.

\begin{tm}
Under conditions \textrm{(A.1)-(A.3)}, assume that the eigenvalues of $V_M$
and $V_n$ are uniformly bounded away from zero and infinity. For dependent
data, if
\begin{equation}
r^3M^{2-2/\gamma}n^{2/\gamma-1}=o(1),~r^3M^3n^{-1}=o(1),~r^{3}pMn^{-1}=o(1)~~%
\text{and}~~r^3p^2n^{-1}=o(1),  \label{eq:condasynor}
\end{equation}
then for any $\boldsymbol{\alpha}_{n}\in\mathbb{R}^{p}$ with unit $L_2$-norm
such that \textrm{(\ref{eq:moment-cond})} holds,
%$\beta_n \in\mathcal{D}(K)$,\footnote{Can we replaced it with the moment condtion ?}
\begin{equation*}
\sqrt{n}\boldsymbol{\alpha}_{n}^{\prime }([{E}\{\nabla_{\boldsymbol{\theta}}
g_t(\boldsymbol{\theta}_{0})\}]^{\prime }V_M^{-1}V_nV_M^{-1}[{E}\{\nabla_{%
\boldsymbol{\theta}} g_t(\boldsymbol{\theta}_{0})\}])^{-1/2}[E\{\nabla_{%
\boldsymbol{\theta}} g_t(\boldsymbol{\theta}_{0})\}]^{\prime }V_M^{-1}[{E}%
\{\nabla_{\boldsymbol{\theta}} g_t(\boldsymbol{\theta}_{0})\}](\widehat{%
\boldsymbol{\theta}}_{n}-\boldsymbol{\theta}_{0})
\end{equation*}
converges to $N(0,1)$ as $n\rightarrow\infty$.
\end{tm}

For finite block size $M$, the above asymptotic distribution holds provided
that
\begin{equation*}
r^3n^{2/\gamma-1}=o(1) \quad \hbox{and} \quad r^3p^2n^{-1}=o(1).
\end{equation*}
Since
\begin{equation*}
\begin{split}
&([E\{\nabla_{\boldsymbol{\theta}} g_t(\boldsymbol{\theta}_{0})\}]^{\prime
}V_M^{-1}[{E}\{\nabla_{\boldsymbol{\theta}} g_t(\boldsymbol{\theta}%
_{0})\}])^{-1}([{E}\{\nabla_{\boldsymbol{\theta}} g_t(\boldsymbol{\theta}%
_{0})\}]^{\prime }V_M^{-1}V_nV_M^{-1}[{E}\{\nabla_{\boldsymbol{\theta}} g_t(%
\boldsymbol{\theta}_{0})\}]) \\
&~~~~~~~~~~~~~~\times([E\{\nabla_{\boldsymbol{\theta}} g_t(\boldsymbol{\theta%
}_{0})\}]^{\prime }V_M^{-1}[{E}\{\nabla_{\boldsymbol{\theta}} g_t(%
\boldsymbol{\theta}_{0})\}])^{-1}\geq([E\{\nabla_{\boldsymbol{\theta}} g_t(%
\boldsymbol{\theta}_{0})\}]^{\prime }V_n^{-1}[{E}\{\nabla_{\boldsymbol{\theta%
}} g_t(\boldsymbol{\theta}_{0})\}])^{-1},
\end{split}%
\end{equation*}
the GEL estimator is asymptotically efficient if $\|V_M-V_n\|_2\rightarrow 0$%
, which implies $V_M^{-1}V_nV_M^{-1}$ is asymptotically equivalent to $%
V_n^{-1}$. This means that if $\{g_t(\boldsymbol{\theta}_0)\}_{t=1}^n$ is a
martingale difference sequence, as $V_M=V_n=E\{g_1(\boldsymbol{\theta}_0)g_1(%
\boldsymbol{\theta}_0)^{\prime }\}$ for any $M\geq1$, selecting $M=1$ will
lead to the efficient GEL estimation. In a general case where the nature of
the dependence in the estimating function is unknown, letting $%
M\rightarrow\infty$ at some suitable diverging rate, so that (\ref%
{eq:condasynor}) is satisfied, will lead to the efficient estimation.
Specifically, as
\begin{equation*}
\|V_M-V_n\|_2\leq CrM^{-1}\sum_{k=1}^Mk\alpha_X(k)^{1-2/\gamma},
\end{equation*}
under (A.1)(i) and (A.2)(ii), choosing $M \to \infty$ such that $r=o(M)$
produces the asymptotically efficient GEL estimator $\widehat{\boldsymbol{%
\theta}}_n$. According to (\ref{eq:condasynor}) and $r=o(M)$, the divergence
rate in $M$ is $M=O(n^{\{(\gamma-2)/(5\gamma-2)\}\wedge{1/6}})$ while $%
r=o(n^{\{(\gamma-2)/(5\gamma-2)\}\wedge {1/6}})$, regardless $p $ being
fixed or diverging. Under such setting, the best growth rate for $r$
%the number of estimating equations
is $r=o(n^{1/6}) $ when $\gamma\geq10$. In comparison with the case of
finite $M$, while letting $M\rightarrow\infty$ can guarantee the efficiency,
it does slows the divergence of $r$.

If the smallest eigenvalues of $V_M$ and $V_n$ decay to zero as $%
r\rightarrow\infty$,
%\footnote{Can we remove this section ? I do not see the AE has asked for this, haven't he ? {\bf He asked this question. See (iii) for details.}}
we assume $\lambda_{\min}(V_M)\asymp r^{-\iota_1}$ and $\lambda_{\min}(V_n)%
\asymp r^{-\iota_2}$ for some positive $\iota_1$ and $\iota_2$. Based on
Corollary \ref{cy1}, by repeating the proof of Proposition 2 in Appendix, it
can be shown that the leading term in the asymptotic expansion of
Proposition 2 remains while the four remainder terms become
\begin{equation*}
\begin{split}
&~O_p(r^{(3+6\iota_1+\iota_2)/2}p^{1/2}M^{1/2}n^{-1/2})+O_p(r^{(3+4\iota_1+%
\iota_2)/2}pn^{-1/2}) \\
+&~O_p(r^{(3+5\iota_1+\iota_2)/2}M^{1-1/\gamma}n^{1/\gamma-1/2})+O_p(r^{(3+5%
\iota_1+\iota_2)/2}M^{3/2}n^{-1/2})
\end{split}%
\end{equation*}
provided that the conditions governing $r, p, M$ and $n$ assumed in
Corollary \ref{cy1} hold.
%, i.e., $r^{2+3\iota_1}M^{2-2/\gamma}n^{2/\gamma-1}=o(1)$, $r^{2+2\iota_1}M^3n^{-1}=o(1)$, $r^{2+3\iota_1}pMn^{-1}=o(1)$ and $r^{2+\iota_1}pM^2n^{-1}=o(1)$.
By the central limit theorem established in \cite{FrancqZakoian_2005_ET},
the leading term in the asymptotic expansion of Proposition 2 converges to $%
N(0,1)$ regardless $\lambda_{\min}(V_M)\rightarrow0$ and $%
\lambda_{\min}(V_n)\rightarrow0$ or not. Hence, the asymptotic normality of
the GEL estimator $\widehat{\boldsymbol{\theta}}_n$ is valid free of the
statue of the eigenvalues of $V_M$ and $V_n$.
%same as of which $\lambda_{\min}(V_m)$ and $\lambda_{\min}(V_n)$ are uniformly bounded.even if $\lambda_{\min}(V_M)\rightarrow0$ and $\lambda_{\min}(V_n)\rightarrow0$,
The difference is that when $\lambda_{\min}(V_M)\rightarrow0$ and $%
\lambda_{\min}(V_n)\rightarrow0$, the growth rate of $r$ and/or $M$ are
reduced.

To put the growth rate of $r$ into perspectives and to highlight the impacts
of data dependence, we consider the independent analogue of Theorem 2 in the
following, whose proof is obtained by assigning $\alpha_X(k)=0$ and $M=1$ in
Proposition 2.

\begin{cy}
Under conditions \textrm{(A.1)-(A.3)}, assume that the eigenvalues of $%
E\{g_1(\boldsymbol{\theta}_0)g_1(\boldsymbol{\theta}_0)^{\prime }\}$ are
uniformly bounded away from zero and infinity. For independent data, if $%
r^3p^2n^{-1}=o(1)$ and $r^3n^{2/\gamma-1}=o(1)$, then for any $\boldsymbol{%
\alpha}_{n}\in\mathbb{R}^{p}$ with unit $L_2$-norm,
\begin{equation*}
\sqrt{n}\boldsymbol{\alpha}_{n}^{\prime }([{E}\{\nabla_{\boldsymbol{\theta}}
g_t(\boldsymbol{\theta}_{0})\}]^{\prime }V_n^{-1}[{E}\{\nabla_{\boldsymbol{%
\theta}} g_t(\boldsymbol{\theta}_{0})\}])^{1/2}(\widehat{\boldsymbol{\theta}}%
_{n}-\boldsymbol{\theta}_{0})\xrightarrow[]{d}N(0,1)~~\text{as}%
~n\rightarrow\infty.
\end{equation*}
\end{cy}

The above corollary shows that under independence, the growth rate for $r$
is $o(n^{1/3-2/(3\gamma)})$ if $p$ is fixed. If $\gamma$ is sufficiently
large, the rate of $r$ can be close to $o(n^{1/3})$. If $p$ grows with $r$
and $p/r\rightarrow y\in(0,1]$, then $r=o(n^{\{1/3-2/(3\gamma)\}\wedge1/5})$%
. In particular, if $\gamma\geq5$, $r=o(n^{1/5})$ which retains Theorem 2 in
\cite{LengTang_2011_Biometrika} for the EL estimator. Comparing the growth
rates for $r$ under the dependent and independent settings, when $M$ is
diverging, we see a slowing down in the rate under dependence from $%
o(n^{1/5})$ to $o(n^{1/6})$ if the best moment conditions hold.

If $p$, the dimension of $\boldsymbol{\theta }$, is fixed, as in a case of
conditional moment restrictions in Example 3, the asymptotic normality of $%
\widehat{\boldsymbol{\theta }}_{n}$ can be attained with some ease. It can
be shown that $\boldsymbol{\beta }_{n}$ is automatically in $\mathcal{D}(K)$
for a large enough $K$, which implies the condition (\ref{eq:moment-cond})
holds for any $\boldsymbol{\alpha }_{n}\in \mathbb{R}^{p}$ with unit $L_{2}$%
-norm. This is summarized in the following corollary.

\begin{cy}
Under conditions \textrm{(A.1)-(A.3)}, assume that the eigenvalues of $V_M$
and $V_n$ are uniformly bounded away from zero and infinity. For dependent
data, if $p$ is fixed, then for any $\boldsymbol{\alpha}_{n}\in\mathbb{R}^p$
with unit $L_2$-norm,
\begin{equation*}
\sqrt{n}\boldsymbol{\alpha}_{n}^{\prime }([{E}\{\nabla_{\boldsymbol{\theta}}
g_t(\boldsymbol{\theta}_{0})\}]^{\prime }V_M^{-1}V_nV_M^{-1}[{E}\{\nabla_{%
\boldsymbol{\theta}} g_t(\boldsymbol{\theta}_{0})\}])^{-1/2}[E\{\nabla_{%
\boldsymbol{\theta}} g_t(\boldsymbol{\theta}_{0})\}]^{\prime }V_M^{-1}[{E}%
\{\nabla_{\boldsymbol{\theta}} g_t(\boldsymbol{\theta}_{0})\}](\widehat{%
\boldsymbol{\theta}}_{n}-\boldsymbol{\theta}_{0})
\end{equation*}
converges to $N(0,1)$ as $n\rightarrow\infty$, provided that $%
r^3M^{2-2/\gamma}n^{2/\gamma-1}=o(1)$ and $r^3M^3n^{-1}=o(1)$.
\end{cy}

This Corollary with $M=1$ recovers that in \cite{DonaldImbensNewey_2003_JOE}
for iid data.

\section{Generalized Empirical Likelihood Ratios}

The EL ratio $w_{n}(\boldsymbol{\theta })=-2\log \{Q^{Q}\mathcal{L}(%
\boldsymbol{\theta })\}$ for $\mathcal{L}(\boldsymbol{\theta })$ defined in (%
\ref{eq:L}) plays an important role in the statistical inference. A
prominent result for fixed dimensional EL is its resembling the parameter
likelihood by have a limiting chi-square distribution under a wide range of
situations, as demonstrated in \cite{Owen_1988_Biometrika}, \cite%
{ChenCui_2003}, \cite{QinLawless_1994_AOS} and \cite%
{ChenVanKeilegom_2009_Test} for independent data, and \cite%
{Kitamura_1997_AOS} for dependent data.

For GEL, we define the GEL ratio as
\begin{equation}
w_n(\boldsymbol{\theta})=\frac{2\rho_{vv}(0)}{\rho_v^2(0)}\bigg\{%
Q\rho(0)-\max_{\lambda\in\widehat{\Lambda}_n(\boldsymbol{\theta}%
)}\sum_{q=1}^Q\rho(\lambda^{\prime }\phi_q(\boldsymbol{\theta}))\bigg\}
\label{eq:gelratio}
\end{equation}
which is the extension of the EL ratio in the GEL framework.

We consider the asymptotic distribution of the GEL ratio $w_{n}(\boldsymbol{%
\theta }_{0})$ when both $r$ and $p$ are diverging. Under such setting, a
natural form of the Wilks' theorem is
\begin{equation}
(2r)^{-1/2}\{w_{n}(\boldsymbol{\theta }_{0})-r\}\xrightarrow {d}N(0,1)~~%
\text{as}~r\rightarrow \infty .  \label{eq:anchisq}
\end{equation}%
For the case of means where $g_{t}(\boldsymbol{\theta })=X_{t}-\boldsymbol{%
\theta }$ with independent observations, \cite{ChenPengQin_2009_Biometrika}
and \cite{HjortMcKeagueVanKeilegom_2009_AOS} evaluated the impact of the
dimensionality on the asymptotic distribution (\ref{eq:anchisq}) for the EL
ratio by providing various diverging rates for $r$. For parameters defined
by general moment restrictions, establishing the limiting distribution of
the GEL ratio is far more challenging. We need the following stronger
version of (A.1)(i):

\medskip

(A.1)'(i) There is some $\eta>8$ such that $\alpha_X(k)^{1-2/\gamma}\asymp
k^{-\eta}$ where $\gamma$ is given in (A.2).

\medskip

Condition (A.1)'(i) is used to guarantee the leading order term of (\ref%
{eq:gelratio}) has the similar probabilistic behavior as the chi-square
distribution. It is automatically satisfied with $\eta =\infty $ if $X_{t}$
is exponentially strong mixing or independent. We also need the following
conditions:
\begin{equation}
r^{3}M^{2-2/\gamma }n^{2/\gamma -1}=o(1),~~r^{3}M^{3}n^{-1}=o(1)~~\text{and}%
~~r^{3/2}M^{-1}\sum_{k=1}^{M}k\alpha _{X}(k)^{1-2/\gamma }=o(1).
\label{eq:cond-ratio}
\end{equation}%
Define
\begin{equation}
\xi =\frac{\eta -8}{4\eta +4}1_{\{8<\eta <32\}}+\frac{2}{11}1_{\{32\leq \eta
\leq \infty \}}+1_{\{\text{indenpendent data}\}}.  \label{eq:xi}
\end{equation}%
The next theorem establishes the asymptotic distribution of $w_{n}(%
\boldsymbol{\theta }_{0})$.

\begin{tm}
Under conditions \textrm{(A.1)'(i)}, \textrm{(A.1)(ii)} and \textrm{%
(A.2)(iii)}, assume that the eigenvalues of $V_n$ are uniformly bounded away
from zero and infinity. If \textrm{(\ref{eq:cond-ratio})} holds and $%
r=o(n^\xi)$ where $\xi$ is defined in \textrm{(\ref{eq:xi})},
% where $\xi$ is defined in (\ref{eq:xi}),
%where $\xi$ is defined in (\ref{eq:xi}),
then
\begin{equation*}
(2r)^{-1/2}\{w_n(\boldsymbol{\theta}_{0})-r\}\xrightarrow[]{d}N(0,1)~~ %
\mbox{as $r\rightarrow\infty$.}
\end{equation*}
\end{tm}

This theorem is new for dependent data and includes some established results
for independent data as special cases. For independent data, this theorem
implies that the asymptotic normality of the GEL ratio is valid if $%
r=o(n^{1/3-2/(3\gamma )})$, which is the same as that in \cite%
{HjortMcKeagueVanKeilegom_2009_AOS} for the EL ratio with independent data.
Our result is more general than theirs since we allow for GEL ratio and for
dependent data. For dependent data, the block size is $M=O(n^{(\gamma
-2)/(4\gamma -2)})$ if $2<\gamma <8$ and $M=O(n^{1/5})$ otherwise, and hence
the asymptotic distribution (\ref{eq:anchisq}) holds if $r=o(n^{\delta })$
with
\begin{equation*}
\delta =\min \left( \frac{\eta -8}{4\eta +4}1_{\{8<\eta <32\}}+\frac{2}{11}%
1_{\{32\leq \eta \leq \infty \}},\frac{\gamma -2}{6\gamma -3}1_{\{2<\gamma
<8\}}+\frac{2}{15}1_{\{\gamma \geq 8\}}\right) .
\end{equation*}%
When $\eta $ and $\gamma $ are sufficiently large, the best diverging rate
is $r=o(n^{2/15})$ for the dependent case, which is slower than the rate of $%
r=o(n^{1/3-2/(3\gamma )})$ for the independence case.

\section{Test for Over-identification}

{ For moment restrictions, it is important to check on the validity of the model by testing the following hypotheses
\begin{equation*}
H_{0}:{E}\{g(X_{t},\boldsymbol{\theta }_{0})\}=\boldsymbol{\mathbf{0}}~~%
\text{for some}~\boldsymbol{\theta }_{0}\in \Theta ~~\text{v.s.}~~H_{1}:{E}%
\{g(X_{t},\boldsymbol{\theta })\}\neq \boldsymbol{\mathbf{0}}~~\text{for any}%
~\boldsymbol{\theta }\in \Theta .
\end{equation*}%
We consider testing the above hypothesis when $r>p$, namely the moment
equation overly identify the parameter $\boldsymbol{\theta }$.}

We formulate the test statistic as the GEL ratio $w_{n}(\widehat{\boldsymbol{%
\theta }}_{n})$. For the EL ratio, it has been demonstrated in the fixed
dimensional case by \cite{QinLawless_1994_AOS} and \cite{Kitamura_1997_AOS}
that
\begin{equation*}
w_{n}(\widehat{\boldsymbol{\theta }}_{EL})\xrightarrow{d}\chi _{r-p}^{2}~~~%
\text{as}~n\rightarrow \infty
\end{equation*}%
under $H_{0}$. This mirrors the J-test of \cite{Hansen_1982_Ecma}'s GMM with
fixed and finite dimensions $r$ and $p$.

To formulate the GEL specification test allowing for increasing dimensions $%
r $ and $p$, we are to study the asymptotic distribution of $w_{n}(\widehat{%
\boldsymbol{\theta }}_{n})$ under $H_{0}$ first. We only need to consider
its leading order $n\bar{g}(\widehat{\boldsymbol{\theta }}_{n})^{\prime }\{M%
\widehat{\Omega }(\widehat{\boldsymbol{\theta }}_{n})\}^{-1}\bar{g}(\widehat{%
\boldsymbol{\theta }}_{n})$ as the remainder terms in the asymptotic
expansion of $w_{n}(\widehat{\boldsymbol{\theta }}_{n})$ can be shown to be
of a smaller order. Since $\widehat{\boldsymbol{\theta }}_{n}$ is consistent
for $\boldsymbol{\theta }_{0}$ under $H_{0}$, Lemma \ref{la19} in Appendix
establishes the relationship between the asymptotic distributions of $n\bar{g%
}(\widehat{\boldsymbol{\theta }}_{n})^{\prime }\{M\widehat{\Omega }(\widehat{%
\boldsymbol{\theta }}_{n})\}^{-1}\bar{g}(\widehat{\boldsymbol{\theta }}_{n})$
and $n\bar{g}(\boldsymbol{\theta }_{0})^{\prime }V_{n}^{-1}\bar{g}(%
\boldsymbol{\theta }_{0})$ under $H_{0}$. We need the following conditions:
\begin{equation}\label{eq:8}
\begin{split}
& ~~~~~~~~r^{3}pn^{-1}=o(1),~~pr^{-1/2}=o(1),~~r^{3}M^{3}n^{-1}=o(1), \\
& r^{3}M^{2-2/\gamma }n^{2/\gamma -1}=o(1)~~\text{and}~~r^{3/2}M^{-1}%
\sum_{k=1}^{M}k\alpha _{X}(k)^{1-2/\gamma }=o(1).
\end{split}%
\end{equation}%
Compared with the conditions for the asymptotic distribution of $w_{n}(%
\boldsymbol{\theta }_{0})$ in (\ref{eq:cond-ratio}), the first two
restrictions in (\ref{eq:8}) are the extra ones used to control the
remainder terms.

\begin{tm}
Under conditions \textrm{(A.1)'(i)}, \textrm{(A.1)(ii)}, \textrm{(A.1)(iii)}%
, \textrm{(A.1)(iv)}, \textrm{(A.2)} and \textrm{(A.3)}, assume that the
eigenvalues of $V_n$ are bounded away from zero and infinity. If \textrm{(%
\ref{eq:8})} holds and $r=o(n^\xi)$, where $\xi$ is defined in \textrm{(\ref%
{eq:xi})}, then
\begin{equation*}
\{2(r-p)\}^{-1/2}\{w_n(\widehat{\boldsymbol{\theta}} _{n})-(r-p)\}%
\xrightarrow{d}N(0,1) \quad \mbox{as $r\rightarrow\infty$.}
\end{equation*}
\end{tm}

The asymptotic normality can be used to derive the over-identification test
under high dimensionality and dependence. Specifically, $H_0$ is rejected if
\begin{equation*}
\{2(r-p)\}^{-1/2}\{w_n(\widehat{\boldsymbol{\theta}}_{n})-(r-p)\} >
z_{1-\alpha}
\end{equation*}
where $z_{1-\alpha}$ is the $1-\alpha$ quantile of $N(0,1)$.

{To show the above GEL test for over-identification is consistent, we assume
that under the alternative hypothesis $H_{1}$,
\begin{equation}
\inf_{\boldsymbol{\theta }\in \Theta }\Vert {E}\{g(X_{t},\boldsymbol{\theta }%
)\}\Vert _{2}\geq \varsigma .  \label{eq:cond61}
\end{equation}%
The following theorem describes the behavior of $\{2(r-p)\}^{-1/2}\{w_{n}(%
\widehat{\boldsymbol{\theta }}_{n})-(r-p)\}$ under $H_{1}$.}

\begin{tm}
Under conditions \textrm{(A.1)(i)}, \textrm{(A.1)(ii)}, \textrm{(A.1)(iv)},
\textrm{(A.2)(ii)} and \textrm{(\ref{eq:cond61})}, if there is a positive
constant $\epsilon$ such that $r^2M^{1-2/\gamma}n^{2/\gamma-1}(\log
n)^\epsilon\varsigma^{-2}=O(1)$, $r^{1/2}Mn^{-1}\varsigma^{-1}=o(1)$ and $%
\Delta_1(r,p)\varsigma^{-1}=O(1)$, then
\begin{equation*}
\{2(r-p)\}^{-1/2}\{w_n(\widehat{\boldsymbol{\theta}}_n)-(r-p)\}%
\xrightarrow{p}\infty \quad \mbox{as $r\rightarrow\infty$}.
\end{equation*}
\end{tm}

This theorem shows that the GEL test is consistent. Unlike Theorem 4, this
theorem does not require the block size $M\rightarrow \infty $ and assumes
the weaker condition (A.1)(i) (instead of the more restrictive one
(A.1)'(i)). From the proof given in the Appendix which follows the technique developed in \cite{ChangTangWu_2013}, the test statistic $%
\{2(r-p)\}^{-1/2}\{w_{n}(\widehat{\boldsymbol{\theta }}_{n})-(r-p)\}$ diverges to infinity at least at the rate of $O(r^{1/2})$ under $H_{1}$.

\section{Penalized Generalized Empirical Likelihood}

In high dimensional data analysis, when the dimension of parameters is
large, i.e., $p\rightarrow\infty$, a more reasonable assumption is that only
a subset of the parameters are nonzero. Write $\boldsymbol{\theta}%
_0=(\theta_{01},\ldots,\theta_{0p})^{\prime }\in\mathbb{R}^p$ and define $%
\mathcal{A}=\{j:\theta_{0j}\neq0\}$ with its cardinality $s=|\mathcal{A}|$.
Without loss of generality, let $\boldsymbol{\theta}=(\boldsymbol{\theta}%
^{(1)^{\prime }},\boldsymbol{\theta}^{(2)^{\prime }})^{\prime }$, where $%
\boldsymbol{\theta}^{(1)}\in\mathbb{R}^s$ and $\boldsymbol{\theta}^{(2)}\in%
\mathbb{R}^{p-s}$ correspond to the nonzero and zero components
respectively, i.e., $\boldsymbol{\theta}_0=(\boldsymbol{\theta}%
_{0}^{(1)^{\prime }},\boldsymbol{\mathbf{0}}^{\prime })^{\prime }$. Under
such sparsity, we can allow the number of parameters is larger than the
number of estimating equations, i.e., $p>r$. However, we still need to
assume $s\leq r$, which means that the ``real" parameters can be uniquely
identified by the moment restrictions (\ref{eq:1}). To carry out the
statistical inference on $\boldsymbol{\theta}$ under the sparsity
assumption, we add a penalty term in (\ref{eq:gel}) and the penalized GEL
estimator is defined as
\begin{equation*}
\widehat{\boldsymbol{\theta}}_n^{(\mathrm{pe})}=\arg\min_{\boldsymbol{\theta}%
\in\Theta}\max_{\lambda\in\widehat{\Lambda}_n(\boldsymbol{\theta})}\bigg\{%
\sum_{q=1}^Q\rho(\lambda^{\prime }\phi_M(B_q,\boldsymbol{\theta}%
))+Q\sum_{j=1}^p p_\tau(|\theta_j|)\bigg\}
\end{equation*}
where $p_{\tau}(\cdot)$ is some penalty function with a tuning parameter $%
\tau$. The following conditions are imposed on the penalty function $%
p_{\tau}(\cdot)$ and the tuning parameter $\tau$. \bigskip

(A.4) $\liminf_{\tau\rightarrow0}\liminf_{\theta\rightarrow0+}p^{\prime
}_{\tau}(\theta)/\tau>0. $

\bigskip

(A.5) There exists a positive constant $C$ such that $\max_{j\in\mathcal{A}%
}p_\tau(|\theta_{0j}|)\leq C\tau$.

\bigskip

Conditions (A.4) and (A.5) hold for many penalty functions such as the one
in \cite{FanLi_2001} and the minimax concave penalty of \cite{Zhang_2010}.
Define
\begin{equation*}
\begin{split}
{\mathbf{S}}(\boldsymbol{\theta}_0)=&~\big(\lbrack E\{\nabla_{\boldsymbol{%
\theta}}g_t(\boldsymbol{\theta}_0)\}]^{\prime }V_M^{-1}[E\{\nabla_{%
\boldsymbol{\theta}}g_t(\boldsymbol{\theta}_0)\}]\big)^{-1}\big(\lbrack
E\{\nabla_{\boldsymbol{\theta}}g_t(\boldsymbol{\theta}_0)\}]^{\prime
}V_M^{-1}V_nV_M^{-1}[E\{\nabla_{\boldsymbol{\theta}}g_t(\boldsymbol{\theta}%
_0)\}]\big) \\
&~~~~~~~~~~~~~~~~~\times\big(\lbrack E\{\nabla_{\boldsymbol{\theta}}g_t(%
\boldsymbol{\theta}_0)\}]^{\prime }V_M^{-1}[E\{\nabla_{\boldsymbol{\theta}%
}g_t(\boldsymbol{\theta}_0)\}]\big)^{-1}.
\end{split}%
\end{equation*}
We correspondingly decompose ${\mathbf{S}}(\boldsymbol{\theta}_0)$ as
\begin{equation}  \label{eq:block}
{\mathbf{S}}(\boldsymbol{\theta}_0)=\left(
\begin{array}{cc}
{\mathbf{S}}_{11}(\boldsymbol{\theta}_0) & {\mathbf{S}}_{12}(\boldsymbol{%
\theta}_0) \\
{\mathbf{S}}_{21}(\boldsymbol{\theta}_0) & {\mathbf{S}}_{22}(\boldsymbol{%
\theta}_0) \\
\end{array}
\right)
\end{equation}
where ${\mathbf{S}}_{11}(\boldsymbol{\theta}_0)$ and ${\mathbf{S}}_{22}(%
\boldsymbol{\theta}_0)$ are $s\times s$ and $(p-s)\times(p-s)$ matrices,
respectively. The following restrictions are needed
\begin{equation}
s\tau r^{-1}nM^{-1} =O(1)~~\text{and}~~
\tau(r^{-1}n)^{1/2}M^{-1}\rightarrow\infty.  \label{eq:cond2}
\end{equation}
Write the penalized GEL estimator $\widehat{\boldsymbol{\theta}}_n^{(\mathrm{%
pe})}=(\widehat{\boldsymbol{\theta}}_n^{(1)^{\prime }},\widehat{\boldsymbol{%
\theta}}_n^{(2)^{\prime }})^{\prime }$ and define
\begin{equation*}
{\mathbf{S}}_p(\boldsymbol{\theta}_0)={\mathbf{S}}_{11}(\boldsymbol{\theta}%
_0)-{\mathbf{S}}_{12} (\boldsymbol{\theta}_0){\mathbf{S}}_{22}^{-1}(%
\boldsymbol{\theta}_0) {\mathbf{S}}_{21}(\boldsymbol{\theta}_0).
\end{equation*}
The following theorem describes the basic properties of the penalized GEL
estimator.

\begin{tm}
\label{tm:pen} Under conditions \textrm{(A.1)-(A.5)}, assume that the
eigenvalues of $V_M$ are uniformly bounded away from zero and infinity. If $%
\max_{j\in\mathcal{A}}p_{\tau}^{\prime }(|\theta_{0j}|)=o(r^{-1/2}n^{-1/2})$%
, $\min_{j\in\mathcal{A}}|\theta_{0j}|/\tau\rightarrow\infty$ and \textrm{(%
\ref{eq:cond2})} holds, the following results hold.

\begin{enumerate}

\item $P\{\widehat{\boldsymbol{\theta}}_n^{(2)}=\mathbf{0}\}\rightarrow1$ as
$n\rightarrow\infty$, provided that \textrm{(\ref{eq:cond1})} holds and $%
r^2pM^2n^{-1}=o(1)$;

\item In addition, if the eigenvalues of $V_n$ are uniformly bounded away
from zero and infinity, then for any $\boldsymbol{\alpha}_n\in\mathbb{R}^s$
with unit $L_2$-norm, then
\begin{equation*}
\sqrt{n}\boldsymbol{\alpha}_n^{\prime }{\mathbf{S}}_p^{-1/2}(\boldsymbol{%
\theta}_0)(\widehat{\boldsymbol{\theta}}_n^{(1)}-\boldsymbol{\theta}%
_{0}^{(1)})\xrightarrow {d} N(0,1)~~\text{as}~n\rightarrow\infty,
\end{equation*}
provided that

\begin{enumerate}
\item[\textrm{(a)}] for independent data, $r^3p^2n^{-1}=o(1)$ and $%
r^3n^{2/\gamma-1}=o(1)$;

\item[\textrm{(b)}] for dependent data, \textrm{(\ref{eq:condasynor})} holds
and $\boldsymbol{\alpha}_n$ satisfies \textrm{(\ref{eq:moment-cond})}
%\footnote{The same here too ? {\bf
%JY: Here should be (4.6) not (4.5). I made a mistake here.}}
with
\begin{equation*}
\begin{split}
\boldsymbol{\beta}_n=&~-V_M^{-1}[E\{\nabla_{\boldsymbol{\theta}}g(%
\boldsymbol{\theta}_0)\}]\big(\lbrack E\{\nabla_{\boldsymbol{\theta}}g_t(%
\boldsymbol{\theta}_0)\}]^{\prime }V_M^{-1}V_nV_M^{-1}[E\{\nabla_{%
\boldsymbol{\theta}}g_t(\boldsymbol{\theta}_0)\}]\big)^{-1} \\
&~~~\times[E\{\nabla_{\boldsymbol{\theta}}g_t(\boldsymbol{\theta}_0)\}]%
^{\prime }V_M^{-1}[E\{\nabla_{\boldsymbol{\theta}^{(1)}}g_t(\boldsymbol{%
\theta}_0)\}]\{{\mathbf{S}}_{11}(\boldsymbol{\theta}_0)-{\mathbf{S}}_{12}(%
\boldsymbol{\theta}_0){\mathbf{S}}_{22}^{-1}(\boldsymbol{\theta}_0){\mathbf{S%
}}_{21}(\boldsymbol{\theta}_0)\}^{1/2}\boldsymbol{\alpha}_n.
\end{split}%
\end{equation*}
\end{enumerate}
\end{enumerate}
\end{tm}

Similar to the consistency of GEL estimator, if the eigenvalues of $E\{g_{t}(%
\boldsymbol{\theta }_{0})g_{t}(\boldsymbol{\theta }_{0})^{\prime }\}$ are
uniformly bounded away from zero and infinity, result (i) still holds
without blocking technique if $r^{2}n^{2/\gamma -1}=o(1)$ and $%
r^{2}pn^{-1}=o(1)$ are satisfied. Comparing Theorem 6 with Theorem 2 and
Corollary 2, since ${\mathbf{S}}_{p}(\boldsymbol{\theta }_{0})\leq {\mathbf{S%
}}_{11}(\boldsymbol{\theta }_{0})$, the penalized GEL estimator is more
efficient in estimating the nonzero components. \cite%
{LengTang_2011_Biometrika} considered the theoretical results of the
penalized EL estimator for independent data by assuming $p/r\rightarrow c\in
(0,1)$. Our results extend theirs to penalized GEL estimator for weakly
dependent data without requiring $p/r\rightarrow c\in (0,1)$.

\section{Simulation Results}

In this section, we present simulation results to compare the finite sample
performance of the GEL estimators with the GMM estimator in the high
dimensional time series setting. Three versions of the GEL estimators were
considered in the simulations: the EL, the ET and the CU estimators. We
experimented two forms of the moment restrictions: one was linear, and the
other was nonlinear. The penalized GEL estimator was also considered in the
non-linear case.

We first conducted simulation for the linear moment restrictions with $%
g(X_{t},\boldsymbol{\theta })=X_{t}-\boldsymbol{\theta }$. The observations $%
\{X_{t}\}_{t=1}^{n}$ were generated according to the vector autoregressive
(VAR) model of order $1$: $X_{t}=\psi X_{t-1}+\varepsilon _{t}$ where $%
\varepsilon _{t}\sim N(\boldsymbol{\mathbf{0}},\Sigma _{\varepsilon })$, $%
\Sigma _{\varepsilon }=(\sigma _{i,j})_{p\times p}$, $\sigma _{i,i}=1-\psi
^{2}$, $\sigma _{i,i\pm 1}=0.5(1-\psi ^{2})$ and $\sigma _{i,j}=0$ for $%
|i-j|>1$. The stationary distribution of $X_{t}$ is $N(\boldsymbol{\mathbf{0}%
},\Sigma _{x})$ where $\Sigma _{x}=(\tilde{\sigma}_{i,j})_{p\times p}$ and $%
\tilde{\sigma}_{i,i}=1$, $\tilde{\sigma}_{i,i\pm 1}=0.5$ and $\tilde{\sigma}%
_{i,j}=0$ for $|i-j|>1$. In this model, $p=r$ and the true parameter $%
\boldsymbol{\theta }_{0}=\mathbf{0}\in \mathbb{R}^{p}$.

The second simulation model was the generalized linear model. The covariates
$\{Z_{t}\}_{t=1}^{n}$ were generated with the same VAR(1) process as the $%
\{X_{t}\}_{t=1}^{n}$ in the first model setting. The response variables $%
\{Y_{t}\}_{t=1}^{n}$ were generated from the Bernoulli distribution such
that $P(Y_{t}=1|Z_{t})=\exp (1+Z_{t}^{\prime }\boldsymbol{\theta }%
_{0})/\{1+\exp (1+Z_{t}^{\prime }\boldsymbol{\theta }_{0})\}$ with the true
parameter $\boldsymbol{\theta }_{0}=(0.8,0.2,0,\ldots ,0)^{\prime }\in
\mathbb{R}^{p}$. Then%
\begin{equation*}
E\bigg\{Y_{t}-\frac{\exp (1+Z_{t}^{\prime }\boldsymbol{\theta }_{0})}{1+\exp
(1+Z_{t}^{\prime }\boldsymbol{\theta }_{0})}\bigg|Z_{t}\bigg\}=0.
\end{equation*}%
In this setting, we have nonlinear moment restrictions
\begin{equation*}
g(X_{t},\boldsymbol{\theta })=\left(
\begin{array}{c}
Z_{t} \\
W_{t} \\
\end{array}%
\right) \bigg\{Y_{t}-\frac{\exp (1+Z_{t}^{\prime }\boldsymbol{\theta })}{%
1+\exp (1+Z_{t}^{\prime }\boldsymbol{\theta })}\bigg\},
\end{equation*}%
where $W_{t}=(Z_{1,t}^{2},\ldots ,Z_{p,t}^{2})^{\prime }$ for $%
Z_{t}=(Z_{1,t},\ldots ,Z_{p,t})^{\prime }$. This model is over-identified.
We considered both non-penalized and penalized estimators under this model
setting.

In both simulation models, we chose $n=500$, $1000$ and $2000$,
respectively. The parameter $\psi $ in the VAR(1) process capturing the
serial dependence was set to be $0.1,0.3$ and $0.5$, respectively. The
dimension $p$ was pegged to the sample size $n$ such that $p=\lfloor
cn^{2/15}\rfloor $, where $c=10$ and $12$ in the first model setting, and $%
c=5$ and $6$ in the second model setting, respectively. Simulations results
were based on $200$ repetitions. For each repetition of each model setting,
we obtained the parameter estimates $\widehat{\boldsymbol{\theta }}$'s based
on the four considered estimation methods: EL, GMM, ET and CU under five
regimes regarding the blocking parameters $L$ and $M$:
\begin{equation*}
\begin{split}
& \text{Regime~(i)}.~~L=M=1;~~~~~~~~~~~~\text{Regime~(ii)}.~~M=\lfloor
n^{1/5}\rfloor ~\text{and}~L=\lfloor 0.5M\rfloor ; \\
& \text{Regime~(iii)}.~~L=M=\lfloor n^{1/5}\rfloor ;~~~~\text{Regime~(iv)}%
.~~M=\lfloor 3n^{1/5}\rfloor ~\text{and}~L=\lfloor 0.5M\rfloor ; \\
& \text{Regime~(v)}.~~L=M=\lfloor 3n^{1/5}\rfloor ~.
\end{split}%
\end{equation*}%
Regime (i) means no blocking. Regimes (ii) and (iv) assigned the block size $%
M$ to be twice of the block separation parameter $L$; and Regimes (iii) and
(v) prescribed $M=L$. For each repetition of the second model setting, we
additionally considered the parameter estimates $\widehat{\boldsymbol{\theta
}}$'s based on the penalized GEL estimation methods. The penalty function $%
p_{\tau }(u)$ used in the simulation satisfied:
\begin{equation*}
p_{\tau }^{\prime }(u)=\tau \bigg\{I(u\leq \tau )+\frac{(a\tau -u)_{+}}{%
(a-1)\tau }I(u>\tau )\bigg\}
\end{equation*}%
for $u>0$, where $a=3.7$, and $s_{+}=s$ for $s>0$ and $0$ otherwise. This
penalty function is given in \cite{FanLi_2001}. We applied the method given
in \cite{LengTang_2011_Biometrika} to determine the penalty parameter $\tau $%
. In each simulation replication, we calculated the $L_{2}$ distance between
$\widehat{\boldsymbol{\theta }}$ and $\boldsymbol{\theta }_{0}$ as $\Vert
\widehat{\boldsymbol{\theta }}-\boldsymbol{\theta }_{0}\Vert _{2}=\{(%
\widehat{\boldsymbol{\theta }}-\boldsymbol{\theta }_{0})^{\prime }(\widehat{%
\boldsymbol{\theta }}-\boldsymbol{\theta }_{0})\}^{1/2}$.

Tables 1 and 2 report empirical medians of the squared estimation errors for
the EL, ET, CU and GMM estimators in the first simulation model with $c=10$
and $c=12$, respectively. And Tables 3 and 4 summarize the empirical median
for the second simulation model with the extra penalized GEL estimators. We
had also collected the average of the squared estimation errors, which
exhibited similar patterns as the empirical median. Hence, we only report
the median of squared estimation errors per the suggestion of one referee.

It is noted that the performance of each estimator at each given blocking
regime was improved when the sample size was increased, which confirms the
convergence of these estimators. For the second nonlinear model, we observed
that the performance of three GEL estimators and their penalized analogues
were improved under the blocking regimes (ii)-(v) which were bona fide
blocking since $L,M>1$. This was not that surprising since dependence was
presence in both simulated models, and applying the blocking can improve the
efficiency of the estimation. However, the performance of the GMM estimator
were largely similar regardless of the blocking regimes used. The empirical
medians of the squared estimation errors of the GMM estimator were much
larger than those of the GEL estimators, which confirmed the existing
research on GMM versus GEL for finite fixed dimensional settings %
\citep{NeweySmith_2004_Ecma, Anatolyev_2005}. Among the three GEL
estimators, we observed that while they were largely similar under the first
simulation model, the EL and the ET estimators performed better than the CU
estimator for the logistic regression model. This might be due to the
multivariate asymmetry in the moment conditions, which makes the bias term
of the CU estimator more pronounced, as shown in %
\cite{NeweySmith_2004_Ecma} and \cite{Anatolyev_2005}. We note that the estimation efficiency among the GEL estimator with respect to the different regimes of the blocking width selection was largely comparable to each other for the simple mean models. However, in the case of the generalized linear model, the regimes (iv) and (v), with the block width $M=\lfloor 3n^{1/5}\rfloor $, led to the best performance. We also observed that under the second model setting where the parameter is sparse, the penalized GEL estimators were much more efficient than their non-penalized counterparts, which confirmed our Theorem 6.

\section{Conclusion}

In this paper, we have investigated the asymptotic properties of the GEL
estimator, the GEL ratio statistic and the over-identification specification
test for high dimensional moment restriction models with increasing number
of parameters and weakly dependent data. We have also investigated a
penalized GEL approach that is designed for the high dimensional sparse
parameter situation with $p>r$, although the true but unknown number of
non-zero parameters is not larger than $r$. We establish the oracle property
of the penalized GEL estimator. Both theoretical and simulation studies find
the penalization leads to efficiency gain for the GEL estimators even for
dependent data.

We establish the consistency and the asymptotic normality of the
high-dimensional GEL and the penalized GEL estimators allowing for fixed
block size $M$ for time series data. However, when the unconditional moment
functions $\{g(X_{t},\boldsymbol{\theta }_{0})\}_{t=1}^{n}$ are
autocorrelated, the simple limiting distributions of the GEL ratio statistic
and the over-identification specification test are established when the
block size $M$ diverges with the sample size $n$. How to practically select $%
M$ is a quite challenging problem. As indicated in \cite%
{HallHorowitzJing_1995} and \cite{Lahari_2003}, although there has been much
research in determining the order of magnitude of $M$, there is in general a
lack of research for selecting the tuning parameter, the coefficient of $M$
for general nonlinear time series models. The simulation study reported in
Section 8 shows that $M=\lfloor 3n^{1/5}\rfloor $ led to satisfactory
performance. Instead of blocking, one could also perform local smoothing of
the unconditional moment functions $\{g(X_{t},\boldsymbol{\theta }%
_{0})\}_{t=1}^{n}$ to reduce temproal dependence
\citep{Smith_1997,
Anatolyev_2005, Kitamura_2007}, which introduces an alternative tuning
parameter, however. We leave it to future research about the performance of
this local smoothing GEL approach for high dimensinal time series models.

\section*{Acknowledgements}

The authors are grateful to the co-editor Jianqing Fan, an Associate Editor,
two anonymous referees and Qiwei Yao for constructive comments. Some results
of this work is the Chapter 2 of Jinyuan Chang's PhD thesis at Peking
University. We thank Jean-Michel Zako\"{\i}an for sharing a longer version
of his paper, and Chenlei Leng and Cheng Yong Tang for sharing their R code.
Jinyuan Chang and Song Xi Chen were supported by National Natural Science
Foundation of China Grants 11131002 and G0113, LMEQF and Center for
Statistical Science at Peking University.

\section*{Appendix}

Throughout the Appendix, $C$ denotes a generic positive finite constant that
may be different in different uses. For any $q=1,\ldots ,Q$ and $k=1,\ldots
,M$, let $\beta _{1}(q,k)=\#\{j<q:X_{(q-1)L+k}\in B_{j}\}$ and $\beta
_{2}(q,k)=\#\{j>q:X_{(q-1)L+k}\in B_{j}\}$. These two quantities denote the
times of the $k$-th element of the $q$-th block occurs in the blocks before
and after the $q$-th block, respectively. Let $\bar{g}(\boldsymbol{\theta })=%
{n}^{-1}\sum_{t=1}^{n}g_{t}(\boldsymbol{\theta })$, $\bar{\phi}(\boldsymbol{%
\theta })={Q}^{-1}\sum_{q=1}^{Q}\phi _{q}(\boldsymbol{\theta })$, $V_{M}=%
\text{Var}\{M^{1/2}{\phi }_{q}(\boldsymbol{\theta }_{0})\},$ $\widehat{%
\Omega }(\boldsymbol{\theta })={Q}^{-1}\sum_{q=1}^{Q}\phi _{q}(\boldsymbol{%
\theta })\phi _{q}(\boldsymbol{\theta })^{\prime }$ and $\Omega (\boldsymbol{%
\theta })={E}\{\phi _{q}(\boldsymbol{\theta })\phi _{q}(\boldsymbol{\theta }%
)^{\prime }\}$.

\subsection*{Some Lemmas I}

The lemmas proposed in this subsection are used to prove Theorem 1.

\begin{la}
\label{la1} $\beta_1(q,k)=(q-1)\wedge\lfloor(M-k)/L\rfloor$ and $%
\beta_2(q,k)=(Q-q)\wedge\lfloor(k-1)/L\rfloor $ .
%as $n$ is sufficiently large.
\end{la}

\noindent\textsc{Proof}: For $t=(q-1)L+k$, suppose $X_t\in B_{\bar{q}}$
where $\bar{q}<q$. Then there exists a positive integer $\bar{k}\in [1,M]$
such that $(q-1)L+k=(\bar{q}-1)L+\bar{k}. $ It means $\bar{q}=q-(\bar{k}%
-k)/L $. From this, we can get $\bar{k}=k+iL$ for some $i\in\{1,\ldots,q-1\}$%
. Note that $\bar{k}\in[1,M]$, then $i\leq\lfloor(M-k)/L\rfloor$. %and $
%K_1\equiv K_2~ (\textrm{mod}~ L).
%$
% Hence, if $K\leq L$, there is only one block $B_{q}$ such that $i\in B_{q}$. Then for any $i=(q(i)-1)L+K$ where $K\leq L$, $\beta_n(i)=1$.
%On the other hand, note $L/M\rightarrow1$, then for sufficient large
%$n$, any integer $K_1\in[L+1,M]$  can be written  as
%$
%K_1=L+l$, where $l=1,\cdots,L-1$.
%Note $
%K_1\equiv K_2~ (\textrm{mod}~ L)
%$, then
%$
%K_2=l\leq L
%$.
%It means for any $i=(q(i)-1)L+K$ where $K\geq L+1$, $\beta_n(i)\leq2$.
Hence, $\beta_1(q,k)=(q-1)\wedge\lfloor(M-k)/L\rfloor $. By the same
argument, $\beta_2(q,k)=(Q-q)\wedge\lfloor(k-1)/L\rfloor $. $\hfill \square$

\begin{la}
\label{la2} Under conditions \textrm{(A.1)(ii)} and \textrm{(A.2)(ii)}, $%
\sup_{\boldsymbol{\theta}\in\Theta}\|\bar{\phi}(\boldsymbol{\theta})-\bar{g}(%
\boldsymbol{\theta})\|_2=O_p(r^{1/2}Mn^{-1}). $
\end{la}

\noindent\textsc{Proof}: %As $
%\bar{\phi}(\theta)-\bar{g}(\theta)=(MQ)^{-1}[\sum_{q=1}^Q\sum_{l\in
%B_q}g_l(\theta)-\sum_{i=1}^ng_i(\theta)]+(n-MQ)(nMQ)^{-1}\sum_{i=1}^ng_i(\theta),
%$ then
By Jensen's inequality,
\begin{equation*}
{E}\bigg\{\sup_{\boldsymbol{\theta}\in\Theta}\|\bar{\phi}(\boldsymbol{\theta}%
)-\bar{g}(\boldsymbol{\theta})\|_2\bigg\}\leq \frac{1}{MQ}\bigg\{%
n-(Q-1)L-M+\sum_{\beta_1(q,k)=0}\beta_2(q,k)+n-MQ\bigg\}\cdot{E}\bigg\{\sup_{%
\boldsymbol{\theta}\in\Theta}\| g_t(\boldsymbol{\theta}) \|_2\bigg\}.
\end{equation*}
From Lemma 1 and (A.1)(ii), $\sum_{\beta_1(q,k)=0}\beta_2(q,k)\leq
(Q-1)(M-L) $ for sufficiently large $n$. Noting that $Q=\lfloor(n-M)/L%
\rfloor+1$, then for sufficiently large $n$
\begin{equation*}
{E}\bigg\{\sup_{\boldsymbol{\theta}\in\Theta}\|\bar{\phi}(\boldsymbol{\theta}%
)-\bar{g}(\boldsymbol{\theta})\|_2\bigg\}\leq2LM^{-1}Q^{-1}\cdot{E}\bigg\{%
\sup_{\boldsymbol{\theta}\in\Theta}\|g_t(\boldsymbol{\theta})\|_2\bigg\}.
\end{equation*}
Hence, (A.1)(ii) and (A.2)(ii) lead to the conclusion. $\hfill \square$

%\bigskip The following lemmas are needed in the proof of Theorem 1 regarding
%the consistency of $\hat{\theta}_{n,p}$.

\begin{la}
\label{la3} Under conditions \textrm{(A.1)(i)} and \textrm{(A.2)(iii)}, $\|%
\widehat{\Omega}(\boldsymbol{\theta}_{0})-\Omega(\boldsymbol{\theta}%
_{0})\|_F=O_p(rM^{1/2}n^{-1/2})$ .
%Furthermore, under Conditions (A.1)-(A.3)
%and (A.5), if $rMn^{-1}=o(1)$, $
%\sup_{\theta\in\Theta}||\hat{\Omega}(\theta)-\Omega(\theta)||_F=o_p\{r^{1/2}\zeta(r,p)\}.
%$
\end{la}

\noindent\textsc{Proof}: Note that
\begin{equation*}
\begin{split}
{E}\{\|\widehat{\Omega}(\boldsymbol{\theta}_{0})-\Omega(\boldsymbol{\theta}%
_{0})\|_F^2\}&=Q^{-1}{E}(\text{tr} \{[\phi_q(\boldsymbol{\theta}_{0})\phi_q(%
\boldsymbol{\theta}_{0})^{\prime }-\Omega(\boldsymbol{\theta}_{0})]^2\}) \\
&~~~+Q^{-2}\sum_{q_1\neq q_2}{E}(\text{tr} \{[\phi_{q_1}(\boldsymbol{\theta}%
_{0})\phi_{q_1}(\boldsymbol{\theta}_{0})^{\prime }-\Omega(\boldsymbol{\theta}%
_{0})][\phi_{q_2}(\boldsymbol{\theta}_{0})\phi_{q_2}(\boldsymbol{\theta}%
_{0})^{\prime }-\Omega(\boldsymbol{\theta}_{0})]\}) \\
&=:A_1+A_2.
\end{split}%
\end{equation*}
As $A_1%
%=&~Q^{-1}\cdot\mathbb{E}\left\{\textrm{tr}\left[\bigg( \phi_{q}(\theta)\phi_{q}(\theta)^T-\Omega(\theta) \bigg)^2\right]\right\}\\
\leq
%Q^{-1}\textrm{tr}\left\{\mathbb{E}\left[\phi_q(\theta)\phi_q(\theta)^T\phi_q(\theta)\phi_q(\theta)^T\right]\right\}
Q^{-1}{E}\{\|\phi_q(\boldsymbol{\theta}_{0})\|_2^4\}, $ by Jensen's
inequality and (A.2)(iii), $A_1=O\left(r^2Mn^{-1}\right).%\label{eq:a1}
$ At the same time,
\begin{equation*}
\begin{split}
A_2&=Q^{-2}\sum_{u,v=1}^r\sum_{q_1\neq q_2} {E}\{ [ \phi_{q_1,u}(\boldsymbol{%
\theta}_{0})\phi_{q_1,v}(\boldsymbol{\theta}_{0})-\Omega_{u,v}(\boldsymbol{%
\theta}_{0}) ][ \phi_{q_2,v}(\boldsymbol{\theta}_{0})\phi_{q_2,u}(%
\boldsymbol{\theta}_{0})-\Omega_{v,u}(\boldsymbol{\theta}_{0}) ]\}, \\
\end{split}%
\end{equation*}
where $\Omega_{u,v}(\boldsymbol{\theta}_{0})$ denotes the $(u,v)$-element of
$\Omega(\boldsymbol{\theta}_{0})$. By Davydov inequality and (A.2)(iii), $%
|A_2|%&\leq8Q^{-2}\sum\limits_{p\neq
%q}\sum\limits_{u,v=1}^r\left\{\mathbb{E}|\phi_{q,u}(\theta)\phi_{q,v}(\theta)-\Omega_{u,v}(\theta) |^{2+\delta}\right\}^{2/(2+\delta)}\cdot \alpha_{\phi}(|p-q|)^{\delta/(2+\delta)}\\
%&=8M^{-4}Q^{-2}\cdot\sum\limits_{i\neq j}\alpha_{X}(|i-j|)^{\delta/(2+\delta)}\sum\limits_{u,v=1}^r\left\{\mathbb{E}|g_{i,u}(\theta)g_{i,v}(\theta)-\mathbb{E}g_{i,u}(\theta)g_{i,v}(\theta)|^{2+\delta}\right\}^{2/(2+\delta)}\\
\leq Cr^{2}Q^{-2}\sum_{q_1\neq q_2}\alpha_{\phi}\{|q_1-q_2|\}^{1-2/\gamma}.$
Hence, by (A.1)(i), $A_2=O(r^{2}Mn^{-1}).%\label{eq:a2}
$ From Markov inequality, $\| \widehat{\Omega}(\boldsymbol{\theta}%
_{0})-\Omega(\boldsymbol{\theta}_{0}) \|_F=O_p(rM^{1/2}n^{-1/2})$. $\hfill
\square$

\begin{la}
\label{la4} Under conditions \textrm{(A.1)(ii)}, \textrm{(A.2)(ii)} and
\textrm{(A.2)(iv)}, then $\sup_{\boldsymbol{\theta}\in\Theta}\lambda_{\max}\{%
\widehat{ \Omega}(\boldsymbol{\theta})\}=O_p(1)%
%+o_p\left(r_n^{\xi/(4+2\delta)}\right).
$ provided that $rMn^{-1}=o(1)$.
\end{la}

\noindent\textsc{Proof}: Using the same approach as in the proof of Lemma %
\ref{la2},
\begin{equation*}
\sup_{\boldsymbol{\theta} \in \Theta }\sup_{\Vert x\Vert _{2}=1}\bigg\{\bigg|%
\frac{1}{MQ}\sum_{q=1}^{Q}\sum_{t\in B_{q}}x^{\prime }g_{t}(\boldsymbol{%
\theta} )g_{t}(\boldsymbol{\theta} )^{\prime }x-\frac{1}{n}%
\sum_{t=1}^{n}x^{\prime }g_{t}(\boldsymbol{\theta} )g_{t}(\boldsymbol{\theta}
)^{\prime }x\bigg|\bigg\}=O_{p}(rMn^{-1}).
\end{equation*}%
%
%
%
%
%
%
%
%
%
%
%
%
%
%
%
%
%
%
%
%
%
%For fixed $\theta\in\Theta$, there exists $x\in\mathbb{R}^r$ with
%unit $L_2$-norm such that
%$\lambda_{\textrm{max}}(\hat{\Omega}(\theta))=x'\hat{\Omega}(\theta)x$.
By Jensen's inequality, for any $\|x\|_2=1$,
\begin{equation*}
\frac{1}{Q}\sum_{q=1}^{Q}x^{\prime }\phi _{q}(\boldsymbol{\theta} )\phi _{q}(%
\boldsymbol{\theta} )^{\prime }x\leq \frac{1}{MQ}\sum_{q=1}^{Q}\sum_{t\in
B_{q}}x^{\prime }g_{t}(\boldsymbol{\theta} )g_{t}(\boldsymbol{\theta}
)^{\prime }x.
\end{equation*}%
Then $\sup_{\boldsymbol{\theta} \in \Theta }\lambda _{\text{max}}\{\widehat{%
\Omega}(\boldsymbol{\theta} )\}\leq \sup_{\boldsymbol{\theta} \in \Theta
}\lambda _{\text{max}}\{n^{-1}\sum_{i=1}^{n}g_{t}(\boldsymbol{\theta} )g_{t}(%
\boldsymbol{\theta} )^{\prime }\}+o_{p}(1)$. The result can be implied by
(A.2)(iv). $\hfill \square $

\begin{la}
\label{la5} Under condition \textrm{(A.2)(ii)}, define $%
\delta_n=o(r^{-1/2}Q^{-1/\gamma})$ and $\Lambda_n=\{\lambda\in\mathbb{R}%
^{r}:\|\lambda\|_2\leq\delta_n\}$, we have $\sup_{1\leq q\leq Q, \boldsymbol{%
\theta}\in\Theta, \lambda\in\Lambda_n}|\lambda^{\prime }\phi_q(\boldsymbol{%
\theta})|\xrightarrow[]{p}0. $ Also w.p.a.1, $\Lambda_n\subset\widehat{%
\Lambda}_n(\boldsymbol{\theta})$ for all $\boldsymbol{\theta}\in\Theta$.
\end{la}

\noindent\textsc{Proof}: From (A.2)(ii) and Markov inequality, $\sup_{1\leq
q\leq Q, \boldsymbol{\theta}\in\Theta}\|\phi_q(\boldsymbol{\theta}%
)\|_2=O_p(r^{1/2}Q^{1/\gamma}). $ Then,
\begin{equation*}
\sup_{1\leq q\leq Q, \boldsymbol{\theta}\in\Theta,
\lambda\in\Lambda_n}|\lambda^{\prime }\phi_q(\boldsymbol{\theta}%
)|\leq\delta_n\cdot\sup_{1\leq q\leq Q, \boldsymbol{\theta}%
\in\Theta}\|\phi_q(\boldsymbol{\theta})\|_2\xrightarrow[]{p}0.
\end{equation*}
It also implies w.p.a.1 $\lambda^{\prime }\phi_q(\boldsymbol{\theta})\in%
\mathcal{V}$ for all $\boldsymbol{\theta}\in\Theta$ and $\|\lambda\|_2\leq%
\delta_n$. $\hfill \square$

\begin{la}
\label{la7} Under conditions \textrm{(A.1)(i)}, \textrm{(A.2)(i)} and
\textrm{(A.2)(iii)}, assume that $\lambda_{\max}(V_M)$ is uniformly bounded
away from infinity. If $r^2M^3n^{-1}=o(1)$, $\|\boldsymbol{\theta}-%
\boldsymbol{\theta}_{0}\|_2=O_p(\tau_n)$ and $rpM\tau_n^2=o(1)$, then $\|%
\widehat{\Omega}(\boldsymbol{\theta})-\widehat{\Omega}(\boldsymbol{\theta}%
_{0})\|_2=O_p( r^{1/2}p^{1/2}M^{-1/2}\tau_n)$.
%and $||\Omega(\theta^*)-\Omega(\tilde{\theta})||_F\leq
%Cr^{1/2}\zeta(r,p)||\theta^*-\tilde{\theta}||_2$.
\end{la}

\noindent\textsc{Proof}: Choose $x\in \mathbb{R}^{r}$ with unit $L_{2}$-norm
such that $\lambda _{\text{max}}\{\widehat{\Omega}(\boldsymbol{\theta} )-%
\widehat{\Omega}(\boldsymbol{\theta} _{0})\}=x^{\prime }\{\widehat{\Omega}(%
\boldsymbol{\theta} )-\widehat{\Omega}(\boldsymbol{\theta} _{0})\}x$. Then,
\begin{equation*}
\begin{split}
&|\lambda _{\text{max}}\{\widehat{\Omega}(\boldsymbol{\theta} )-\widehat{%
\Omega}(\boldsymbol{\theta} _{0})\}| \\
&~~~~~~\leq \frac{1}{Q}\sum_{q=1}^{Q}|x^{\prime }\phi _{q}(\boldsymbol{\theta%
} )-x^{\prime }\phi _{q}(\boldsymbol{\theta} _{0})|\cdot \Vert \phi _{q}(%
\boldsymbol{\theta} )-\phi _{q}(\boldsymbol{\theta} _{0})\Vert _{2}+\frac{2}{%
Q}\sum_{q=1}^{Q}|x^{\prime }\phi _{q}(\boldsymbol{\theta} _{0})|\cdot \Vert
\phi _{q}(\boldsymbol{\theta} )-\phi _{q}(\boldsymbol{\theta} _{0})\Vert _{2}
\\
&~~~~~~ \leq \frac{1}{Q}\sum_{q=1}^{Q}\Vert \phi _{q}(\boldsymbol{\theta}
)-\phi _{q}(\boldsymbol{\theta} _{0})\Vert _{2}^{2}+2[ \lambda _{\text{max}%
}\{\widehat{\Omega}(\boldsymbol{\theta} _{0})\}] ^{1/2}\bigg\{\frac{1}{Q}%
\sum_{q=1}^{Q}\Vert \phi _{q}(\boldsymbol{\theta} )-\phi _{q}(\boldsymbol{%
\theta} _{0})\Vert _{2}^{2}\bigg\}^{1/2}. \\
\end{split}%
\end{equation*}%
Note that $r^{2}M^{3}n^{-1}=o(1)$, by Lemmas \ref{la3} and $%
\lambda_{\max}(V_M)$ is uniformly bounded away from infinity, $\lambda _{%
\text{max}}\{\widehat{\Omega}(\boldsymbol{\theta} _{0})\}=O_{p}(M^{-1})$.
From (A.2)(i), $Q^{-1}\sum_{q=1}^{Q}\Vert \phi _{q}(\boldsymbol{\theta}
)-\phi _{q}(\boldsymbol{\theta} _{0})\Vert _{2}^{2}=rp\cdot O_{p}(\Vert
\boldsymbol{\theta} -\boldsymbol{\theta} _{0}\Vert _{2}^{2})$. If $rpM\tau
_{n}^{2}=o(1)$, then $|\lambda _{\text{max}}\{\widehat{\Omega}(\boldsymbol{%
\theta} )-\widehat{\Omega}(\boldsymbol{\theta} _{0})%
\}|=O_{p}(r^{1/2}p^{1/2}M^{-1/2}\tau _{n})$. Using the same argument, $%
|\lambda _{\text{min}}\{\widehat{\Omega}(\boldsymbol{\theta} )-\widehat{%
\Omega}(\boldsymbol{\theta} _{0})\}|=O_{p}(r^{1/2}p^{1/2}M^{-1/2}\tau _{n})$%
. This completes the proof.$\hfill \square $

\begin{la}
\label{la9} Under conditions \textrm{(A.1)(i)}, \textrm{(A.1)(ii)}, \textrm{%
(A.2)(i)}, \textrm{(A.2)(ii)} and \textrm{(A.2)(iii)}, assume that the
eigenvalues of $V_M$ are uniformly bounded away from zero and infinity. If $%
r^2M^{2-2/\gamma}n^{2/\gamma-1}=o(1)$, $r^2M^{3}n^{-1}=o(1)$, $\|\widetilde{%
\boldsymbol{\theta}}-\boldsymbol{\theta}_{0}\|_2=O_p(\tau_n)$, $%
rpM\tau_n^2=o(1)$ and $\|\bar{g}(\widetilde{\boldsymbol{\theta}}%
)\|_2=O_p(r^{1/2}n^{-1/2})$, then $\widehat{\lambda}(\widetilde{\boldsymbol{%
\theta}})=\arg\max_{\lambda\in\widehat{\Lambda}_n(\widetilde{\boldsymbol{%
\theta}})}\widehat{S}_n(\widetilde{\boldsymbol{\theta}},\lambda)$ exists
w.p.a.1, $\sup_{\lambda\in\widehat{\Lambda}_n(\widetilde{\boldsymbol{\theta}}%
)}\widehat{S}_n(\widetilde{\boldsymbol{\theta}},\lambda)=%
\rho(0)+O_p(rMn^{-1})$ and $\|\widehat{\lambda}(\widetilde{\boldsymbol{\theta%
}})\|_2=O_p(r^{1/2}Mn^{-1/2})$ where $\widehat{S}_n(\boldsymbol{\theta}%
,\lambda)$ is defined in \textrm{(\ref{eq:sn})}.
\end{la}

\noindent\textsc{Proof}: Pick $\delta_n=o(r^{-1/2}Q^{-1/\gamma})$ and $%
r^{1/2}Mn^{-1/2}=o(\delta_n)$, which is guaranteed by $r^2M^{2-2/%
\gamma}n^{2/\gamma-1}=o(1)$. From Lemma \ref{la2} and Triangle inequality,
then $\|\bar{\phi}(\widetilde{\boldsymbol{\theta}})\|_2\leq \|\bar{g}(%
\widetilde{\boldsymbol{\theta}}) \|_2+O_p(r^{1/2}Mn^{-1})$ which implies $\|%
\bar{\phi}(\widetilde{\boldsymbol{\theta}})\|_2=O_p(r^{1/2}n^{-1/2})$. Let $%
\bar{\lambda}=\arg\max_{\lambda\in\Lambda_n}\widehat{S}_n(\widetilde{%
\boldsymbol{\theta}},\lambda)$, where $\Lambda_n$ is defined in Lemma \ref%
{la5}. By Lemmas \ref{la3}, \ref{la5} and \ref{la7}, noting $\rho_{vv}(0)<0$%
,
\begin{equation*}
\begin{split}
\rho(0)=\widehat{S}_n(\widetilde{\boldsymbol{\theta}},0)&\leq\widehat{S}_n(%
\widetilde{\boldsymbol{\theta}},\bar{\lambda})=\rho(0)+\rho_v(0)\bar{\lambda}%
^{\prime }\bar{\phi}(\widetilde{\boldsymbol{\theta}})+\frac{1}{2}\bar{\lambda%
}^{\prime }\bigg\{\frac{1}{Q}\sum\limits_{q=1}^Q\rho_{vv}(\dot{\lambda}%
^{\prime }\phi_q(\widetilde{\boldsymbol{\theta}}))\phi_q(\widetilde{%
\boldsymbol{\theta}})\phi_q(\widetilde{\boldsymbol{\theta}})^{\prime }\bigg\}%
\bar{\lambda} \\
&\leq\rho(0)+|\rho_v(0)|\cdot\|\bar{\lambda}\|_2\cdot\|\bar{\phi}(\widetilde{%
\boldsymbol{\theta}})\|_2-C\|\bar{\lambda}\|_2^2\cdot\{M^{-1}+o_p(M^{-1})\}.
\end{split}%
\end{equation*}
where $\dot{\lambda}$ lies on the jointing line between $\boldsymbol{\mathbf{%
0}}$ and $\bar{\lambda}$. Hence, $\|\bar{\lambda}\|_2\leq C\cdot M\cdot\|%
\bar{\phi}(\widetilde{\boldsymbol{\theta}})\|_2\cdot\{1+o_p(1)%
\}=O_p(r^{1/2}Mn^{-1/2})=o_p(\delta_n). $
%The last step is based on $rM^{1-\zeta}n^{\zeta-1/2}=o(1)$.
Thus $\bar{\lambda}\in\text{int}(\Lambda_n)$ w.p.a.1. Since $\Lambda_n\subset%
\widehat{\Lambda}_n(\widetilde{\boldsymbol{\theta}})$ w.p.a.1, $\widehat{%
\lambda}(\widetilde{\boldsymbol{\theta}})=\bar{\lambda}$ w.p.a.1 by the
concavity of $\widehat{S}_n(\widetilde{\boldsymbol{\theta}},\lambda)$ and $%
\widehat{\Lambda}_n(\widetilde{\boldsymbol{\theta}})$. Then,
\begin{equation*}
\widehat{S}_n(\widetilde{\boldsymbol{\theta}},\widehat{\lambda}(\widetilde{%
\boldsymbol{\theta}}))\leq\rho(0)+|\rho_v(0)|\cdot\|\widehat{\lambda}(%
\widetilde{\boldsymbol{\theta}})\|_2\cdot\|\bar{\phi}(\widetilde{\boldsymbol{%
\theta}})\|_2-C\cdot M^{-1}\cdot\|\widehat{\lambda}(\widetilde{\boldsymbol{%
\theta}})\|_2^2\cdot\{1+o_p(1)\}
\end{equation*}
leads to $\sup_{\lambda\in\widehat{\Lambda}_n(\widetilde{\boldsymbol{\theta}}%
)}\widehat{S}_n(\widetilde{\boldsymbol{\theta}},\lambda)=%
\rho(0)+O_p(rMn^{-1})$. $\hfill \square$

\subsection*{Proof of Theorem 1}

Choose $\delta _{n}=o(r^{-1/2}Q^{-1/\gamma })$ and $r^{1/2}Mn^{-1/2}=o(%
\delta _{n})$. Let $\bar{\lambda}=\text{sign}\{\rho_v(0)\}\cdot\delta _{n}%
\bar{\phi}(\widehat{\boldsymbol{\theta}}_{n})/\Vert \bar{\phi}(\widehat{%
\boldsymbol{\theta}}_{n})\Vert _{2}$, then $\bar{\lambda}\in \Lambda _{n}$.
By Taylor expansion, Lemmas \ref{la4} and \ref{la5}, noting $\rho_{vv}(0)<0$%
,
\begin{equation*}
\begin{split}
\widehat{S}_n(\widehat{\boldsymbol{\theta}}_{n},\bar{\lambda})&
=\rho(0)+\rho_v(0)\bar{\lambda}^{\prime }\bar{\phi}(\widehat{\boldsymbol{%
\theta}}_{n})+\frac{1}{2}\bar{\lambda}^{\prime }\bigg\{\frac{1}{Q}%
\sum\limits_{q=1}^{Q}\rho_{vv}(\dot{\lambda}^{\prime }\phi_q(\widehat{%
\boldsymbol{\theta}}_n))\phi_{q}(\widehat{\boldsymbol{\theta}}_{n})\phi_{q}(%
\widehat{\boldsymbol{\theta}}_{n})^{\prime }\bigg\}\bar{\lambda} \\
& \geq \rho(0)+|\rho_v(0)|\cdot\delta _{n}\cdot \Vert \bar{\phi}(\widehat{%
\boldsymbol{\theta}}_{n})\Vert _{2}-C\cdot O_{p}(1)\cdot \Vert \bar{\lambda}%
\Vert _{2}^{2}.
\end{split}%
\end{equation*}%
%
%
%
%
%
%
%
%
%
%
%
%
%
%
%
%
%
%
%
%
%
%The last step is based on Lemma 4.
Meanwhile, by the same way in the proof of Lemma \ref{la3}, $\Vert \bar{g}(%
\boldsymbol{\theta} _{0})-{E}\{g_{t}(\boldsymbol{\theta} _{0})\}\Vert
_{2}=O_{p}(r^{1/2}n^{-1/2})$. Since ${E}\{g_{t}(\boldsymbol{\theta} _{0})\}=%
\boldsymbol{\mathbf{0}}$, $\Vert \bar{g}(\boldsymbol{\theta} _{0})\Vert
_{2}=O_{p}(r^{1/2}n^{-1/2})$. Then, from Lemma \ref{la9},
\begin{equation*}
\widehat{S}_n(\widehat{\boldsymbol{\theta}}_{n},\bar{\lambda})\leq
\sup\limits_{\lambda \in \widehat{\Lambda}_{n}(\widehat{\boldsymbol{\theta}}%
_{n})}\widehat{S}_n(\widehat{\boldsymbol{\theta}}_{n},\lambda )\leq
\sup\limits_{\lambda \in \widehat{\Lambda}_{n}(\boldsymbol{\theta} _{0})}%
\widehat{S}_n(\boldsymbol{\theta} _{0},\lambda )=\rho(0)+O_{p}(rMn^{-1}).
\end{equation*}%
Hence, $\Vert \bar{\phi}(\widehat{\boldsymbol{\theta}}_{n})\Vert
_{2}=O_{p}(\delta _{n}).$ Consider any $\varepsilon _{n}\rightarrow 0$ and
let $\tilde{\lambda}=\text{sign}\{\rho_v(0)\}\cdot\varepsilon _{n}\bar{\phi}(%
\widehat{\boldsymbol{\theta}}_{n})$, then $\Vert \widetilde{\lambda}\Vert
_{2}=o_{p}(\delta _{n})$. Using the same way above, we can obtain
\begin{equation*}
|\rho_v(0)|\cdot\varepsilon _{n}\cdot \Vert \bar{\phi}(\widehat{\boldsymbol{%
\theta}}_{n})\Vert _{2}^{2}-C\cdot O_{p}(1)\cdot \varepsilon _{n}^{2}\cdot
\Vert \bar{\phi}(\widehat{\boldsymbol{\theta}}_{n})\Vert
_{2}^{2}=O_{p}(rMn^{-1}).
\end{equation*}%
Then, $\varepsilon _{n}\Vert \bar{\phi}(\widehat{\boldsymbol{\theta}}%
_{n})\Vert _{2}^{2}=O_{p}(rMn^{-1})$. Thus, $\Vert \bar{\phi}(\widehat{%
\boldsymbol{\theta}}_{n})\Vert _{2}^{2}=O_{p}(rMn^{-1})$. From Lemma \ref%
{la2}, $\Vert \bar{g}(\widehat{\boldsymbol{\theta}}_{n})\Vert
_{2}=O_{p}(r^{1/2}M^{1/2}n^{-1/2})$. If $\Vert \widehat{\boldsymbol{\theta}}%
_{n}-\boldsymbol{\theta} _{0}\Vert _{2}$ does not converge to zero in
probability, then there exists a subsequence $\{(n_{\ast},M_{\ast},r_{%
\ast},p_{\ast })\}$ such that $\Vert \widehat{\boldsymbol{\theta}}_{n_{\ast
}}-\boldsymbol{\theta} _{0}\Vert _{2}\geq \varepsilon$ a.s. for some
positive constant $\varepsilon$. By (A.1)(iv), $\Vert {E}\{g_{t}(\widehat{%
\boldsymbol{\theta}}_{n_{\ast }})\}\Vert _{2}=o_{p}\{\Delta _{1}(r_{\ast
},p_{\ast })\}+O_{p}(r_{\ast}^{1/2}M_{\ast}^{1/2}n_{\ast}^{-1/2})$. On the
other hand, from (A.1)(iii), $\Vert {E}\{g_{t}(\widehat{\boldsymbol{\theta}}%
_{n_{\ast }})\}\Vert _{2}\geq \Delta _{1}(r_{\ast },p_{\ast })\Delta _{2}(\varepsilon )$. As $\liminf_{r,p\rightarrow\infty}\Delta_1(r,p)>0$,
it is a contradiction. Hence, $\Vert \widehat{\boldsymbol{\theta}}_{n}-%
\boldsymbol{\theta} _{0}\Vert _{2}\xrightarrow[]{p}0$. By (A.2)(iv), $\|\bar{%
g}(\widehat{\boldsymbol{\theta}}_{n})-\bar{g}(\boldsymbol{\theta}%
_{0})\|_2\geq C\|\widehat{\boldsymbol{\theta}}_{n}-\boldsymbol{\theta}%
_{0}\|_2$ w.p.a.1. Then, $\|\widehat{\boldsymbol{\theta}}_{n}-\boldsymbol{%
\theta}_{0}\|_2=O_p(r^{1/2}M^{1/2}n^{-1/2})$. In addition, if $%
r^2pM^{2}n^{-1}=o(1)$, from Lemmas \ref{la3} and \ref{la7}, $\lambda _{\text{%
max}}\{\widehat{\Omega}(\widehat{\boldsymbol{\theta}}_{n})\}\leq CM^{-1}$
w.p.a.1. By repeating the above arguments, we can obtain $\Vert \bar{\phi}(%
\widehat{\boldsymbol{\theta}}_{n})\Vert _{2}=O_{p}(r^{1/2}n^{-1/2})$ and $%
\Vert \widehat{\boldsymbol{\theta}}_{n}-\boldsymbol{\theta} _{0}\Vert
_{2}=O_{p}(r^{1/2}n^{-1/2})$. From Lemma \ref{la9}, $\|\widehat{\lambda}(%
\widehat{\boldsymbol{\theta}}_{n})\|_2=O_p(r^{1/2}Mn^{-1/2})$. Therefore, we
complete the proof of Theorem 1. $\hfill \square $

\subsection*{Proof of Corollary \protect\ref{cy1}}

To construct Corollary \ref{cy1}, we need the analogue of Lemma \ref{la9}
listed below.

\begin{la}
\label{newla1} Under conditions \textrm{(A.1)(i)}, \textrm{(A.1)(ii)},
\textrm{(A.2)(i)}, \textrm{(A.2)(ii)} and \textrm{(A.2)(iii)}, assume that $%
\lambda_{\min}(V_M)\asymp r^{-\iota_1}$ for some $\iota_1>0$.

\begin{enumerate}
\item[\textrm{(a)}.] If $r^{2+3\iota_1}M^{2-2/\gamma}n^{2/\gamma-1}=o(1)$, $%
r^{2+2\iota_1}M^{3}n^{-1}=o(1)$, $\|\widetilde{\boldsymbol{\theta}}-%
\boldsymbol{\theta}_{0}\|_2=O_p(\tau_n)$, $r^{1+2\iota_1}pM\tau_n^2=o(1)$
and $\|\bar{g}(\widetilde{\boldsymbol{\theta}})\|_2=O_p(r^{(1+%
\iota_1)/2}n^{-1/2})$, then $\widehat{\lambda}(\widetilde{\boldsymbol{\theta}%
})=\arg\max_{\lambda\in\widehat{\Lambda}_n(\widetilde{\boldsymbol{\theta}})}%
\widehat{S}_n(\widetilde{\boldsymbol{\theta}},\lambda)$ exists w.p.a.1, $%
\sup_{\lambda\in\widehat{\Lambda}_n(\widetilde{\boldsymbol{\theta}})}%
\widehat{S}_n(\widetilde{\boldsymbol{\theta}},\lambda)=\rho(0)+O_p(r^{1+2%
\iota_1}Mn^{-1})$ and $\|\widehat{\lambda}(\widetilde{\boldsymbol{\theta}}%
)\|_2=O_p(r^{(1+3\iota_1)/2}Mn^{-1/2})$ where $\widehat{S}_n(\boldsymbol{%
\theta},\lambda)$ is defined in \textrm{(\ref{eq:sn})}.

\item[\textrm{(b)}.] If $r^{2+2\iota_1}M^{2-2/\gamma}n^{2/\gamma-1}=o(1)$
and $r^{2+2\iota_1}M^{3}n^{-1}=o(1)$, then $\widehat{\lambda}({\boldsymbol{%
\theta}}_0)=\arg\max_{\lambda\in\widehat{\Lambda}_n({\boldsymbol{\theta}}_0)}%
\widehat{S}_n({\boldsymbol{\theta}}_0 ,\lambda)$ exists w.p.a.1, $%
\sup_{\lambda\in\widehat{\Lambda}_n({\boldsymbol{\theta}}_0)}\widehat{S}_n({%
\boldsymbol{\theta}}_0,\lambda)=\rho(0)+O_p(r^{1+\iota_1}Mn^{-1})$.%
% and $\|\widehat{\lambda}(\widetilde{\btheta})%
%\|_2=O_p(r^{(1+2\iota_1)/2}Mn^{-1/2})$.
\end{enumerate}
\end{la}

\noindent\textsc{Proof}: We first prove (a). Pick $\delta_n=o(r^{-1/2}Q^{-1/%
\gamma})$ and $r^{(1+3\iota_1)/2}Mn^{-1/2}=o(\delta_n)$, which is guaranteed
by $r^{2+3\iota_1}M^{2-2/\gamma}n^{2/\gamma-1}=o(1)$. From Lemma \ref{la2}
and Triangle inequality, then $\|\bar{\phi}(\widetilde{\boldsymbol{\theta}}%
)\|_2\leq \|\bar{g}(\widetilde{\boldsymbol{\theta}})
\|_2+O_p(r^{1/2}Mn^{-1}) $ which implies $\|\bar{\phi}(\widetilde{%
\boldsymbol{\theta}})\|_2=O_p(r^{(1+\iota_1)/2}n^{-1/2})$. Let $\bar{\lambda}%
=\arg\max_{\lambda\in\Lambda_n}\widehat{S}_n(\widetilde{\boldsymbol{\theta}}%
,\lambda)$, where $\Lambda_n$ is defined in Lemma \ref{la5}. By Lemmas \ref%
{la3}, \ref{la5} and \ref{la7}, noting $\rho_{vv}(0)<0$,
\begin{equation*}
\begin{split}
\rho(0)=\widehat{S}_n(\widetilde{\boldsymbol{\theta}},0)&\leq\widehat{S}_n(%
\widetilde{\boldsymbol{\theta}},\bar{\lambda})=\rho(0)+\rho_v(0)\bar{\lambda}%
^{\prime }\bar{\phi}(\widetilde{\boldsymbol{\theta}})+\frac{1}{2}\bar{\lambda%
}^{\prime }\bigg\{\frac{1}{Q}\sum\limits_{q=1}^Q\rho_{vv}(\dot{\lambda}%
^{\prime }\phi_q(\widetilde{\boldsymbol{\theta}}))\phi_q(\widetilde{%
\boldsymbol{\theta}})\phi_q(\widetilde{\boldsymbol{\theta}})^{\prime }\bigg\}%
\bar{\lambda} \\
&\leq\rho(0)+|\rho_v(0)|\cdot\|\bar{\lambda}\|_2\cdot\|\bar{\phi}(\widetilde{%
\boldsymbol{\theta}})\|_2-C\|\bar{\lambda}\|_2^2\cdot\{M^{-1}r^{-%
\iota_1}+o_p(M^{-1}r^{-\iota_1})\}.
\end{split}%
\end{equation*}
where $\dot{\lambda}$ lies on the jointing line between $\boldsymbol{\mathbf{%
0}}$ and $\bar{\lambda}$. Therefore, $\|\bar{\lambda}\|_2\leq C\cdot
Mr^{\iota_1}\cdot\|\bar{\phi}(\widetilde{\boldsymbol{\theta}}%
)\|_2\cdot\{1+o_p(1)\}=O_p(r^{(1+3\iota_1)/2}Mn^{-1/2})=o_p(\delta_n). $
%The last step is based on $rM^{1-\zeta}n^{\zeta-1/2}=o(1)$.
Thus $\bar{\lambda}\in\text{int}(\Lambda_n)$ w.p.a.1. Since $\Lambda_n\subset%
\widehat{\Lambda}_n(\widetilde{\boldsymbol{\theta}})$ w.p.a.1, $\widehat{%
\lambda}(\widetilde{\boldsymbol{\theta}})=\bar{\lambda}$ w.p.a.1 by the
concavity of $\widehat{S}_n(\widetilde{\boldsymbol{\theta}},\lambda)$ and $%
\widehat{\Lambda}_n(\widetilde{\boldsymbol{\theta}})$. Then,
\begin{equation*}
\widehat{S}_n(\widetilde{\boldsymbol{\theta}},\widehat{\lambda}(\widetilde{%
\boldsymbol{\theta}}))\leq\rho(0)+|\rho_v(0)|\cdot\|\widehat{\lambda}(%
\widetilde{\boldsymbol{\theta}})\|_2\cdot\|\bar{\phi}(\widetilde{\boldsymbol{%
\theta}})\|_2-C\cdot M^{-1}r^{-\iota_1}\cdot\|\widehat{\lambda}(\widetilde{%
\boldsymbol{\theta}})\|_2^2\cdot\{1+o_p(1)\}
\end{equation*}
leads to $\sup_{\lambda\in\widehat{\Lambda}_n(\widetilde{\boldsymbol{\theta}}%
)}\widehat{S}_n(\widetilde{\boldsymbol{\theta}},\lambda)=\rho(0)+O_p(r^{1+2%
\iota_1}Mn^{-1})$. The proof of (b) is similar to that of (a). $\hfill
\square$

\medskip

Here, we begin to prove Corollary \ref{cy1}. Choose $\delta
_{n}=o(r^{-1/2}Q^{-1/\gamma })$ and $r^{(1+3\iota_1)/2}Mn^{-1/2}=o(\delta
_{n})$. Let $\bar{\lambda}=\text{sign}\{\rho_v(0)\}\cdot\delta _{n}\bar{\phi}%
(\widehat{\boldsymbol{\theta}}_{n})/\Vert \bar{\phi}(\widehat{\boldsymbol{%
\theta}}_{n})\Vert _{2}$, then $\bar{\lambda}\in \Lambda _{n}$. By Taylor
expansion, Lemmas \ref{la4} and \ref{la5}, noting $\rho_{vv}(0)<0$,
\begin{equation*}
\begin{split}
\widehat{S}_n(\widehat{\boldsymbol{\theta}}_{n},\bar{\lambda})&
=\rho(0)+\rho_v(0)\bar{\lambda}^{\prime }\bar{\phi}(\widehat{\boldsymbol{%
\theta}}_{n})+\frac{1}{2}\bar{\lambda}^{\prime }\bigg\{\frac{1}{Q}%
\sum\limits_{q=1}^{Q}\rho_{vv}(\dot{\lambda}^{\prime }\phi_q(\widehat{%
\boldsymbol{\theta}}_n))\phi_{q}(\widehat{\boldsymbol{\theta}}_{n})\phi_{q}(%
\widehat{\boldsymbol{\theta}}_{n})^{\prime }\bigg\}\bar{\lambda} \\
& \geq \rho(0)+|\rho_v(0)|\cdot\delta _{n}\cdot \Vert \bar{\phi}(\widehat{%
\boldsymbol{\theta}}_{n})\Vert _{2}-C\cdot O_{p}(1)\cdot \Vert \bar{\lambda}%
\Vert _{2}^{2}.
\end{split}%
\end{equation*}%
%
%
%
%
%
%
%
%
%
%
%
%
%
%
%
%
%
%
%
%
%
%The last step is based on Lemma 4.
Meanwhile, by the same way in the proof of Lemma \ref{la3}, $\Vert \bar{g}(%
\boldsymbol{\theta} _{0})-{E}\{g_{t}(\boldsymbol{\theta} _{0})\}\Vert
_{2}=O_{p}(r^{1/2}n^{-1/2})$. Since ${E}\{g_{t}(\boldsymbol{\theta} _{0})\}=%
\boldsymbol{\mathbf{0}}$, $\Vert \bar{g}(\boldsymbol{\theta} _{0})\Vert
_{2}=O_{p}(r^{1/2}n^{-1/2})$. Then, from Lemma \ref{newla1}(b),
\begin{equation*}
\widehat{S}_n(\widehat{\boldsymbol{\theta}}_{n},\bar{\lambda})\leq
\sup\limits_{\lambda \in \widehat{\Lambda}_{n}(\widehat{\boldsymbol{\theta}}%
_{n})}\widehat{S}_n(\widehat{\boldsymbol{\theta}}_{n},\lambda )\leq
\sup\limits_{\lambda \in \widehat{\Lambda}_{n}(\boldsymbol{\theta} _{0})}%
\widehat{S}_n(\boldsymbol{\theta} _{0},\lambda
)=\rho(0)+O_{p}(r^{1+\iota_1}Mn^{-1}).
\end{equation*}%
Hence, $\Vert \bar{\phi}(\widehat{\boldsymbol{\theta}}_{n})\Vert
_{2}=O_{p}(\delta _{n}).$ Consider any $\varepsilon _{n}\rightarrow 0$ and
let $\tilde{\lambda}=\text{sign}\{\rho_v(0)\}\cdot\varepsilon _{n}\bar{\phi}(%
\widehat{\boldsymbol{\theta}}_{n})$, then $\Vert \widetilde{\lambda}\Vert
_{2}=o_{p}(\delta _{n})$. Using the same way above, we can obtain
\begin{equation*}
|\rho_v(0)|\cdot\varepsilon _{n}\cdot \Vert \bar{\phi}(\widehat{\boldsymbol{%
\theta}}_{n})\Vert _{2}^{2}-C\cdot O_{p}(1)\cdot \varepsilon _{n}^{2}\cdot
\Vert \bar{\phi}(\widehat{\boldsymbol{\theta}}_{n})\Vert
_{2}^{2}=O_{p}(r^{1+\iota_1}Mn^{-1}).
\end{equation*}%
Then, $\varepsilon _{n}\Vert \bar{\phi}(\widehat{\boldsymbol{\theta}}%
_{n})\Vert _{2}^{2}=O_{p}(r^{1+\iota_1}Mn^{-1})$. Thus, $\Vert \bar{\phi}(%
\widehat{\boldsymbol{\theta}}_{n})\Vert _{2}^{2}=O_{p}(r^{1+\iota_1}Mn^{-1})$%
. From Lemma \ref{la2}, $\Vert \bar{g}(\widehat{\boldsymbol{\theta}}%
_{n})\Vert _{2}=O_{p}(r^{(1+\iota_1)/2}M^{1/2}n^{-1/2})$. If $\Vert \widehat{%
\boldsymbol{\theta}}_{n}-\boldsymbol{\theta} _{0}\Vert _{2}$ does not
converge to zero in probability, then there exists a subsequence $%
\{(n_{\ast},M_{\ast},r_{\ast},p_{\ast })\}$ such that $\Vert \widehat{%
\boldsymbol{\theta}}_{n_{\ast }}-\boldsymbol{\theta} _{0}\Vert _{2}\geq
\varepsilon$ a.s. for some positive constant $\varepsilon$. By (A.1)(iv), $%
\Vert {E}\{g_{t}(\widehat{\boldsymbol{\theta}}_{n_{\ast }})\}\Vert
_{2}=o_{p}\{\Delta _{1}(r_{\ast },p_{\ast
})\}+O_{p}(r_{\ast}^{(1+\iota_1)/2}M_{\ast}^{1/2}n_{\ast}^{-1/2})$. On the
other hand, from (A.1)(iii), $\Vert {E}\{g_{t}(\widehat{\boldsymbol{\theta}}%
_{n_{\ast }})\}\Vert _{2}\geq \Delta _{1}(r_{\ast },p_{\ast })\Delta _{2}(\varepsilon )$. As $\liminf_{r,p\rightarrow\infty}\Delta_1(r,p)>0$,
it is a contradiction. Hence, $\Vert \widehat{\boldsymbol{\theta}}_{n}-%
\boldsymbol{\theta} _{0}\Vert _{2}\xrightarrow[]{p}0$. By (A.2)(iv), $\|\bar{%
g}(\widehat{\boldsymbol{\theta}}_{n})-\bar{g}(\boldsymbol{\theta}%
_{0})\|_2\geq C\|\widehat{\boldsymbol{\theta}}_{n}-\boldsymbol{\theta}%
_{0}\|_2$ w.p.a.1. Then, $\|\widehat{\boldsymbol{\theta}}_{n}-\boldsymbol{%
\theta}_{0}\|_2=O_p(r^{(1+\iota_1)/2}M^{1/2}n^{-1/2})$. In addition, if $%
r^{2+\iota_1}pM^{2}n^{-1}=o(1)$, from Lemmas \ref{la3} and \ref{la7}, $%
\lambda _{\text{max}}\{\widehat{\Omega}(\widehat{\boldsymbol{\theta}}%
_{n})\}\leq CM^{-1}$ w.p.a.1. By repeating the above arguments, we can
obtain $\Vert \bar{\phi}(\widehat{\boldsymbol{\theta}}_{n})\Vert
_{2}=O_{p}(r^{(1+\iota_1)/2}n^{-1/2})$ and $\Vert \widehat{\boldsymbol{\theta%
}}_{n}-\boldsymbol{\theta} _{0}\Vert _{2}=O_{p}(r^{(1+\iota_1)/2}n^{-1/2})$.
From Lemma \ref{newla1}, $\|\widehat{\lambda}(\widehat{\boldsymbol{\theta}}%
_{n})\|_2=O_p(r^{(1+3\iota_1)/2}Mn^{-1/2})$. Therefore, we complete the
proof of Corollary \ref{cy1}. $\hfill \square $

\subsection*{Some Lemmas II}

The lemmas proposed in this subsection are used to establish Proposition 1,
Proposition 2 and Theorem 2. The proof of Proposition 1 is based on the
asymptotic expansion given in Proposition 2, so we will first construct the
proof of Proposition 2 later.
%Multiplying $M^{-1}$ on both sides of (\ref{eq:key}) and expanding $\bar{\phi}(\hat{\theta}%

\begin{la}
\label{la10} Under conditions \textrm{(A.1)-(A.2)}, assume that $%
\lambda_{\max}(V_M)$ is uniformly bounded away from infinity. If $%
r^2M^{2-2/\gamma}n^{2/\gamma-1}=o(1)$, $r^2pM^2n^{-1}=o(1)$ and $%
r^2M^3n^{-1}=o(1)$, then for any $x\in\mathbb{R}^p$, $y,z\in\mathbb{R}^r$,
\begin{equation*}
\begin{split}
&~~\bigg\|\frac{1}{Q}\sum_{q=1}^Q\rho_v(\widehat{\lambda}(\widehat{%
\boldsymbol{\theta}}_n)^{\prime }\phi_q(\widehat{\boldsymbol{\theta}}%
_n))\nabla_{\boldsymbol{\theta}}\phi_q(\widehat{\boldsymbol{\theta}}_{n}) x-%
\frac{1}{Q}\sum_{q=1}^Q\rho_v(0)\nabla_{\boldsymbol{\theta}}\phi_q(\widehat{%
\boldsymbol{\theta}}_{n})x\bigg\|_2=O_p(rp^{1/2}M^{1/2}n^{-1/2})\cdot\|x\|_2,
\\
&~~~~~~~~~~~\bigg|\frac{M}{Q}\sum_{q=1}^Qy^{\prime }\rho_{vv}(\widetilde{%
\lambda}^{\prime }\phi_q(\widehat{\boldsymbol{\theta}}_{n}))\phi_q(\widehat{%
\boldsymbol{\theta}}_{n})\phi_q(\widehat{\boldsymbol{\theta}}_{n})^{\prime
}z-\frac{M}{Q}\sum_{q=1}^Q\rho_{vv}(0)y^{\prime }\phi_q(\widehat{\boldsymbol{%
\theta}}_{n})\phi_q(\widehat{\boldsymbol{\theta}}_{n})^{\prime }z\bigg| \\
&~~~~~~~~~~~~~~~~~~~~~=O_p(rM^{1-1/\gamma}n^{1/\gamma-1/2})\cdot\|y\|_2\cdot%
\|z\|_2,
\end{split}%
\end{equation*}
where $\widehat{\lambda}(\widehat{\boldsymbol{\theta}}_{n})$ and $\widetilde{%
\lambda}$ are defined in \textrm{(\ref{eq:key})}.
\end{la}

\noindent\textsc{Proof}: From Theorem 1, both $\widehat{\lambda}(\widehat{%
\boldsymbol{\theta}}_{n})$ and $\widetilde{\lambda}$ are $%
O_p(r^{1/2}Mn^{-1/2})=o_p(\delta_n)$ where $\delta_n$ is defined in Lemma %
\ref{la5}. By Taylor expansion and Cauchy-Schwarz inequality,
\begin{equation*}
\begin{split}
&\bigg\|\frac{1}{Q}\sum_{q=1}^Q\rho_v(\widehat{\lambda}(\widehat{\boldsymbol{%
\theta}}_n)^{\prime }\phi_q(\widehat{\boldsymbol{\theta}}_n))\nabla_{%
\boldsymbol{\theta}}\phi_q(\widehat{\boldsymbol{\theta}}_{n}) x-\frac{1}{Q}%
\sum_{q=1}^Q\rho_v(0)\nabla_{\boldsymbol{\theta}}\phi_q(\widehat{\boldsymbol{%
\theta}}_{n})x\bigg\|_2^2 \\
&~~~~~~~~~~\leq\bigg[\frac{1}{Q}\sum_{q=1}^Q\rho_{vv}^2(\dot{\lambda}%
^{\prime }\phi_q(\widehat{\boldsymbol{\theta}}_{n}))\{\widehat{\lambda}(%
\widehat{\boldsymbol{\theta}}_{n})^{\prime }\phi_q(\widehat{\boldsymbol{%
\theta}}_{n})\}^2\bigg]\bigg[\frac{1}{Q}\sum_{q=1}^Qx^{\prime }\{\nabla_{%
\boldsymbol{\theta}}\phi_q(\widehat{\boldsymbol{\theta}}_{n})\}^{\prime
}\{\nabla_{\boldsymbol{\theta}}\phi_q(\widehat{\boldsymbol{\theta}}_{n})\}x %
\bigg]
\end{split}%
\end{equation*}
where $\dot{\lambda}$ lies on the jointing line between $\boldsymbol{\mathbf{%
0}}$ and $\widehat{\lambda}(\widehat{\boldsymbol{\theta}}_n)$. From Lemma %
\ref{la5} and $\lambda_{\text{max}}\{\widehat{\Omega}(\widehat{\boldsymbol{%
\theta}}_{n})\}=O_p(M^{-1})$ which is provided by Lemmas \ref{la3} and \ref%
{la7}, we obtain
\begin{equation*}
\sum_{q=1}^Q\rho_{vv}^2(\dot{\lambda}^{\prime }\phi_q(\widehat{\boldsymbol{%
\theta}}_{n}))\big\{\widehat{\lambda}(\widehat{\boldsymbol{\theta}}%
_{n})^{\prime }\phi_q(\widehat{\boldsymbol{\theta}}_{n})\big\}^2\leq
C\sum_{q=1}^Q\big\{\widehat{\lambda}(\widehat{\boldsymbol{\theta}}%
_{n})^{\prime }\phi_q(\widehat{\boldsymbol{\theta}}_{n})\big\}%
^2\cdot\{1+o_p(1)\}=O_p(rMn^{-1}).
\end{equation*}
% If $\zeta_1(r,p)rMn^{-1}=o(1)$,
%$r^2M^3n^{-1}=o(1)$ and
%$\zeta_2^2(r)M^{-1}\sum_{k=1}^Mk\alpha_X(k)^{1-2/\gamma}=o(1)$, by
%From Lemmas 3, 6 and 7,
%$\lambda_{\textrm{max}}(\hat{\Omega}(\hat{\theta}_{EL}))=O_p(M^{-1})$.
%Then,
%\[
%\frac{1}{Q}\sum_{q=1}^Q\frac{(\hat{\lambda}'\hat{\phi}_q)^2}{(1-c_q\hat{\lambda}'\hat{\phi}_q)^4}=O_p\{M^{-1}||\hat{\lambda}||_2^2\}=O_p(rMn^{-1}).
%\]
On the other hand,
\begin{equation*}
\frac{1}{Q}\sum_{q=1}^Qx^{\prime }\big\{\nabla_{\boldsymbol{\theta}}\phi_q(%
\widehat{\boldsymbol{\theta}}_{n})\big\}^{\prime }\big\{\nabla_{\boldsymbol{%
\theta}}\phi_q(\widehat{\boldsymbol{\theta}}_{n})\big\}x\leq \frac{1}{MQ}%
\sum_{q=1}^Q\sum_{t\in B_q}\big\|\nabla_{\boldsymbol{\theta}} g_t(\widehat{%
\boldsymbol{\theta}}_{n})\cdot x\big\|_2^2=O_p(rp)\cdot\|x\|_2^2.
\end{equation*}
Using the same argument, we can obtain the other result. $\hfill \square$

\begin{la}
\label{la11} Under conditions \textrm{(A.1)(i)}, \textrm{(A.1)(ii)} and
\textrm{(A.3)}, then $\|\nabla_{\boldsymbol{\theta}}\bar{\phi}(\boldsymbol{%
\theta})-\nabla_{\boldsymbol{\theta}}\bar{\phi}(\boldsymbol{\theta}%
^*)\|_F=O_p(r^{1/2}p\cdot\|\boldsymbol{\theta}-\boldsymbol{\theta}^*\|_2)$
for any $\boldsymbol{\theta}$, $\boldsymbol{\theta}^*$ in a neighborhood of $%
\boldsymbol{\theta}_{0}$, and $\|\nabla_{\boldsymbol{\theta}}\bar{\phi}(%
\boldsymbol{\theta}_{0})-{E}\{ \nabla_{\boldsymbol{\theta}} g_t(\boldsymbol{%
\theta}_{0})\}\|_F=O_p(r^{1/2}p^{1/2}n^{-1/2})$ provided that $M=o(n^{1/2})$.
\end{la}

\noindent\textsc{Proof}: Using Taylor expansion and noting (A.3), the first
conclusion holds. Using the same method in the proof of Lemma \ref{la3}, $%
\|\nabla_{\boldsymbol{\theta}} \bar{g}(\boldsymbol{\theta}_{0})-{E}\{\nabla_{%
\boldsymbol{\theta}} g_t(\boldsymbol{\theta}_{0})\}%
\|_F=O_p(r^{1/2}p^{1/2}n^{-1/2}). $ By the same way in the proof of Lemma %
\ref{la2}, $\|\nabla_{\boldsymbol{\theta}}\bar{\phi}(\boldsymbol{\theta}%
_{0})-\nabla_{\boldsymbol{\theta}} \bar{g}(\boldsymbol{\theta}%
_{0})\|_F=O_p(r^{1/2}p^{1/2}Mn^{-1}). $ Hence, by Triangle inequality, we
can obtain $\|\nabla_{\boldsymbol{\theta}}\bar{\phi}(\boldsymbol{\theta}%
_{0})-{E}\{\nabla_{\boldsymbol{\theta}} g_t(\boldsymbol{\theta}%
_{0})\}\|_F=O_p(r^{1/2}p^{1/2}n^{-1/2}). $ $\hfill \square$

%\begin{la}\label{la12}
%Under conditions (A.1)(i), (A.2)(ii) and (A.2)(iv), then for any
%$x\in\mathbb{R}^r$,
%\begin{equation*}
%\big\|M^{-1}\Omega^{-1}(\btheta_{0})x-V_n^{-1}x\big\|_2\leq CrM^{-1}%
%\sum_{k=1}^Mk\alpha_X(k)^{1-2/\gamma}\cdot \|x\|_2
%\end{equation*}
%provided that $rM^{-1}%
%\sum_{k=1}^Mk\alpha_X(k)^{1-2/\gamma}=o(1)$.
%\end{la}
%
%\noindent\textsc{Proof}: Let $u=M^{-1}\Omega^{-1}(\btheta_{0})x$ and $%
%v=V_n^{-1}M^{-1}\Omega^{-1}(\btheta_{0})\{M\Omega(\btheta_{0})-V_n\}V_n^{-1}x
%$, then
%\begin{equation*}
%\big\|M^{-1}\Omega^{-1}(\btheta_{0})x-V_n^{-1}x\big\|_2^2 =u^{\prime
%}\big\{M\Omega(\btheta_{0})-V_n\big\}v.
%\end{equation*}
%By Lemma \ref{la6}, we can obtain the result. $\hfill \square$
%
%
%

\subsection*{Proof of Proposition 2}

Define
\begin{equation*}
\boldsymbol{\beta}=([{E}\{\nabla_{\boldsymbol{\theta}} g_t(\boldsymbol{\theta%
}_{0})\}]^{\prime }V_M^{-1}V_nV_M^{-1}[{E}\{\nabla_{\boldsymbol{\theta}} g_t(%
\boldsymbol{\theta}_{0})\}])^{-1/2} \boldsymbol{\alpha}_{n},
\end{equation*}
then
\begin{equation*}
\begin{split}
\|{E}\{\nabla_{\boldsymbol{\theta}} g_t(\boldsymbol{\theta}_{0})\}\cdot%
\boldsymbol{\beta}\|_2^2&=\boldsymbol{\alpha}_{n}^{\prime }(U^{\prime}U)^{
-1/2}U^{\prime }V_n^{-1/2}V_M^2V_n^{-1/2}U(U^{\prime}U)^{ -1/2}\boldsymbol{%
\alpha}_{n} \\
&\leq \lambda_{\text{max}}(V_n^{-1/2}V_M^2V_n^{-1/2})\cdot\|U(U^{\prime}U)^{
-1/2}\boldsymbol{\alpha}_{n}\|_2^2= \lambda_{\text{max}}^2(V_M)\lambda_{%
\min}^{-1}(V_n),
\end{split}%
\end{equation*}
where $U=V_n^{1/2}V_M^{-1}[{E}\{\nabla_{\boldsymbol{\theta}} g_t(\boldsymbol{%
\theta}_{0})\}]$. Therefore, $\|{E}\{\nabla_{\boldsymbol{\theta}} g_t(%
\boldsymbol{\theta}_{0})\}\cdot\boldsymbol{\beta}\|_2=O(1)$. Meanwhile,
\begin{equation*}
\begin{split}
\|\boldsymbol{\beta}\|_2^2\leq&~ \lambda_{\text{min}}^{-1}([{E}\{\nabla_{%
\boldsymbol{\theta}} g_t(\boldsymbol{\theta}_{0})\}]^{\prime
}V_M^{-1}V_nV_M^{-1}[{E}\{\nabla_{\boldsymbol{\theta}} g_t(\boldsymbol{\theta%
}_{0})\}]) \\
\leq&~\lambda_{\text{max}}^2(V_M)\lambda_{\text{min}}^{-1}([{E}\{\nabla_{%
\boldsymbol{\theta}} g_t(\boldsymbol{\theta}_{0})\}]^{\prime }[{E}\{\nabla_{%
\boldsymbol{\theta}} g_t(\boldsymbol{\theta}_{0})\}])\lambda_{%
\min}^{-1}(V_n).
\end{split}%
\end{equation*}
Hence, $\|\boldsymbol{\beta}\|_2\leq C$. From Lemma \ref{la5},
\begin{equation*}
\frac{M}{Q}\sum_{q=1}^Q\rho_{vv}(\widetilde{\lambda}^{\prime }\phi_q(%
\widehat{\boldsymbol{\theta}}_{n}))\phi_q(\widehat{\boldsymbol{\theta}}%
_{n})\phi_q(\widehat{\boldsymbol{\theta}}_{n})^{\prime }=\rho_{vv}(0)M%
\widehat{\Omega}(\widehat{\boldsymbol{\theta}}_{n})\cdot\{1+o_p(1)\}.
\end{equation*}
Noting Lemmas \ref{la3} and \ref{la7}, we know the eigenvalues of $M\widehat{%
\Omega}(\widehat{\boldsymbol{\theta}}_{n})$ are uniformly bounded away from
zero and infinity w.p.a.1. Hence, the eigenvalues of $MQ^{-1}\sum_{q=1}^Q%
\rho_{vv}(\widetilde{\lambda}^{\prime }\phi_q(\widehat{\boldsymbol{\theta}}%
_{n}))\phi_q(\widehat{\boldsymbol{\theta}}_{n})\phi_q(\widehat{\boldsymbol{%
\theta}}_{n})^{\prime }$ are uniformly bounded away from zero and infinity
w.p.a.1. By Lemma \ref{la10} and (\ref{eq:key}),
\begin{equation*}
\boldsymbol{\beta}^{\prime }\{\nabla_{\boldsymbol{\theta}}\bar{\phi}(%
\widehat{\boldsymbol{\theta}}_{n})\}^{\prime }\bigg\{\frac{M}{Q}%
\sum_{q=1}^Q\rho_{vv}(\widetilde{\lambda}^{\prime }\phi_q(\widehat{%
\boldsymbol{\theta}}_{n}))\phi _q(\widehat{\boldsymbol{\theta}}_{n})\phi_q(%
\widehat{\boldsymbol{\theta}}_{n})^{\prime } \bigg\}^{-1}\bar{\phi}(\widehat{%
\boldsymbol{\theta}}_n)=O_p(r^{3/2}p^{1/2}M^{1/2}n^{-1}).
\end{equation*}
From Lemmas \ref{la11} and \ref{la10},
\begin{equation*}
\begin{split}
&~\boldsymbol{\beta}^{\prime }[{E}\{\nabla_{\boldsymbol{\theta}} g_t(%
\boldsymbol{\theta}_{0})\}]^{\prime }\{M\widehat{\Omega}(\widehat{%
\boldsymbol{\theta}}_{n})\}^{-1}\bar{\phi}(\widehat{\boldsymbol{\theta}}_n)
\\
=&~O_p(
r^{3/2}p^{1/2}M^{1/2}n^{-1})+O_p(r^{3/2}M^{1-1/\gamma}n^{1/%
\gamma-1})+O_p(r^{3/2}pn^{-1}).
\end{split}%
\end{equation*}
Note that Lemmas \ref{la3} and \ref{la7},
\begin{equation*}
\begin{split}
&~\boldsymbol{\beta}^{\prime }[{E}\{\nabla_{\boldsymbol{\theta}} g_t(%
\boldsymbol{\theta}_{0})\}]^{\prime }V_M^{-1}\bar{\phi}(\widehat{\boldsymbol{%
\theta}}_n) \\
=&~O_p(r^{3/2}p^{1/2}M^{1/2}n^{-1})+O_p(r^{3/2}M^{1-1/\gamma}n^{1/%
\gamma-1})+O_p(r^{3/2}M^{3/2}n^{-1})+O_p(r^{3/2}pn^{-1}).
\end{split}%
\end{equation*}
Expanding $\bar{\phi}(\widehat{\boldsymbol{\theta}}_n)$ around $\boldsymbol{%
\theta}=\boldsymbol{\theta}_0$, by Lemmas \ref{la11} and \ref{la2},
\begin{equation*}
\begin{split}
&~\boldsymbol{\beta}^{\prime }[{E}\{\nabla_{\boldsymbol{\theta}} g_t(%
\boldsymbol{\theta}_{0})\}]^{\prime }V_M^{-1}[{E}\{\nabla_{\boldsymbol{\theta%
}} g_t(\boldsymbol{\theta}_{0})\}](\widehat{\boldsymbol{\theta}}_{n}-%
\boldsymbol{\theta}_{0}) \\
=&-\boldsymbol{\beta}^{\prime }[{E}\{\nabla_{\boldsymbol{\theta}} g_t(%
\boldsymbol{\theta}_{0})\}]^{\prime }V_M^{-1}\bar{g}(\boldsymbol{\theta}%
_{0})+O_p(
r^{3/2}p^{1/2}M^{1/2}n^{-1})+O_p(r^{3/2}M^{1-1/\gamma}n^{1/\gamma-1}) \\
&+O_p(r^{3/2}pn^{-1})+O_p(r^{3/2}M^{3/2}n^{-1}).
\end{split}%
\end{equation*}
Hence, we obtain Proposition 2.$\hfill \square$

\subsection*{Proof of Proposition 1}

From Proposition 2, if we pick $\boldsymbol{\alpha}_{n}=(\widehat{%
\boldsymbol{\theta}}_{n}-\boldsymbol{\theta}_{0})/\|\widehat{\boldsymbol{%
\theta}}_{n}-\boldsymbol{\theta}_{0}\|_2$, we can obtain that
\begin{equation*}
\begin{split}
&~\sqrt{n}\cdot\lambda_{\max}^{-1/2}\{[{E}\{\nabla_{\boldsymbol{\theta}} g_t(%
\boldsymbol{\theta}_{0})\}]^{\prime }V_M^{-1}V_nV_M^{-1}[{E}\{\nabla_{%
\boldsymbol{\theta}} g_t(\boldsymbol{\theta}_{0})\}]\} \\
&~~~~~~~~~~~~~~~~~~~~~~\cdot\lambda_{\text{min}}\{[{E}\{\nabla_{\boldsymbol{%
\theta}} g_t(\boldsymbol{\theta}_{0})\}]^{\prime }V_M^{-1}[{E}\{\nabla_{%
\boldsymbol{\theta}} g_t(\boldsymbol{\theta}_{0})\}]\}\cdot\|\widehat{%
\boldsymbol{\theta}}_{n}-\boldsymbol{\theta}_{0}\|_2 \\
=&~O_p\{\|\sqrt{n}\boldsymbol{\alpha}_{n}^{\prime }([{E}\{\nabla_{%
\boldsymbol{\theta}} g_t(\boldsymbol{\theta}_{0})\}]^{\prime
}V_M^{-1}V_nV_M^{-1}[{E}\{\nabla_{\boldsymbol{\theta}} g_t(\boldsymbol{\theta%
}_{0})\}])^{-1/2}[{E}\{\nabla_{\boldsymbol{\theta}} g_t(\boldsymbol{\theta}%
_{0})\}]^{\prime }V_M^{-1}\bar{g}(\boldsymbol{\theta}_{0})\|_2\} \\
&+O_p(r^{3/2}p^{1/2}M^{1/2}n^{-1/2})+O_p(r^{3/2}pn^{-1/2})+O_p(r^{3/2}M^{1-1/\gamma}n^{1/\gamma-1/2})+O_p(r^{3/2}M^{3/2}n^{-1/2}).
\label{eq:theta}
\end{split}%
\end{equation*}
Note that
\begin{equation*}
\lambda_{\max}^{-1/2}\{[{E}\{\nabla_{\boldsymbol{\theta}} g_t(\boldsymbol{%
\theta}_{0})\}]^{\prime }V_M^{-1}V_nV_M^{-1}[{E}\{\nabla_{\boldsymbol{\theta}%
} g_t(\boldsymbol{\theta}_{0})\}]\}\lambda_{\text{min}}\{[{E}\{\nabla_{%
\boldsymbol{\theta}} g_t(\boldsymbol{\theta}_{0})\}]^{\prime }V_M^{-1}[{E}%
\{\nabla_{\boldsymbol{\theta}} g_t(\boldsymbol{\theta}_{0})\}]\}>C,
\end{equation*}
which is assumed in (A.2)(iv) and the eigenvalues of $V_M$ and $V_n$ are
uniformly bounded away from zero and infinity. Therefore, by
\begin{equation*}
{E}(\|\sqrt{n}\boldsymbol{\alpha}_{n}^{\prime }([{E}\{\nabla_{\boldsymbol{%
\theta}} g_t(\boldsymbol{\theta}_{0})\}]^{\prime }V_M^{-1}V_nV_M^{-1}[{E}%
\{\nabla_{\boldsymbol{\theta}} g_t(\boldsymbol{\theta}_{0})\}])^{-1/2}[{E}%
\{\nabla_{\boldsymbol{\theta}} g_t(\boldsymbol{\theta}_{0})\}]^{\prime }V_M^{-1}\bar{g}(\boldsymbol{\theta}_{0})\|_2^2)=p,
\end{equation*}
we complete the proof of Proposition 1.$\hfill \square$

\subsection*{Proof of Theorem 2}

From Proposition 2, it is only need to show
\begin{equation*}
S_n:=-\sqrt{n}\boldsymbol{\alpha}_{n}^{\prime }([{E}\{\nabla_{\boldsymbol{%
\theta}} g_t(\boldsymbol{\theta}_{0})\}]^{\prime }V_M^{-1}V_nV_M^{-1}[{E}%
\{\nabla_{\boldsymbol{\theta}} g_t(\boldsymbol{\theta}_{0})\}])^{-1/2}[{E}%
\{\nabla_{\boldsymbol{\theta}} g_t(\boldsymbol{\theta}_{0})\}]^{\prime
}V_M^{-1}\bar{g}(\boldsymbol{\theta}_{0})\xrightarrow{d}N(0,1).
\end{equation*}
Let
\begin{equation*}
x_{n,t}=-\boldsymbol{\alpha}_{n}^{\prime }([{E}\{\nabla_{\boldsymbol{\theta}%
} g_t(\boldsymbol{\theta}_{0})\}]^{\prime }V_M^{-1}V_nV_M^{-1}[{E}\{\nabla_{%
\boldsymbol{\theta}} g_t(\boldsymbol{\theta}_{0})\}])^{-1/2}[{E}\{\nabla_{%
\boldsymbol{\theta}} g_t(\boldsymbol{\theta}_{0})\}]^{\prime }V_M^{-1}\bar{g}%
(\boldsymbol{\theta}_{0})=:\boldsymbol{\beta}_n^{\prime }g_t(\boldsymbol{%
\theta}_{0}),
\end{equation*}
then $S_n=n^{-1/2}\sum_{t=1}^nx_{n,t}$. As restriction (\ref{eq:moment-cond}%
) holds, $\sup_{n}\sup_{1\leq t\leq n}{E}(|x_{n,t}|^\gamma)<\infty$. On the
other hand, $\text{Var}(S_n)=1$. Note that (A.1)(i), then by the central
limit theorem proposed in Francq and Zako\"{\i}an (2005), we have $S_n%
\xrightarrow{d}N(0,1)$.$\hfill \square$

\subsection*{Some Lemmas III}

To prove Theorem 3, we employ the blocking technique by splitting the
observations to big blocks of length $h$ and small blocks of length $b$.
Suppose that $\tilde{B}_i=(X_{(i-1)(h+b)+1},\ldots,X_{i(h+b)})=(\tilde{B}%
_{i1},\tilde{B}_{i2})$, where $\tilde{B}_{i1}=(X_{(i-1)(h+b)+1},%
\ldots,X_{(i-1)(h+b)+h})$ , $\tilde{B}_{i2}=(X_{(i-1)(h+b)+(h+1)},%
\ldots,X_{i(h+b)})$ and $b<h$. Then $n=T(h+b)+m$, where $m<h+b$. Later, we
will discuss the selection of $b$ and $h$. By a similar argument to those in
finding the order of $\|\bar{g}(\boldsymbol{\theta}_{0})\|_2$ and the proof
of Lemma \ref{la2}, we can obtain $\|\bar{g}(\boldsymbol{\theta}%
_{0})-(Th)^{-1}\sum_{i=1}^T\sum_{t\in\tilde{B}_{i1}}g_t(\boldsymbol{\theta}%
_{0})\|_2=O_p(r^{1/2}T^{1/2}b^{1/2}n^{-1}) $. Furthermore, define $%
\widetilde{V}_n=\text{Var}\{h^{-1/2}\sum_{t=1}^hg_t(\boldsymbol{\theta}%
_{0})\}$, $Z_{T,i}=h^{-1/2}\widetilde{V}_n^{-1/2}\sum_{t\in\tilde{B}%
_{i1}}g_t(\boldsymbol{\theta}_{0})$ and $G_{T,k}=\sum_{i=1}^kZ_{T,i}$, then $%
{E}(Z_{T,i})=\boldsymbol{\mathbf{0}}$ and ${E}( Z_{T,i}Z_{T,i}^{\prime
})=I_r $. It can be shown that $|x^{\prime }(V_n-\widetilde{V}_n)y|\leq
Crh^{-1}\sum_{k=1}^hk\alpha_{X}(k)^{1-2/\gamma}\cdot\|x\|_2\|y\|_2 $. Define
$\mathscr{G}_{T,0}=\{\varnothing,\Omega\}$, $\mathscr{G}_{T,k}=%
\sigma(Z_{T,1},\ldots,Z_{T,k})$, $k=1,\ldots,T$, and $%
S_{T,k}=T^{-1}(2r)^{-1/2}\{(\sum_{i=1}^kZ_{T,i}^{\prime
})(\sum_{i=1}^kZ_{T,i})-kr\}$. ${E}_{T,k}(\cdot)$ denote the conditional
expectation given $\mathscr{G}_{T,k}$. Let $%
D_{T,k}=S_{T,k}-S_{T,k-1}=T^{-1}(2r)^{-1/2}(2Z_{T,k}^{\prime
}G_{T,k-1}+\|Z_{T,k}\|_2^2-r)$. The following lemmas are used to establish
Theorem 3.

\begin{la}
\label{la13} Under conditions \textrm{(A.1)(i)} and \textrm{(A.2)(iii)},
assume that the eigenvalues of $V_n$ are uniformly bounded away from zero
and infinity. Then
\begin{equation*}
\begin{split}
{E}\{(\|Z_{T,k}\|_2^2-r)^2\}&\leq Cr^2h^2, \\
{E}\{|{E}_{T,k-1}(Z_{T,k}^{\prime }G_{T,k-1})|\}&\leq Crk^{1/2}\bigg\{%
\sum_{l=1}^h\alpha_X(b+l)^{1-2/\gamma}\bigg\}^{1/2}, \\
{E}\{{E}_{T,k-1}^2(\|Z_{T,k}\|_2^2-r)\}&\leq
Cr^2h^2\alpha_X(b+1)^{1-2/\gamma},
\end{split}%
\end{equation*}
and for any $i\neq j$,
\begin{equation*}
|{E}\{(\|Z_{T,i}\|_2^2-r)(\|Z_{T,j}\|_2^2-r)\}|\leq
Cr^2h^2\alpha_X\{(b+h)|i-j|-h+1\}^{1-2/\gamma}, \\
\end{equation*}
provided that $rh^{-1}\sum_{k=1}^hk\alpha_X(k)^{1-2/\gamma}=o(1)$. %\[
%\left|\mathbb{E}\{[G_{T,i-1}'Z_{T,i}-\mathbb{E}_{T,i-1}G_{T,i-1}'Z_{T,i}][G_{T,j-1}'Z_{T,j}-\mathbb{E}_{T,j-1}G_{T,j-1}'Z_{T,j}]\}\right|\leq
%Cr^2h^2ij\alpha_X\{(b+h)|i-j|-h+1\}^{1-2/\gamma}.
%\]
\end{la}

\noindent\textsc{Proof}: As $rh^{-1}\sum_{k=1}^hk\alpha_X(k)^{1-2/%
\gamma}=o(1)$, $\lambda_{\text{max}}(\widetilde{V}_n^{-1})\leq C$. Then, ${E}%
(\|Z_{T,k}\|_2^4)\leq Ch^{-2}{E}\{ \|\sum_{t=1}^hg_t(\boldsymbol{\theta}%
_{0})\|_2^4\}$. By Triangle and Jensen's inequalities, ${E}%
\{\|\sum_{t=1}^hg_t(\boldsymbol{\theta}_{0})\|_2^4\}\leq Cr^2h^4$. Hence, ${E%
}\{(\|Z_{T,k}\|_2^2-r)^2\}\leq Cr^2h^2$. By Cauchy-Schwarz inequality,
\begin{equation*}
{E}\{|{E}_{T,k-1}(Z_{T,k}^{\prime }G_{T,k-1}) |\}\leq\{{E}%
(\|G_{T,k-1}\|_2^2)\}^{1/2}[{E}\{ \|{E}_{T,k-1}(Z_{T,k})\|_2^2\}]^{1/2}.
\end{equation*}
Using the same method in the proof of Lemma \ref{la3}, ${E}(
\|G_{T,k-1}\|_2^2)\leq Crk$. On the other hand, note that $\|{E}%
_{T,k-1}(Z_{T,k})\|_2^2\leq C\sum_{t\in\tilde{B}_{k1}}\|{E}_{T,k-1}\{g_t(%
\boldsymbol{\theta}_{0})\}\|_2^2$. Hence,
\begin{equation}
\begin{split}
{E}\{\|{E}_{T,k-1}(Z_{T,k})\|_2^2\}&~\leq C\sum_{t\in\tilde{B}%
_{k1}}\sum_{j=1}^r{E}[{E}^2_{T,k-1}\{g_{tj}(\boldsymbol{\theta}_{0})\}] \\
&~\leq C\sum_{t\in\tilde{B}_{k1}}\sum_{j=1}^r\left[{E}\{|g_{tj}(\boldsymbol{%
\theta}_{0})|^\gamma\}\right]^{2/\gamma}\alpha_X\{t+b-(k-1)(h+b)\}^{1-2/%
\gamma} \\
&~= Cr\sum_{l=1}^h\alpha_X(b+l)^{1-2/\gamma}.  \label{eq:4}
\end{split}%
\end{equation}
%\end{spacing}
%\begin{spacing}{1.4}
\noindent This is based on the fact that if ${E}(X)=0$, then $({E}[\{{E}(X|%
\mathcal{F})\}^2])^{1/2}=\sup\{{E}(XY):Y\in\mathcal{F},{E}(Y^2)=1\}$ for any
$\sigma$-field $\mathcal{F}$ (details can be found in Durrett (2010)), and
Davydov inequality. Then, ${E}\{|{E}_{T,k-1}(Z_{T,k}^{\prime
}G_{T,k-1})|\}\leq Crk^{1/2}\{\sum_{l=1}^h\alpha_X(b+l)^{1-2/\gamma}\}^{1/2}$%
. Using the same argument above, we can obtain ${E}(\|Z_{T,k}\|_2^{2\gamma})%
\leq Cr^{\gamma}h^{\gamma}$. Then, by the same argument of (\ref{eq:4}),
\begin{equation}
{E}\{{E}_{T,k-1}^2(\|Z_{T,k}\|_2^2-r)\}\leq
Cr^2h^2\alpha_X(b+1)^{1-2/\gamma}.  \label{eq:t}
\end{equation}
For any $i\neq j$, by Davydov inequality,
\begin{equation*}
|{E}\{(\|Z_{T,i}\|_2^2-r)(\|Z_{T,j}\|_2^2-r)\}|\leq C\{{E}%
(|\|Z_{T,i}\|_2^2-r|^{\gamma})\}^{2/\gamma}\alpha_X\{(b+h)|i-j|-h+1\}^{1-2/%
\gamma}.
\end{equation*}
Hence, we complete the proof of this lemma. $\hfill \square$

\begin{la}
\label{la14} Under conditions \textrm{(A.1)(i)} and \textrm{(A.2)(iii)},
assume that the eigenvalues of $V_n$ are uniformly bounded away from zero
and infinity. Then
\begin{equation*}
\begin{split}
{E}\{\|{E}_{T,k-1}(Z_{T,k}Z_{T,k}^{\prime })-I_r\|_F^2\}&\leq
Cr^2h^2\alpha_X(b+1)^{1-2/\gamma}, \\
{E}[\{G_{T,k-1}^{\prime }{E}_{T,k-1}(Z_{T,k})\}^2]&\leq
Cr^2h^2k^2\alpha_X(b+1)^{1/2-1/\gamma}.
\end{split}%
\end{equation*}
provided that $rh^{-1}\sum_{k=1}^hk\alpha_X(k)^{1-2/\gamma}=o(1)$.
\end{la}

\noindent\textsc{Proof}: Note that
\begin{equation*}
\begin{split}
&~{E}\{\|{E}_{T,k-1}(Z_{T,k}Z_{T,k}^{\prime })-I_r\|_F^2\} \\
\leq&~C{E}\bigg(\bigg\|{E}_{T,k-1}\bigg\{\bigg[h^{-1/2}\sum_{t\in\tilde{B}%
_{k1}}g_t(\boldsymbol{\theta}_{0})\bigg]\bigg[h^{-1/2}\sum_{t\in\tilde{B}%
_{k1}}g_t(\boldsymbol{\theta}_{0})\bigg]^{\prime }\bigg\}-\widetilde{V}_n%
\bigg\|_F^2\bigg) \\
\leq&~C\sum_{u,v=1}^r{E}\bigg(\bigg|{E}_{T,k-1}\bigg\{\bigg[%
h^{-1/2}\sum_{t\in\tilde{B}_{k1}}g_{tu}(\boldsymbol{\theta}_{0})\bigg]\bigg[%
h^{-1/2}\sum_{t\in\tilde{B}_{k1}}g_{tv}(\boldsymbol{\theta}_{0})\bigg]\bigg\}%
-\widetilde{V}_n(u,v)\bigg|^2\bigg) \\
\leq&~Cr^2h^2\alpha_X(b+1)^{1-2/\gamma},
\end{split}%
\end{equation*}
where $\widetilde{V}_n(u,v)$ denotes the $(u,v)$-element of $\widetilde{V}_n$%
. The last step is similar to (\ref{eq:t}). Then we obtain the first
conclusion. As
\begin{equation*}
\begin{split}
\|{E}_{T,k-1}(Z_{T,k})\|_2^4\leq&~\{{E}_{T,k-1}(Z_{T,k})\}^{\prime }\{{E}%
_{T,k-1}(Z_{T,k}Z_{T,k}^{\prime })\}\{{E}_{T,k-1}(Z_{T,k})\} \\
\leq&~ \|{E}_{T,k-1}(Z_{T,k})\|_2^2\|{E}_{T,k-1}(Z_{T,k}Z_{T,k}^{\prime
})-I_r\|_F+\|{E}_{T,k-1}(Z_{T,k})\|_2^2,
\end{split}%
\end{equation*}
then by Cauchy-Schwarz inequality,
\begin{equation*}
\begin{split}
&~{E}\{\|{E}_{T,k-1}(Z_{T,k})\|_2^4\} \\
\leq&~[{E}\{\|{E}_{T,k-1}(Z_{T,k})\|_2^4\}]^{1/2}[{E}\{\|{E}%
_{T,k-1}(Z_{T,k}Z_{T,k}^{\prime })-I_r\|_F^2\}]^{1/2}+{E}\{\|{E}%
_{T,k-1}(Z_{T,k})\|_2^2\}.
\end{split}%
\end{equation*}
Let $U=[{E}\{\|{E}_{T,k-1}(Z_{T,k})\|_2^4\}]^{1/2}$, from (\ref{eq:4}) and
the first result in this lemma,
\begin{equation*}
U^2\leq CUrh\alpha_{X}(b+1)^{1/2-1/\gamma}+Crh\alpha_X(b+1)^{1-2/\gamma}.
\end{equation*}
Then, $U\leq Crh\alpha_X(b+1)^{1/2-1/\gamma}$. Hence, ${E}\{\|{E}%
_{T,k-1}(Z_{T,k})\|_2^4\}\leq Cr^2h^2\alpha_X(b+1)^{1-2/\gamma}$. Also,
\begin{equation*}
{E}[\{G_{T,k-1}^{\prime }{E}_{T,k-1}(Z_{T,k})\}^2]\leq \{{E}%
(\|G_{T,k-1}\|_2^4)\}^{1/2}[{E}\{\|{E}_{T,k-1}(Z_{T,k})\|_2^4\} ]^{1/2}\leq
Cr^2h^2k^2\alpha_X(b+1)^{1/2-1/\gamma}.
\end{equation*}
We complete the proof.$\hfill \square$

\begin{la}
\label{la15} Under conditions \textrm{(A.1)(i)} and \textrm{(A.2)(iii)},
assume that the eigenvalues of $V_n$ are uniformly bounded away from zero
and infinity. Then
%$r^{-1}n^{1+\kappa_2}\alpha_X(b)^{1/2-1/\gamma}=o(1)$ and
%$r^{-2}n^{1+\kappa_2}\sum_{k=1}^T\alpha_X\{(b+h)k-h\}^{1-2/\gamma}=o(1)$,
$r^{-1}T^{-2}\sum_{j=2}^T(T-j)G_{T,j-1}^{\prime }Z_{T,j}=o_p(1) $ provided
that $rh^{-1}\sum_{k=1}^hk\alpha_X(k)^{1-2/\gamma}=o(1)$ and $%
n^2h^2\alpha_X(b+1)^{1-2/\gamma}=o(1)$.
\end{la}

\noindent\textsc{Proof}: %Note that
%\[
%\sum_{j=2}^T(T-j)G_{T,j-1}'Z_{T,j}=\sum_{j=2}^T(T-j)[G_{T,j-1}'Z_{T,j}-\mathbb{E}_{T,j-1}G_{T,j-1}'Z_{T,j}]+\sum_{j=2}^T(T-j)\mathbb{E}_{T,j-1}G_{T,j-1}'Z_{T,j}.
%\]
Note that
\begin{equation*}
\begin{split}
&{E}\bigg(\bigg[\frac{1}{rT^2}\sum_{j=2}^T(T-j)\{G_{T,j-1}^{\prime }Z_{T,j}-{%
E}_{T,j-1}(G_{T,j-1}^{\prime }Z_{T,j})\}\bigg]^2\bigg) \\
&~~~~~~~~~~~~ =\frac{1}{r^2T^4}\sum_{j=2}^T(T-j)^2{E}[\{G_{T,j-1}^{\prime
}Z_{T,j}-{E}_{T,j-1}(G_{T,j-1}^{\prime }Z_{T,j})\}^2].
\end{split}%
\end{equation*}
By the first result of Lemma \ref{la14},
\begin{equation*}
\begin{split}
{E}\{(G_{T,j-1}^{\prime }Z_{T,j})^2\} &={E}(\|G_{T,j-1}\|_2^2)+{E}%
[G_{T,j-1}^{\prime }\{{E}_{T,j-1}(Z_{T,j}Z_{T,j}^{\prime })-I_r\}G_{T,j-1}]
\\
&\leq Crj+\{{E}(\|G_{T,j-1}\|_2^4)\}^{1/2}[{E}\{\|{E}%
_{T,j-1}(Z_{T,j}Z_{T,j}^{\prime })-I_r\|_F^2\}]^{1/2} \\
&\leq Crj+Cr^2h^2j^2\alpha_X(b+1)^{1/2-1/\gamma}.
\end{split}%
\end{equation*}
Using Cauchy-Schwarz inequality and the fact ${E}\{|{E}_{T,j-1}(G_{T,j-1}^{%
\prime }Z_{T,j})|^2\}\leq {E}\{(G_{T,j-1}^{\prime }Z_{T,j})^2\}$,
\begin{equation*}
{E}\bigg[\frac{1}{r^2T^4}\sum_{j=2}^T(T-j)^2\{G_{T,j-1}^{\prime }Z_{T,j}-{E}%
_{T,j-1}(G_{T,j-1}^{\prime }Z_{T,j})\}^2\bigg]\leq
Cr^{-1}+Cnh\alpha_X(b+1)^{1/2-1/\gamma}\rightarrow0.
\end{equation*}
Then, $r^{-1}T^{-2}\sum_{j=2}^T(T-j)\{G_{T,j-1}^{\prime }Z_{T,j}-{E}%
_{T,j-1}(G_{T,j-1}^{\prime }Z_{T,j})\}= o_p(1)$. From Lemma \ref{la13}, we
have $r^{-1}T^{-2}\sum_{j=2}^T(T-j){E}_{T,j-1}(G_{T,j-1}^{\prime }Z_{T,j})=
o_p(1)$. Hence, we complete the proof. $\hfill \square$

\begin{la}
\label{la16} Under conditions \textrm{(A.1)(i)} and \textrm{(A.2)(iii)},
assume that the eigenvalues of $V_n$ are uniformly bounded away from zero
and infinity. If $rh^{-1}\sum_{k=1}^hk\alpha_X(k)^{1-2/\gamma}=o(1)$, $%
rT\sum_{l=1}^h\alpha_X(b+l)^{1-2/\gamma}=o(1)$ and $rh^2\alpha_X(b+1)^{1-2/%
\gamma}=o(1)$, then $\sum_{k=1}^T{E}_{T,k-1}(D_{T,k})= o_p(1)$.
\end{la}

\noindent\textsc{Proof}: Note that
\begin{equation*}
\begin{split}
\sum_{k=1}^T{E}_{T,k-1}(D_{T,k})&=\frac{2}{(2r)^{1/2}T}\sum_{k=1}^T{E}%
_{T,k-1}(Z_{T,k}^{\prime }G_{T,k-1})+\frac{1}{(2r)^{1/2}T}\sum_{k=1}^T{E}%
_{T,k-1}(\|Z_{T,k}\|_2^2-r) \\
&=:I_1+I_2.
\end{split}%
\end{equation*}
From Lemma \ref{la13},
\begin{equation*}
{E}(|I_1|)\leq\frac{2}{(2r)^{1/2}T}\sum_{k=1}^T{E}\{|{E}_{T,k-1}(Z_{T,k}^{%
\prime }G_{T,k-1})|\}\leq C r^{1/2}T^{1/2}\bigg\{\sum_{l=1}^h%
\alpha_X(b+l)^{1-2/\gamma}\bigg\}^{1/2}\rightarrow0
\end{equation*}
and
\begin{equation*}
{E}(|I_2|)\leq\frac{1}{(2r)^{1/2}T}\sum_{k=1}^T[{E}\{{E}^2_{T,k-1}(\|Z_{T,k}%
\|_2^2-r)\}]^{1/2}\leq Cr^{1/2}h\alpha_X(b+1)^{1/2-1/\gamma}\rightarrow0.
\end{equation*}
Then, we complete the proof of this lemma. $\hfill \square$

\begin{la}
\label{la17} Under conditions \textrm{(A.1)(i)} and \textrm{(A.2)(iii)},
assume that the eigenvalues of $V_n$ are uniformly bounded away from zero
and infinity. If $rh^{-1}\sum_{k=1}^hk\alpha_X(k)^{1-2/\gamma}=o(1)$, $%
r^2n^{2}h^2\alpha_X(b+1)^{1-2/\gamma}=o(1)$ and $rh^3n^{-1}=o(1)$, then $%
S_{T,T}\xrightarrow[]{d}N(0,1)$.
\end{la}

\noindent\textsc{Proof}: We will use the martingale central limit theorem to
show $S_{T,T}\xrightarrow[]dN(0,1)$. Note that
\begin{equation*}
S_{T,T}=\sum_{k=1}^TD_{T,k}=\sum_{k=1}^T\{D_{T,k}-{E}_{T,k-1}(D_{T,k})\}+%
\sum_{k=1}^T{E}_{T,k-1}(D_{T,k}).
\end{equation*}
The first part on the right hand of above equation are the sum of a sequence
of martingale difference with respect to $\{\mathscr{G}_{T,k}\}_{k=0}^T$.
From Lemma \ref{la16}, $S_{T,T}=\sum_{k=1}^T\{D_{T,k}-{E}_{T,k-1}(D_{T,k})%
\}+o_p(1). $

By the martingale central limit theorem (Billingsley, 1995), in order to
show the conclusion, it is sufficient to show that, letting $\sigma_{T,k}^2={%
E}_{T,k-1}[\{D_{T,k}-{E}_{T,k-1}(D_{T,k})\}^2]$, as $T\rightarrow\infty$,
\begin{equation*}
V_{T,T}:=\sum_{k=1}^T\sigma_{T,k}^2\xrightarrow[]p1~~~~\text{and}%
~~~~\sum_{k=1}^T{E}\{D_{T,k}-{E}_{T,k-1}(D_{T,k})\}^4\rightarrow0.
\end{equation*}

For the first part,
\begin{equation*}
\begin{split}
V_{T,T}=&~\frac{2}{rT^2}\sum_{k=1}^T(G_{T,k-1}^{\prime }[{E}%
_{T,k-1}(Z_{T,k}Z_{T,k}^{\prime })-\{{E}_{T,k-1}(Z_{T,k})\}\{{E}%
_{T,k-1}(Z_{T,k}^{\prime })\}]G_{T,k-1}) \\
&+\frac{2}{rT^2}\sum_{k=1}^TG_{T,k-1}^{\prime }[{E}_{T,k-1}\{Z_{T,k}(%
\|Z_{T,k}\|_2^2-r)\}-{E}_{T,k-1}(Z_{T,k})\cdot{E}_{T,k-1}(\|Z_{T,k}\|_2^2-r)]
\\
&+\frac{1}{2rT^2}\sum_{k=1}^T[{E}_{T,k-1}\{(\|Z_{T,k}\|_2^2-r)^2\}-{E}%
_{T,k-1}^2(\|Z_{T,k}\|_2^2-r)] \\
=:&~I_1+I_2+I_3.
\end{split}%
\end{equation*}
We will show that $I_1\xrightarrow[]p1, I_2\xrightarrow[]p0$ and $I_3%
\xrightarrow[]p0$.

Note that $0\leq I_3 \leq (2r)^{-1}T^{-2}\sum_{k=1}^T{E}_{T,k-1}\{(\|Z_{T,k}%
\|_2^2-r)^2\} $ and ${E}\{(\|Z_{T,k}\|_2^2-r)^2\}\leq Cr^2h^2$, then $I_3%
\xrightarrow[]p0$. Using Cauchy-Schwarz, Triangle and Jensen's inequalities,
\begin{equation*}
\begin{split}
&|{E}_{T,k-1}\{G_{T,k-1}^{\prime }Z_k(\|Z_{T,k}\|_2^2-r)\}| \\
&~~~~~~~~\leq\{{E}_{T,k-1}(G_{T,k-1}^{\prime }Z_{T,k}Z_{T,k}^{\prime
}G_{T,k-1})\}^{1/2}[{E}_{T,k-1}\{(\|Z_{T,k}\|_2^2-r)^2\}]^{1/2} \\
&~~~~~~~~\leq\{\|G_{T,k-1}\|_2^2+\|G_{T,k-1}\|_2^2\cdot\|{E}%
_{T,k-1}(Z_{T,k}Z_{T,k}^{\prime })-I_r\|_F\}^{1/2}[{E}_{T,k-1}\{(\|Z_{T,k}%
\|_2^2-r)^2\}]^{1/2} \\
&~~~~~~~~\leq C\|G_{T,k-1}\|_2[{E}_{T,k-1}\{(\|Z_{T,k}\|_2^2-r)^2\}]^{1/2}%
\{1+\|{E}_{T,k-1}(Z_{T,k}Z_{T,k}^{\prime })-I_r\|_F^{1/2}\}.
\end{split}%
\end{equation*}
Then, by Cauchy-Schwarz inequality and Lemmas \ref{la13} and \ref{la14},
\begin{equation*}
\begin{split}
&{E}[|{E}_{T,k-1}\{G_{T,k-1}^{\prime }Z_k(\|Z_{T,k}\|_2^2-r)\}|] \\
&~~~~~~~~~~~~~\leq C\{{E}(\|G_{T,k-1}\|_2^2)\}^{1/2}[{E}\{(\|Z_{T,k}%
\|_2^2-r)^2\}]^{1/2} \\
&~~~~~~~~~~~~~~~~~+C\{{E}(\|G_{T,k-1}\|_2^2)\}^{1/2}[{E}\{(\|Z_{T,k}%
\|_2^2-r)^4\}]^{1/4}[{E}\{\|{E}_{T,k-1}(Z_{T,k}Z_{T,k}^{\prime
})-I_r\|_F^2\}]^{1/4} \\
&~~~~~~~~~~~~~\leq
Cr^{3/2}hk^{1/2}+Cr^2h^{3/2}k^{1/2}\alpha_X(b+1)^{1/4-1/(2\gamma)}.
\end{split}%
\end{equation*}
Hence, $r^{-1}T^{-2}\sum_{k=1}^T{E}[|G_{T,k-1}^{\prime }{E}
_{T,k-1}\{Z_{T,k}(\|Z_{T,k}\|_2^2-r)\}|]\rightarrow0. $ By the same
argument, we can obtain %\[
%\left|\mathbb{E}_{T,k-1}G_{T,k-1}'Z_{T,k}\cdot\mathbb{E}_{T,k-1}[||Z_{T,k}||_2^2-r]\right|\leq\{\mathbb{E}_{T,k-1}(G_{T,k-1}'Z_{T,k})^2\}^{1/2}\{\mathbb{E}_{T,k-1}[||Z_{T,k}||_2^2-r]^2\}^{1/2}.
%\]
%Repeating preceding proof,
$r^{-1}T^{-2}\sum_{k=1}^T{E}\{|{E}_{T,k-1}(G_{T,k-1}^{\prime }Z_{T,k})\cdot{E%
}_{T,k-1}(\|Z_{T,k}\|_2^2-r)|\}\rightarrow0$. Thus, $I_2\xrightarrow[]p0$.
Note that
\begin{equation*}
\begin{split}
I_1=&~\frac{2}{rT^2}\sum_{k=1}^T\|G_{T,k-1}\|_2^2 \\
&+\frac{2}{rT^2}\sum_{k=1}^TG_{T,k-1}^{\prime }\{{E}%
_{T,k-1}(Z_{T,k}Z_{T,k}^{\prime })-{E}_{T,k-1}(Z_{T,k})\cdot{E}%
_{T,k-1}(Z_{T,k}^{\prime })-I_r\}G_{T,k-1}.
\end{split}%
\end{equation*}
By Triangle inequality,
\begin{equation*}
\begin{split}
&\bigg|\frac{2}{rT^2}\sum_{k=1}^TG_{T,k-1}^{\prime }\{{E}%
_{T,k-1}(Z_{T,k}Z_{T,k}^{\prime })-{E}_{T,k-1}(Z_{T,k})\cdot{E}%
_{T,k-1}(Z_{T,k}^{\prime })-I_r\}G_{T,k-1}\bigg| \\
&~~~~~~~~\leq\frac{2}{rT^2}\sum_{k=1}^T\|G_{T,k-1}\|_2^2\|{E}%
_{T,k-1}(Z_{T,k}Z_{T,k}^{\prime })-I_r\|_F+\frac{2}{rT^2}\sum_{k=1}^T%
\{G_{T,k-1}^{\prime }{E}_{T,k-1}(Z_{T,k})\}^2.
\end{split}%
\end{equation*}
By Cauchy-Schwarz inequality and Lemma \ref{la14},
\begin{equation*}
\begin{split}
&~{E}\{\|G_{T,k-1}\|_2^2\|{E}_{T,k-1}(Z_{T,k}Z_{T,k}^{\prime })-I_r\|_F\} \\
\leq&~\{{E}(\|G_{T,k-1}\|_2^4)\}^{1/2}[{E}\{\|{E}_{T,k-1}(Z_{T,k}Z_{T,k}^{%
\prime })-I_r\|_F^2\}]^{1/2} \\
\leq&~ Cr^2h^2k^2\alpha_X(b+1)^{1/2-1/\gamma}.
\end{split}%
\end{equation*}
Then,
\begin{equation*}
{E}\bigg\{\frac{2}{rT^2}\sum_{k=1}^T\|G_{T,k-1}\|_2^2\|{E}%
_{T,k-1}(Z_{T,k}Z_{T,k}^{\prime })-E_r\|_F\bigg\}\leq
Crnh\alpha_X(b+1)^{1/2-1/\gamma}\rightarrow0.
\end{equation*}
On the other hand, by Lemma \ref{la14},
\begin{equation*}
{E}\bigg[\frac{2}{rT^2}\sum_{k=1}^T\{G_{T,k-1}^{\prime }{E}%
_{T,k-1}(Z_{T,k})\}^2\bigg]\leq Crnh\alpha_X(b+1)^{1/2-1/\gamma}\rightarrow0.
\end{equation*}
Then, $I_1=2r^{-1}T^{-2}\sum_{k=1}^T\|G_{T,k-1}\|_2^2+o_p(1). $ % Note that

From Lemma \ref{la15},
\begin{equation*}
\begin{split}
\frac{2}{rT^2}\sum_{k=1}^T\|G_{T,k-1}\|_2^2&=\frac{2}{rT^2}%
\sum_{i=1}^T(T-i)\|Z_{T,i}\|_2^2+\frac{4}{rT^2}\sum_{j=2}^T(T-j)G_{T,j-1}^{%
\prime }Z_{T,j} \\
&=\frac{2}{rT^2}\sum_{i=1}^T(T-i)\|Z_{T,i}\|_2^2+o_p(1).
\end{split}%
\end{equation*}
In order to prove $I_1\xrightarrow[]p1$, it is only need to show $%
2r^{-1}T^{-2}\sum_{i=1}^T(T-i)(\|Z_{T,i}\|_2^2-r)\xrightarrow[]p0. $ Note
that
\begin{equation*}
{E}\bigg\{\frac{2}{rT^2}\sum_{i=1}^T(T-i)(\|Z_{T,i}\|_2^2-r)\bigg\}=0,
\end{equation*}
%\end{spacing}
%\begin{spacing}{1.3}
\noindent it is sufficient to show
\begin{equation*}
\frac{4}{r^2T^4}\bigg[\sum_{i=1}^T(T-i)^2{E}\{(\|Z_{T,i}\|_2^2-r)^2\}+%
\sum_{i\neq j}(T-i)(T-j){E}\{(\|Z_{T,i}\|_2^2-r)(\|Z_{T,j}\|_2^2-r)\}\bigg]%
\rightarrow0,
\end{equation*}
which can be derived from Lemma \ref{la13}. Hence, $I_1\xrightarrow[]p1$.

For the second part, we only need to prove $\sum_{k=1}^T{E}%
(D_{T,k}^4)\rightarrow0$. Note that
\begin{equation*}
D_{T,k}^4\leq Cr^{-2}T^{-4}\{(G_{T,k-1}^{\prime
}Z_{T,k})^4+(\|Z_{T,k}\|_2^2-r)^4\}
\end{equation*}
and
\begin{equation*}
(G_{T,k-1}^{\prime
}Z_{T,k})^4=\sum_{i_1,\ldots,i_4=1}^{k-1}\sum_{j_1,%
\ldots,j_4=1}^rZ_{T,i_1,j_1}Z_{T,i_2,j_2}Z_{T,i_3,j_3}Z_{T,i_4,j_4}Z_{T,k,j_1}Z_{T,k,j_2}Z_{T,k,j_3}Z_{T,k,j_4},
\end{equation*}
where $Z_{T,i,j}$ denotes the $j$th component of $Z_{T,i}$. By the same way
of the Lemma 15 in Francq and Zako\"{\i}an (2007), $r^{-2}h^{-2}{E}%
\{(G_{T,k-1}^{\prime }Z_{T,k})^4\}\leq Ck^2. $ Then, $\sum_{k=1}^T{E}%
(D_{T,k}^4)\leq Ch^2T^{-1}+Cr^2h^{4}T^{-3}\rightarrow0. $ Hence, we complete
the proof. $\hfill \square$

\begin{la}
\label{la18} Under conditions \textrm{(A.1)(i)} and \textrm{(A.2)(iii)},
assume that the eigenvalues of $V_n$ are uniformly bounded away from zero
and infinity. Then $(2r)^{-1/2}\{n\bar{g}(\boldsymbol{\theta}_{0})^{\prime
}V_n^{-1}\bar{g}(\boldsymbol{\theta}_{0})-r\}\xrightarrow[]{d}N(0,1)$
provided that $r^{3/2}h^{-1}\sum_{k=1}^hk\alpha_X(k)^{1-2/\gamma}=o(1)$, $%
rbh^{-1}=o(1)$, $r^2n^{2}h^2\alpha_X(b+1)^{1-2/\gamma}=o(1)$ and $%
rh^3n^{-1}=o(1)$.
\end{la}

\noindent\textsc{Proof}: Note that
\begin{equation*}
\begin{split}
&~(2r)^{-1/2}\{n\bar{g}(\boldsymbol{\theta}_{0})^{\prime }V_n^{-1}\bar{g}(%
\boldsymbol{\theta}_{0})-r\} \\
=&~\frac{n}{Th}S_{T,T}+O_p\bigg\{r^{3/2}h^{-1}\sum_{k=1}^hk\alpha_X(k)^{1-2/%
\gamma}\bigg\}+O_p(r^{1/2}b^{1/2}h^{-1/2}).
\end{split}%
\end{equation*}
Then, by Lemma \ref{la17}, we have $(2r)^{-1/2}\{n\bar{g}(\boldsymbol{\theta}%
_{0})^{\prime }V_n^{-1}\bar{g}(\boldsymbol{\theta}_{0})-r\}\xrightarrow[]{d}%
N(0,1)$. $\hfill \square$

\subsection*{Proof of Theorem 3}

Let $\widehat{\lambda}(\boldsymbol{\theta}_0)=\arg\max_{\lambda\in\widehat{%
\Lambda}_n(\boldsymbol{\theta}_0)}\sum_{q=1}^Q\rho(\lambda^{\prime }\phi_q(%
\boldsymbol{\theta}_0))$. From Lemma \ref{la2}, $\|\bar{\phi}(\boldsymbol{%
\theta}_{0})-\bar{g}(\boldsymbol{\theta}_{0})\|_2=O_p(r^{1/2}Mn^{-1})$.
Hence, $\|\bar{\phi}(\boldsymbol{\theta}_{0})\|_2=O_p(r^{1/2}n^{-1/2})$.
Then, by Lemma \ref{la9}, $\|\widehat{\lambda}(\boldsymbol{\theta}%
_{0})\|_2=O_p(r^{1/2}Mn^{-1/2})$. Meanwhile, $\sup_{1\leq q\leq Q}|\widehat{%
\lambda}(\boldsymbol{\theta}_{0})^{\prime }\phi_q(\boldsymbol{\theta}%
_{0})|=o_p(1)$. Expanding $\max_{\lambda\in\widehat{\Lambda}_n(\boldsymbol{%
\theta}_0)}\sum_{q=1}^Q\rho(\lambda^{\prime }\phi_q(\boldsymbol{\theta}_0))$
around $\lambda=\boldsymbol{\mathbf{0}}$,
\begin{equation*}
\begin{split}
\max_{\lambda\in\widehat{\Lambda}_n(\boldsymbol{\theta}_0)}\sum_{q=1}^Q\rho(%
\lambda^{\prime }\phi_q(\boldsymbol{\theta}_0)) =&~\sum_{q=1}^Q\bigg[%
\rho(0)+\rho_v(0)\widehat{\lambda}(\boldsymbol{\theta}_0)^{\prime }\phi_q(%
\boldsymbol{\theta}_0)+\frac{1}{2}\rho_{vv}(\dot{\lambda}^{\prime }\phi_q(%
\boldsymbol{\theta}_0))\{\widehat{\lambda}(\boldsymbol{\theta}_0)^{\prime
}\phi_q(\boldsymbol{\theta}_0)\}^2\bigg]
\end{split}%
\end{equation*}
where $\dot{\lambda}$ lies on the line joining $\widehat{\lambda}(%
\boldsymbol{\theta}_0)$ and $\boldsymbol{\mathbf{0}}$. On the other hand,
from $\nabla_\lambda\widehat{S}_n(\boldsymbol{\theta}_0,\widehat{\lambda}(%
\boldsymbol{\theta}_0))=\boldsymbol{\mathbf{0}}$, we have
\begin{equation*}
\widehat{\lambda}(\boldsymbol{\theta}_0)=-\bigg\{\frac{1}{Q}%
\sum_{q=1}^Q\rho_{vv}(\ddot{\lambda}^{\prime }\phi_q(\boldsymbol{\theta}%
_0))\phi_q(\boldsymbol{\theta}_0)\phi_q(\boldsymbol{\theta}_0)^{\prime }%
\bigg\}^{-1}\bigg\{\frac{1}{Q}\sum_{q=1}^Q\rho_v(0)\phi_q(\boldsymbol{\theta}%
_0)\bigg\}
\end{equation*}
for some $\ddot{\lambda}$ lies on the line joining $\widehat{\lambda}(%
\boldsymbol{\theta}_0)$ and $\boldsymbol{\mathbf{0}}$. Hence,
\begin{equation*}
\begin{split}
\max_{\lambda\in\widehat{\Lambda}_n(\boldsymbol{\theta}_0)}\sum_{q=1}^Q\rho(%
\lambda^{\prime }\phi_q(\boldsymbol{\theta}_0))=Q\rho(0)-Q\rho^2_v(0)\bar{%
\phi}(\boldsymbol{\theta}_0)^{\prime }\ddot{\Omega}^{-1}\bar{\phi}(%
\boldsymbol{\theta}_0)+\frac{1}{2}Q\rho_v^2(0)\bar{\phi}(\boldsymbol{\theta}%
_0)^{\prime }\ddot{\Omega}^{-1}\dot{\Omega}\ddot{\Omega}^{-1}\bar{\phi}(%
\boldsymbol{\theta}_0) \\
\end{split}%
\end{equation*}
where $\dot{\Omega}=Q^{-1}\sum_{q=1}^Q\rho_{vv}(\dot{\lambda}\phi_q(%
\boldsymbol{\theta}_0))\phi_q(\boldsymbol{\theta}_0)\phi_q(\boldsymbol{\theta%
}_0)^{\prime }$ and $\ddot{\Omega}=Q^{-1}\sum_{q=1}^Q\rho_{vv}(\ddot{\lambda}%
\phi_q(\boldsymbol{\theta}_0))\phi_q(\boldsymbol{\theta}_0)\phi_q(%
\boldsymbol{\theta}_0)^{\prime }$. Then, the generalized empirical
likelihood ratio can be written as
\begin{equation*}
\begin{split}
w_n(\boldsymbol{\theta}_0)=&~2Q\rho_{vv}(0)\bar{\phi}(\boldsymbol{\theta}%
_0)^{\prime }\ddot{\Omega}^{-1}\bar{\phi}(\boldsymbol{\theta}%
_0)-Q\rho_{vv}(0)\bar{\phi}(\boldsymbol{\theta}_0)^{\prime }\ddot{\Omega}%
^{-1}\dot{\Omega}\ddot{\Omega}^{-1}\bar{\phi}(\boldsymbol{\theta}_0) \\
=&~Q\bar{\phi}(\boldsymbol{\theta}_0)^{\prime }\widehat{\Omega}^{-1}(%
\boldsymbol{\theta}_0)\bar{\phi}(\boldsymbol{\theta}_0)+Q\bar{\phi}(%
\boldsymbol{\theta}_0)^{\prime }\{2\rho_{vv}(0)\ddot{\Omega}^{-1}-\widehat{%
\Omega}^{-1}(\boldsymbol{\theta}_0)-\rho_{vv}(0)\ddot{\Omega}^{-1}\dot{\Omega%
}\ddot{\Omega}^{-1}\}\bar{\phi}(\boldsymbol{\theta}_0).
\end{split}%
\end{equation*}
By the same argument of Lemma \ref{la10},
\begin{equation*}
\|\dot{\Omega}-\rho_{vv}(0)\widehat{\Omega}(\boldsymbol{\theta}%
_0)\|_2=O_p(rM^{-1/\gamma}n^{1/\gamma-1/2})=\|\ddot{\Omega}-\rho_{vv}(0)%
\widehat{\Omega}(\boldsymbol{\theta}_0)\|_2.
\end{equation*}
From Lemmas \ref{la3} and $\|M\Omega(\boldsymbol{\theta}_0)-V_n%
\|_2=O(rM^{-1})$, we know the eigenvalues of $M\widehat{\Omega}(\boldsymbol{%
\theta}_0)$ are uniformly bounded away from zero and infinity. Hence,
\begin{equation*}
w_n(\boldsymbol{\theta}_0)=Q\bar{\phi}(\boldsymbol{\theta}_0)^{\prime }%
\widehat{\Omega}^{-1}(\boldsymbol{\theta}_0)\bar{\phi}(\boldsymbol{\theta}%
_0)+O_p(r^2M^{1-1/\gamma}n^{1/\gamma-1/2}).
\end{equation*}
By Lemmas \ref{la2} and \ref{la3}, we have
\begin{equation}
\begin{split}
(2r)^{-1/2}\{w_n(\boldsymbol{\theta}_{0})-r\}=&~(2r)^{-1/2}\{n\bar{g}(%
\boldsymbol{\theta}_{0})^{\prime }V_n^{-1}\bar{g}(\boldsymbol{\theta}%
_{0})-r\}+O_p(r^{5/2}M^{2-2/\gamma}n^{2/\gamma-1}) \\
&+O_p(r^{3/2}M^{1-1/\gamma}n^{1/\gamma-1/2})+O_p(r^{3/2}M^{3/2}n^{-1/2}) \\
&+O_p\bigg\{r^{3/2}M^{-1}\sum_{k=1}^Mk\alpha_X(k)^{1-2/\gamma}\bigg\}.
\label{eq:wn1}
\end{split}%
\end{equation}
%\end{spacing}
%\begin{spacing}{1.36}
\noindent The key step is to show $(2r)^{-1/2}\{n\bar{g}(\boldsymbol{\theta}%
_{0})^{\prime }V_n^{-1}\bar{g}(\boldsymbol{\theta}_{0})-r\}\xrightarrow{d}%
N(0,1)$. In the independent case, the requirements in Lemma \ref{la18} can
be simplified as $rbh^{-1}=o(1)$ and $rh^3n^{-1}=o(1)$. We can pick $b=0$
and $h=1$, then $r=o(n)$. In this case, we can regard $\eta=\infty$. In the
dependent case with $\eta<\infty$, suppose $b\asymp n^{\kappa_1}$ and $%
h\asymp n^{\kappa_2}$, where $0<\kappa_1<\kappa_2<1$. Note that (A.1)'(i),
the requirements in Lemma \ref{la18} turn to
\begin{equation*}
r=o(n^{2\kappa_2/3}),~~r=o(n^{\kappa_2-\kappa_1}),~~r=o(n^{(\eta\kappa_1-2-2%
\kappa_2)/2})~~\text{and}~~r=o(n^{1-3\kappa_2}),
\end{equation*}
where $\eta\kappa_1-2-2\kappa_2>0$ and $1-3\kappa_2>0$. In the following, we
will consider the selection of $(\kappa_1,\kappa_2)$ to satisfy these
inequalities. From $2\kappa_2+2-\eta\kappa_1< 0$, $3\kappa_2-1<0$ and $%
\kappa_1<\kappa_2$, we can get $\frac{2\kappa_2+2}{\eta}<\kappa_1<\kappa_2<%
\frac{1}{3}$. In order to guarantee there exists the solution for above
inequalities in $(0,1)^2$, it is necessary to require $\eta>8$. If $%
8<\eta<\infty$,
\begin{equation*}
\begin{split}
\xi:=&\sup_{\substack{ \frac{2}{\eta-2}<\kappa_2<\frac{1}{3}  \\ \frac{%
2\kappa_2+2}{\eta}<\kappa_1<\kappa_2}}\min\bigg(\frac{2\kappa_2}{3}%
,\kappa_2-\kappa_1, \frac{\eta\kappa_1-2-2\kappa_2}{2},1-3\kappa_2\bigg) \\
=&\sup_{\frac{2}{\eta-2}<\kappa_2<\frac{1}{3}}\min\bigg\{\frac{2\kappa_2}{3},%
\frac{(\eta-2)\kappa_2-2}{\eta+2},1-3\kappa_2\bigg\} \\
=&~\frac{\eta-8}{4\eta+4}1_{(8<\eta<32)}+\frac{2}{11}1_{(32\leq\eta<\infty)}.
\end{split}%
\end{equation*}
In the dependent case with $\eta=\infty$ where $X_t$ is exponentially strong
mixing. The requirements in Lemma \ref{la18} turn to $r^{3/2}h^{-1}=o(1)$, $%
rbh^{-1}=o(1)$ and $rh^3n^{-1}=o(1)$. Then,
\begin{equation*}
r=o(n^{2\kappa_2/3}),~~r=o(n^{\kappa_2-\kappa_1})~~\text{and}%
~~r=o(n^{1-3\kappa_2}).
\end{equation*}
In this setting,
\begin{equation*}
\xi:=\sup_{\substack{ 0<\kappa_2<\frac{1}{3}  \\ 0<\kappa_1<\kappa_2}}\min%
\bigg(\frac{2\kappa_2}{3},\kappa_2-\kappa_1,1-3\kappa_2\bigg)=\frac{2}{11}.
\end{equation*}
Define
\begin{equation*}
\xi=\frac{\eta-8}{4\eta+4}1_{\{8<\eta<32\}}+\frac{2}{11}1_{\{32\leq\eta\leq%
\infty\}}+1_{\{\text{indenpendent data}\}}.
\end{equation*}
Hence, if $r=o(n^{\xi})$, then $(2r)^{-1/2}\{n\bar{g}(\boldsymbol{\theta}%
_{0})^{\prime }V_n^{-1}\bar{g}(\boldsymbol{\theta}_{0})-r\}\xrightarrow[]{d}%
N(0,1)$. If (\ref{eq:cond-ratio}) holds, the other terms in (\ref{eq:wn1})
are $o_p(1)$. We complete the proof of Theorem 3. $\hfill \square$

%\end{spacing}
%\begin{spacing}{1.3}

\subsection*{Proof of Theorem 4}

In order to establish Theorem 4, we need the following lemma.

\begin{la}
\label{la19} For any $\widetilde{\boldsymbol{\theta}}\in\Theta$ and $r\times
r$ matrix $\widehat{V}_n$ such that $\|\widetilde{\boldsymbol{\theta}}-%
\boldsymbol{\theta}_{0}\|_2=O_p(p^{1/2}n^{-1/2})$ and $\|\widehat{V}%
_n-V_n\|_2=o_p(r^{-1/2})$, if $\|\nabla_{\boldsymbol{\theta}} \bar{g}(\dot{%
\boldsymbol{\theta}})-{E}\{\nabla_{\boldsymbol{\theta}} g_t(\boldsymbol{%
\theta}_{0})\}\|_2=o_p(p^{-1/2})$ for any $\dot{\boldsymbol{\theta}}$ with $%
\|\dot{\boldsymbol{\theta}}-\boldsymbol{\theta}_{0}\|_2\leq \|\widetilde{%
\boldsymbol{\theta}}-\boldsymbol{\theta}_{0}\|_2$, and the eigenvalues of $[{%
E}\{\nabla_{\boldsymbol{\theta}} g_t(\boldsymbol{\theta}_{0})\}]^{\prime }[{E%
}\{\nabla_{\boldsymbol{\theta}} g_t(\boldsymbol{\theta}_{0})\}]$ and $V_n$
are uniformly bounded away from zero and infinity, then $(2r)^{-1/2}\{n\bar{g%
}(\widetilde{\boldsymbol{\theta}})^{\prime }\widehat{V}_n^{-1}\bar{g}(%
\widetilde{\boldsymbol{\theta}})-n\bar{g}(\boldsymbol{\theta}_{0})^{\prime
}V_n^{-1}\bar{g}(\boldsymbol{\theta}_{0})\}\xrightarrow{p}0 $ provided that $%
p=o(r^{1/2})$.
\end{la}

\noindent\textsc{Proof}: Note that
\begin{equation*}
\begin{split}
&~~~~(2r)^{-1/2}|n\bar{g}(\widetilde{\boldsymbol{\theta}})^{\prime }\widehat{%
V}_n^{-1}\bar{g}(\widetilde{\boldsymbol{\theta}})-n\bar{g}(\boldsymbol{\theta%
}_{0})^{\prime }V_n^{-1}\bar{g}(\boldsymbol{\theta}_{0})| \\
&=(2r)^{-1/2}|n\bar{g}(\widetilde{\boldsymbol{\theta}})^{\prime }\widehat{V}%
_n^{-1}\bar{g}(\widetilde{\boldsymbol{\theta}})-n\bar{g}(\boldsymbol{\theta}%
_{0})^{\prime }\widehat{V}_n^{-1}\bar{g}(\boldsymbol{\theta}%
_{0})|+(2r)^{-1}|n\bar{g}(\boldsymbol{\theta}_{0})^{\prime }\widehat{V}%
_n^{-1}\bar{g}(\boldsymbol{\theta}_{0})-n\bar{g}(\boldsymbol{\theta}%
_{0})^{\prime }V_n^{-1}\bar{g}(\boldsymbol{\theta}_{0})| \\
&=:I_1+I_2.
\end{split}%
\end{equation*}
We only need to show $I_1\xrightarrow{p}0$ and $I_2\xrightarrow{p}0$.

For $I_1$, by Taylor expansion, $\bar{g}(\widetilde{\boldsymbol{\theta}})=%
\bar{g}(\boldsymbol{\theta}_{0})+\nabla_{\boldsymbol{\theta}} \bar{g}(\dot{%
\boldsymbol{\theta}})\cdot(\widetilde{\boldsymbol{\theta}}-\boldsymbol{\theta%
}_{0})$. Then,
\begin{equation*}
I_1\leq (2r)^{-1/2}|2n(\widetilde{\boldsymbol{\theta}}-\boldsymbol{\theta}%
_{0})^{\prime }\{\nabla_{\boldsymbol{\theta}} \bar{g}(\dot{\boldsymbol{\theta%
}})\}^{\prime }\widehat{V}_n^{-1}\bar{g}(\boldsymbol{\theta}%
_{0})|+(2r)^{-1/2}|n(\widetilde{\boldsymbol{\theta}}-\boldsymbol{\theta}%
_{0})^{\prime }\{\nabla_{\boldsymbol{\theta}} \bar{g}(\dot{\boldsymbol{\theta%
}})\}^{\prime }\widehat{V}_n^{-1}\{\nabla_{\boldsymbol{\theta}}\bar{g}(\dot{%
\boldsymbol{\theta}})\}(\widetilde{\boldsymbol{\theta}}-\boldsymbol{\theta}%
_{0})|.
\end{equation*}
As the eigenvalues of $V_n$ are uniformly bounded away from zero and
infinity, and $\|\widehat{V}_n-V_n\|_2=o_p(r^{-1/2})$, then the eigenvalues
of $\widehat{V}_n$ are uniformly bounded away from zero and infinity
w.p.a.1. Hence,
\begin{equation*}
\begin{split}
&~\|\{\nabla_{\boldsymbol{\theta}} \bar{g}(\dot{\boldsymbol{\theta}}%
)\}^{\prime }\widehat{V}_n^{-1}-[{E}\{\nabla_{\boldsymbol{\theta}} g(%
\boldsymbol{\theta}_{0})\}]^{\prime }V_n^{-1}\|_2 \\
\leq&~\|(\{\nabla_{\boldsymbol{\theta}}\bar{g}(\dot{\boldsymbol{\theta}}%
)\}^{\prime }-[{E}\{\nabla_{\boldsymbol{\theta}} g(\boldsymbol{\theta}%
_{0})\}])\widehat{V}_n^{-1}\|_2+\|[{E}\{\nabla_{\boldsymbol{\theta}} g(%
\boldsymbol{\theta}_{0})\}]^{\prime }(\widehat{V}_n^{-1}-V_n^{-1})\|_2 \\
=&~o_p(p^{-1/2})+o_p(r^{-1/2})=o_p(p^{-1/2}).
\end{split}%
\end{equation*}
On the other hand,
\begin{equation*}
\begin{split}
{E}(\|[{E}\{\nabla_{\boldsymbol{\theta}} g(\boldsymbol{\theta}%
_{0})\}]^{\prime }V_n^{-1}\bar{g}(\boldsymbol{\theta}_{0})\|_2^2)=&~{E}\{%
\text{tr}(V_n^{-1}[{E}\{\nabla_{\boldsymbol{\theta}} g(\boldsymbol{\theta}%
_{0})\}][{E}\{\nabla_{\boldsymbol{\theta}} g(\boldsymbol{\theta}%
_{0})\}]^{\prime }V_n^{-1}\bar{g}(\boldsymbol{\theta}_{0})\bar{g}(%
\boldsymbol{\theta}_{0}))^{\prime }\} \\
=&~n^{-1}\text{tr}([{E}\{\nabla_{\boldsymbol{\theta}} g(\boldsymbol{\theta}%
_{0})\}]^{\prime }V_n^{-1}[{E}\{\nabla_{\boldsymbol{\theta}} g(\boldsymbol{%
\theta}_{0})\}]) \\
\leq&~Cpn^{-1}.
\end{split}%
\end{equation*}
Then,
\begin{equation*}
\|\{\nabla_{\boldsymbol{\theta}}\bar{g}(\dot{\boldsymbol{\theta}})\}^{\prime
}\widehat{V}_n^{-1}\bar{g}(\boldsymbol{\theta}_{0})%
\|_2=O_p(p^{1/2}n^{-1/2})+o_p(r^{1/2}p^{-1/2}n^{-1/2}).
\end{equation*}
Therefore,
\begin{equation*}
I_1\leq O_p(pr^{-1/2})+o_p(1)\xrightarrow{p}0
\end{equation*}
provided that $p=o(r^{1/2})$.

For $I_2$,
\begin{equation*}
I_2=(2r)^{-1/2}|n\bar{g}(\boldsymbol{\theta}_{0})^{\prime }(\widehat{V}%
_n^{-1}-V_n^{-1})\bar{g}(\boldsymbol{\theta}%
_{0})|=O(r^{-1/2}n)o_p(r^{-1/2})O_p(rn^{-1})=o_p(1).
\end{equation*}
Hence, we complete the proof of this lemma. $\hfill \square$ \bigskip

\noindent\textbf{Remark}: This lemma is similar to the Lemma 6.1 of Donald,
Imbens and Newey (2003). However, we work on the operator-norm in
establishing the consistency results, whereas Donald, Imbens and Newey
(2003) employed the Frobenius-norm. The matrix $\widehat{V}_n$ and $%
\widetilde{\boldsymbol{\theta}}$ are the consistency estimators of $V_n$ and
$\boldsymbol{\theta}_{0}$ respectively.

\bigskip

\noindent Here, we begin to establish Theorem 4. From Proposition 1, we know
$\|\widehat{\boldsymbol{\theta}}_{n}-\boldsymbol{\theta}_{0}%
\|_2=O_p(p^{1/2}n^{-1/2})$. By the same argument of the proof of Theorem 3,
we have
\begin{equation*}
w_n(\widehat{\boldsymbol{\theta}}_n)=Q\bar{\phi}(\widehat{\boldsymbol{\theta}%
}_n)^{\prime }\widehat{\Omega}^{-1}(\widehat{\boldsymbol{\theta}}_n)\bar{\phi%
}(\widehat{\boldsymbol{\theta}}_n)+O_p(r^2M^{1-1/\gamma}n^{1/\gamma-1/2}).
\end{equation*}
Note Lemma \ref{la2},
\begin{equation*}
\begin{split}
w_n(\widehat{\boldsymbol{\theta}}_{n})=&~n\bar{g}(\widehat{\boldsymbol{\theta%
}}_{n})^{\prime }\{M\widehat{\Omega}(\widehat{\boldsymbol{\theta}}%
_{n})\}^{-1}\bar{g}(\widehat{\boldsymbol{\theta}}_{n})+O_p(r^{2}M^{1-1/%
\gamma}n^{1/\gamma-1/2})+O_p(rMn^{-1/2})+O_p(rM^2n^{-1}).
\end{split}%
\end{equation*}
By Lemmas \ref{la3} and \ref{la7}, it yields that
\begin{equation*}
\begin{split}
\big\|M\widehat{\Omega}(\widehat{\boldsymbol{\theta}}_{n})-V_n\big\|%
_2&=O_p(r^{1/2}pM^{1/2}n^{-1/2}+rM^{3/2}n^{-1/2})+O_p\bigg\{%
rM^{-1}\sum_{k=1}^Mk\alpha_X(k)^{1-2/\gamma}\bigg\} \\
&=o_p(r^{-1/2}).
\end{split}%
\end{equation*}
Noting Lemma \ref{la11}, for any $\dot{\boldsymbol{\theta}}$ such that $\|%
\dot{\boldsymbol{\theta}}-\boldsymbol{\theta}_{0}\|_2\leq\|\widehat{%
\boldsymbol{\theta}}_{n}-\boldsymbol{\theta}_{0}\|_2=O_p(p^{1/2}n^{-1/2})$,
\begin{equation*}
\|\nabla_{\boldsymbol{\theta}} \bar{g}(\dot{\boldsymbol{\theta}})-{E}%
\{\nabla_{\boldsymbol{\theta}} g_t(\boldsymbol{\theta}_{0})\}%
\|_F=O_p(r^{1/2}p^{3/2}n^{-1/2})=o_p(p^{-1/2}).
\end{equation*}
By Lemma \ref{la19}, we can get $n\bar{g}(\widehat{\boldsymbol{\theta}}%
_{n})^{\prime }\{M\widehat{\Omega}(\widehat{\boldsymbol{\theta}}_{n})\}^{-1}%
\bar{g}(\widehat{\boldsymbol{\theta}}_{n})-n\bar{g}(\boldsymbol{\theta}%
_{0})^{\prime }V_n^{-1}\bar{g}(\boldsymbol{\theta}_{0})=o_p(r^{1/2})$. Then,
by Lemma \ref{la18}, we complete the proof of Theorem 4. $\hfill \square$

\subsection*{Proof of Theorem 5}

We only need to prove that for some $c>1$, $P\{w_n(\widehat{\boldsymbol{%
\theta}}_n)> cr\}\rightarrow1$. To prove this, we use the technique for the
proof of Theorem 1 in \cite{ChangTangWu_2013}. Let
\begin{equation*}
\widetilde{\lambda}=\frac{-\rho_v(0)}{2\rho_{vv}(0)Q^\omega}\frac{e}{%
\max_{1\leq q\leq Q}\|\phi_q(\widehat{\boldsymbol{\theta}}_n)\|_2}
\end{equation*}
where $e$ is a $r$-dimensional vector with unit $L_2$-norm, and $\omega>0$
will be determined later. Then, $\widetilde{\lambda}\in\widehat{\Lambda}_n(%
\widehat{\boldsymbol{\theta}}_n)$ when $Q$ is sufficiently large. Note that $%
\rho_{vv}(0)<0$, by Taylor expansion, we have
\begin{equation*}
\begin{split}
w_n(\widehat{\boldsymbol{\theta}}_n)=&~\frac{2\rho_{vv}(0)}{\rho_v^2(0)}%
\bigg\{Q\rho(0)-\max_{\lambda\in\widehat{\Lambda}_n(\widehat{\boldsymbol{%
\theta}}_n)}\sum_{q=1}^Q\rho(\lambda^{\prime }\phi_q(\widehat{\boldsymbol{%
\theta}}_n))\bigg\} \\
\geq&~\frac{1}{Q^\omega}\sum_{q=1}^Q\frac{e^{\prime }\phi_q(\widehat{%
\boldsymbol{\theta}}_n)}{\max_{1\leq q\leq Q}\|\phi_q(\widehat{\boldsymbol{%
\theta}}_n)\|_2} \\
&~~~~~~~~~~~~~~~~~~~~~-\frac{1}{4Q^{2\omega}\rho_{vv}(0)}\sum_{q=1}^Q\frac{%
\rho_{vv}(\dot{\lambda}\phi_q(\widehat{\boldsymbol{\theta}}_n))e^{\prime
}\phi_q(\widehat{\boldsymbol{\theta}}_n)\phi_q(\widehat{\boldsymbol{\theta}}%
_n)^{\prime }e}{\max_{1\leq q\leq Q}\|\phi_q(\widehat{\boldsymbol{\theta}}%
_n)\|_2^2}
\end{split}%
\end{equation*}
where $\dot{\lambda}$ lies on the jointing line between $\widetilde{\lambda}$
and $\boldsymbol{\mathbf{0}}$. By the definition of $\widetilde{\lambda}$,
we have
\begin{equation*}
\frac{1}{4Q^{2\omega}\rho_{vv}(0)}\sum_{q=1}^Q\frac{\rho_{vv}(\dot{\lambda}%
\phi_q(\widehat{\boldsymbol{\theta}}_n))e^{\prime }\phi_q(\widehat{%
\boldsymbol{\theta}}_n)\phi_q(\widehat{\boldsymbol{\theta}}_n)^{\prime }e}{%
\max_{1\leq q\leq Q}\|\phi_q(\widehat{\boldsymbol{\theta}}_n)\|_2^2}\leq
\frac{1}{2}Q^{1-2\omega}~~\text{w.p.a. 1}.
\end{equation*}
Hence, for any $c>1$,
\begin{equation*}
P\{w_n(\widehat{\boldsymbol{\theta}}_n)\leq cr\}\leq P\bigg\{\sum_{q=1}^Q%
\frac{e^{\prime }\phi_q(\widehat{\boldsymbol{\theta}}_n)}{\max_{1\leq q\leq
Q}\|\phi_q(\widehat{\boldsymbol{\theta}}_n)\|_2}\leq crQ^{\omega}+\frac{1}{2}%
Q^{1-\omega}\bigg\}+o(1).
\end{equation*}
From (A.2)(ii), we have $\|\phi_q(\widehat{\boldsymbol{\theta}}_n)\|_2\leq
r^{1/2}M^{-1}\sum_{t\in B_q}B_n(X_t)$. Then, by Markov inequality,
\begin{equation*}
P\bigg\{\max_{1\leq q\leq Q}\|\phi_q(\widehat{\boldsymbol{\theta}}%
_n)\|_2>(cK)^{-1}r^{1/2}Q^{1/\gamma}(\log Q)^{\epsilon/2}\bigg\}\rightarrow0
\end{equation*}
for each fixed $K>0$, which implies that
\begin{equation*}
P\{w_n(\widehat{\boldsymbol{\theta}}_n)\leq cr\}\leq P\bigg\{%
\sum_{q=1}^Qe^{\prime }\phi_q(\widehat{\boldsymbol{\theta}}_n)\leq
K^{-1}(rQ^{\omega}+Q^{1-\omega})r^{1/2}Q^{1/\gamma}(\log Q)^{\epsilon/2}%
\bigg\}+o(1).
\end{equation*}
Let $rQ^{\omega}=Q^{1-\omega}$, i.e., $Q^{\omega}=Q^{1/2}r^{-1/2}$, then
\begin{equation*}
P\{w_n(\widehat{\boldsymbol{\theta}}_n)\leq cr\}\leq P\bigg\{%
\sum_{q=1}^Qe^{\prime }\phi_q(\widehat{\boldsymbol{\theta}}_n)\leq
2K^{-1}rQ^{1/\gamma+1/2}(\log Q)^{\epsilon/2}\bigg\}+o(1).
\end{equation*}

On the other hand, by Lemma \ref{la2} and (A.1)(iv),
\begin{equation*}
\sum_{q=1}^Qe^{\prime }\phi_q(\widehat{\boldsymbol{\theta}}_n)=Qe^{\prime }%
\bar{g}(\widehat{\boldsymbol{\theta}}_n)+O_p(r^{1/2})=Qe^{\prime }E\{g_t(%
\widehat{\boldsymbol{\theta}}_n)\}+O_p(r^{1/2})+o_p\{Q\Delta_1(r,p)\}.
\end{equation*}
Select $e=E\{g_t(\widehat{\boldsymbol{\theta}}_n)\}/\|E\{g_t(\widehat{%
\boldsymbol{\theta}}_n)\}\|_2$. Then,
\begin{equation*}
\begin{split}
&~P\{w_n(\widehat{\boldsymbol{\theta}}_n)\leq cr\} \\
\leq&~ P\big[\|E\{g_t(\widehat{\boldsymbol{\theta}}_n)\}\|_2\leq
2K^{-1}rQ^{1/\gamma-1/2}(\log
Q)^{\epsilon/2}+O_p(r^{1/2}Q^{-1})+o_p\{\Delta_1(r,p)\}\big]+o(1) \\
\leq&~P\big[\varsigma\leq 2K^{-1}rQ^{1/\gamma-1/2}(\log
Q)^{\epsilon/2}+O_p(r^{1/2}Q^{-1})+o_p\{\Delta_1(r,p)\}\big]+o(1). \\
\end{split}%
\end{equation*}
As $r^2M^{1-2/\gamma}n^{2/\gamma-1}(\log n)^{\epsilon}\varsigma^{-2}=O(1)$, $%
r^{1/2}Mn^{-1}\varsigma^{-1}=o(1)$ and $\Delta_1(r,p)\varsigma^{-1}=O(1)$,
we can choose sufficiently large $K$ to guarantee
\begin{equation*}
P\big[\varsigma\leq 2K^{-1}rQ^{1/\gamma-1/2}(\log
Q)^{\epsilon/2}+O_p(r^{1/2}Q^{-1})+o_p\{\Delta_1(r,p)\}\big]\rightarrow0,
\end{equation*}
which leads to $P\{w_n(\widehat{\boldsymbol{\theta}}_n)\leq cr\}\rightarrow0$
for any $c>1$. Hence, we complete the proof.$\hfill \square$

\section*{Proof of Theorem \protect\ref{tm:pen}}

Let
\begin{equation*}
\widehat{S}_n^{(\mathrm{pe})}(\boldsymbol{\theta},\lambda)=\frac{1}{Q}%
\sum_{q=1}^Q\rho(\lambda^{\prime }\phi_q(\boldsymbol{\theta}%
))+\sum_{j=1}^pp_{\tau}(|\theta_j|)~~%
\mbox{for
any $\btheta\in\Theta$ and
$\lambda\in\widehat{\Lambda}_n(\btheta)$}.
\end{equation*}
Then,
\begin{equation*}
\widehat{\boldsymbol{\theta}}_n^{(\mathrm{pe})}=\arg\min_{\boldsymbol{\theta}%
\in\Theta}\max_{\lambda\in\widehat{\Lambda}_n(\boldsymbol{\theta})}\widehat{S%
}_n^{(\mathrm{pe})}(\boldsymbol{\theta},\lambda)~~\text{and}~~\widehat{%
\boldsymbol{\theta}}_n=\arg\min_{\boldsymbol{\theta}\in\Theta}\max_{\lambda%
\in\widehat{\Lambda}_n(\boldsymbol{\theta})}\widehat{S}_n(\boldsymbol{\theta}%
,\lambda).
\end{equation*}
The following lemma will be used to construct Theorem \ref{tm:pen}.

\begin{la}
Under conditions \textrm{(A.1)}, \textrm{(A.2)} and \textrm{(A.5)}, assume
that the eigenvalues of $V_M$ are uniformly bounded away from zero and
infinity. If \textrm{(\ref{eq:cond1})} holds, $r^2pM^2n^{-1}=o(1)$ and $%
s\tau r^{-1}M^{-1} n=O(1)$, then $\|\widehat{\boldsymbol{\theta}}_n^{(%
\mathrm{pe})}-\boldsymbol{\theta}_0\|_2=O_p(r^{1/2}n^{-1/2})$.
\end{la}

\noindent \textbf{Proof}: Choose $\delta_n=o(r^{-1/2}Q^{-1/\gamma})$ and $%
r^{1/2}Mn^{-1/2}=o(\delta_n)$. Let $\bar{\lambda}=\text{sign}%
\{\rho_v(0)\}\delta_n\bar{\phi}(\widehat{\boldsymbol{\theta}}_n^{(\mathrm{pe}%
)})/\|\bar{\phi}(\widehat{\boldsymbol{\theta}}_n^{(\mathrm{pe})})\|_2$, then
$\bar{\lambda}\in\Lambda_n$ where $\Lambda_n$ is defined in Lemma \ref{la5}.
By Taylor expansion, Lemmas \ref{la4} and \ref{la5}, noting $\rho_{vv}(0)<0$%
, we have
\begin{equation*}
\begin{split}
\widehat{S}_n(\widehat{\boldsymbol{\theta}}_n^{(\mathrm{pe})},\bar{\lambda}%
)=&~\rho(0)+\rho_v(0)\bar{\lambda}^{\prime }\bar{\phi}(\widehat{\boldsymbol{%
\theta}}_n^{(\mathrm{pe})})+\frac{1}{2}\bar{\lambda}^{\prime }\bigg\{\frac{1%
}{Q}\sum_{q=1}^Q\rho_{vv}(\dot{\lambda}^{\prime }\phi_q(\widehat{\boldsymbol{%
\theta}}_n^{(\mathrm{pe})}))\phi_q(\widehat{\boldsymbol{\theta}}_n^{(\mathrm{%
pe})})\phi_q(\widehat{\boldsymbol{\theta}}_n^{(\mathrm{pe})})^{\prime }%
\bigg\}\bar{\lambda} \\
\geq&~ \rho(0)+|\rho_v(0)|\delta_n\|\bar{\phi}(\widehat{\boldsymbol{\theta}}%
_n^{(\mathrm{pe})})\|_2-C\|\bar{\lambda}\|_2^2\cdot O_p(1).
\end{split}%
\end{equation*}
On the other hand,
\begin{equation*}
\widehat{S}_n^{(\mathrm{pe})}(\widehat{\boldsymbol{\theta}}_n^{(\mathrm{pe}%
)},\bar{\lambda})\leq \sup_{\lambda\in\widehat{\Lambda}_n(\widehat{%
\boldsymbol{\theta}}_n^{(\mathrm{pe})})}\widehat{S}_n^{(\mathrm{pe})}(%
\widehat{\boldsymbol{\theta}}_n^{(\mathrm{pe})},\lambda)\leq\sup_{\lambda\in%
\widehat{\Lambda}_n(\boldsymbol{\theta}_0)}\widehat{S}_n^{(\mathrm{pe})}(%
\boldsymbol{\theta}_0,\lambda).
\end{equation*}
By Lemma 7 and (A.5), as $sr^{-1}\tau n=O(1)$,
\begin{equation*}
\begin{split}
\sup_{\lambda\in\widehat{\Lambda}_n(\boldsymbol{\theta}_0)}\widehat{S}_n^{(%
\mathrm{pe})}(\boldsymbol{\theta}_0,\lambda)=&~\sup_{\lambda\in\widehat{%
\Lambda}_n(\boldsymbol{\theta}_0)}\widehat{S}_n(\boldsymbol{\theta}%
_0,\lambda)+\sum_{j=1}^pp_{\tau}(|\theta_{0j}|) \\
=&~\rho(0)+O_p(rMn^{-1}+s\tau)=\rho(0)+O_p(rMn^{-1}).
\end{split}%
\end{equation*}
Note that $\widehat{S}_n^{(\mathrm{pe})}(\boldsymbol{\theta},\lambda)\geq%
\widehat{S}_n(\boldsymbol{\theta},\lambda)$ for any $\boldsymbol{\theta}%
\in\Theta$ and $\lambda\in\widehat{\Lambda}_n(\boldsymbol{\theta})$, it
yields $\|\bar{\phi}(\widehat{\boldsymbol{\theta}}_n^{(\mathrm{pe}%
)})\|_2=O_p(\delta_n)$. Consider any $\varepsilon_n\rightarrow0$ and let $%
\widetilde{\lambda}=\text{sign}\{\rho_v(0)\}\varepsilon_n\bar{\phi}(\widehat{%
\boldsymbol{\theta}}_n^{(\mathrm{pe})})$, then $\|\widetilde{\lambda}%
\|_2=o_p(\delta_n)$. Using the same way above, we can obtain
\begin{equation*}
|\rho_v(0)|\cdot\varepsilon_n\|\bar{\phi}(\widehat{\boldsymbol{\theta}}_n^{(%
\mathrm{pe})})\|_2^2-O_p(1)\cdot\varepsilon_n^2\|\bar{\phi}(\widehat{%
\boldsymbol{\theta}}_n^{(\mathrm{pe})})\|_2^2=O_p(rMn^{-1}).
\end{equation*}
Then, $\varepsilon_n\|\bar{\phi}(\widehat{\boldsymbol{\theta}}_n^{(\mathrm{pe%
})})\|_2^2=O_p(rMn^{-1})$. Thus, $\|\bar{\phi}(\widehat{\boldsymbol{\theta}}%
_n^{(\mathrm{pe})})\|_2=O_p(r^{1/2}M^{1/2}n^{-1/2})$. Following the same
arguments given in the proof of Theorem 1, we can obtain $\|\widehat{%
\boldsymbol{\theta}}_n^{(\mathrm{pe})}-\boldsymbol{\theta}%
_0\|_2=O_p(r^{1/2}n^{-1/2})$. $\hfill\Box$

\bigskip

Here, we begin to prove Theorem \ref{tm:pen}. $\widehat{\boldsymbol{\theta}}%
_n^{(\mathrm{pe})}$ and its Lagrange multiplier $\widehat{\lambda}^{(\mathrm{%
pe})}$ satisfy the score equation
\begin{equation*}
\boldsymbol{\mathbf{0}}=\nabla_\lambda \widehat{S}_n^{(\mathrm{pe})}(%
\widehat{\boldsymbol{\theta}}_n^{(\mathrm{pe})},\widehat{\lambda}^{(\mathrm{%
pe})})=\nabla_\lambda \widehat{S}_n(\widehat{\boldsymbol{\theta}}_n^{(%
\mathrm{pe})},\widehat{\lambda}^{(\mathrm{pe})}).
\end{equation*}
By the implicit theorem (Theorem 9.28 of Rudin, 1976), for all $\boldsymbol{%
\theta}$ in a $\|\cdot\|_2$-neighborhood of $\widehat{\boldsymbol{\theta}}%
_n^{(\mathrm{pe})}$, there is a $\widehat{\lambda}(\theta)$ such that $%
\nabla_\lambda \widehat{S}_n^{(\mathrm{pe})}(\boldsymbol{\theta},\widehat{%
\lambda}(\boldsymbol{\theta}))=\boldsymbol{\mathbf{0}}$ and $\widehat{\lambda%
}(\boldsymbol{\theta})$ is continuously differentiable in $\boldsymbol{\theta%
}$. By the concavity of $\widehat{S}_n^{(\mathrm{pe})}(\boldsymbol{\theta}%
,\lambda)$ with respect to $\lambda$, $\widehat{S}_n^{(\mathrm{pe})}(%
\boldsymbol{\theta},\widehat{\lambda}(\boldsymbol{\theta}))=\max_{\lambda\in%
\widehat{\Lambda}_n(\boldsymbol{\theta})}\widehat{S}_n(\boldsymbol{\theta}%
,\lambda)$. From the envelope theorem,
\begin{equation}  \label{eq:enve}
\begin{split}
\boldsymbol{\mathbf{0}}=&~\nabla_{\boldsymbol{\theta}}\widehat{S}_n^{(%
\mathrm{pe})}(\boldsymbol{\theta},\widehat{\lambda}(\boldsymbol{\theta}))%
\big|_{\boldsymbol{\theta}=\widehat{\boldsymbol{\theta}}_n^{(\mathrm{pe})}}
\\
=&~\frac{1}{Q}\sum_{q=1}^Q\rho_v(\widehat{\lambda}(\widehat{\boldsymbol{%
\theta}}_n^{(\mathrm{pe})})^{\prime }\phi_q(\widehat{\boldsymbol{\theta}}%
_n^{(\mathrm{pe})}))\big\{\nabla_{\boldsymbol{\theta}}\phi_q(\widehat{%
\boldsymbol{\theta}}_n^{(\mathrm{pe})})\big\}^{\prime }\widehat{\lambda}(%
\widehat{\boldsymbol{\theta}}_n^{(\mathrm{pe})})+\sum_{q=1}^Q\nabla_{%
\boldsymbol{\theta}}p_{\tau}(|\theta_j|)\big|_{\boldsymbol{\theta}=\widehat{%
\boldsymbol{\theta}}_n^{(\mathrm{pe})}}. \\
\end{split}%
\end{equation}

For any $\boldsymbol{\theta}$ such that $\|\boldsymbol{\theta}-\boldsymbol{%
\theta}_0\|_2=O_p(r^{1/2}n^{-1/2})$ and $\|\bar{g}(\boldsymbol{\theta}%
)\|_2=O_p(r^{1/2}n^{-1/2})$, define
\begin{equation*}
\mathbf{h}(\boldsymbol{\theta})=\frac{1}{Q}\sum_{q=1}^Q\rho_v(\widehat{%
\lambda}(\boldsymbol{\theta})^{\prime }\phi_q({\boldsymbol{\theta}}))\big\{%
\nabla_{\boldsymbol{\theta}}\phi_q({\boldsymbol{\theta}})\big\}^{\prime }%
\widehat{\lambda}({\boldsymbol{\theta}})+\sum_{q=1}^Q\nabla_{\boldsymbol{%
\theta}}p_{\tau}(|\theta_j|).
\end{equation*}
Write $\mathbf{h}(\boldsymbol{\theta})=(h_1(\boldsymbol{\theta}),\ldots,h_p(%
\boldsymbol{\theta}))^{\prime }$. From Lemma 7, it yields that $\|\widehat{%
\lambda}(\boldsymbol{\theta})\|_2=O_p(r^{1/2}Mn^{-1/2})$ which implies $%
\sup_{1\leq q\leq Q}|\widehat{\lambda}(\boldsymbol{\theta})^{\prime }\phi_q(%
\boldsymbol{\theta})|=o_p(1)$. For each $j=1,\ldots,p$,
\begin{equation*}
\begin{split}
h_j(\boldsymbol{\theta})=&~\frac{1}{Q}\sum_{q=1}^Q\rho_v(0)\widehat{\lambda}(%
\boldsymbol{\theta}_0)^{\prime }\frac{\partial\phi_q(\boldsymbol{\theta}_0)}{%
\partial\theta_j}+\frac{1}{Q}\sum_{q=1}^Q\rho_v(0)\widehat{\lambda}(%
\boldsymbol{\theta}_0)^{\prime }\frac{\partial^2 \phi_q(\boldsymbol{\theta}%
_0)}{\partial\theta_j\partial \boldsymbol{\theta}^{\prime }}(\boldsymbol{%
\theta}-\boldsymbol{\theta}_0)+p_{\tau}^{\prime }(|\theta_j|)\text{sign}%
(\theta_j) \\
&+\text{higher order terms}.
\end{split}%
\end{equation*}
From (A.4), there exists a positive constant $C$ such that $p_{\tau}^{\prime
}(|\theta_j|)\geq C\tau$. On the other hand, as $\tau
(r^{-1}n)^{1/2}M^{-1}\rightarrow\infty$,
\begin{equation*}
\max_{j\notin\mathcal{A}}\bigg|\frac{1}{Q}\sum_{q=1}^Q\rho_v(0)\widehat{%
\lambda}(\boldsymbol{\theta}_0)^{\prime }\frac{\partial\phi_q(\boldsymbol{%
\theta}_0)}{\partial\theta_j}\bigg|= O_p(r^{1/2}Mn^{-1/2})=o_p(\tau).
\end{equation*}
Similarly, we can show
\begin{equation*}
\max_{j\notin\mathcal{A}}\bigg|\frac{1}{Q}\sum_{q=1}^Q\rho_v(0)\widehat{%
\lambda}(\boldsymbol{\theta}_0)^{\prime }\frac{\partial^2 \phi_q(\boldsymbol{%
\theta}_0)}{\partial\theta_j\partial \boldsymbol{\theta}^{\prime }}(%
\boldsymbol{\theta}-\boldsymbol{\theta}_0)\bigg|=o_p(\tau).
\end{equation*}
Hence, $p_{\tau}^{\prime }(|\theta_j|)\text{sign}(\theta_j)$ dominates the
sign of $h_j(\boldsymbol{\theta})$ uniformly for all $j\notin\mathcal{A}$.
If $\widehat{\boldsymbol{\theta}}_n^{(2)}\neq \boldsymbol{\mathbf{0}}$,
there exists some $j\notin\mathcal{A}$ such that $\widehat{\theta}%
_{n,j}\neq0 $. Under our above arguments, we can find
\begin{equation*}
P\big\{h_j(\widehat{\boldsymbol{\theta}}_n^{(\mathrm{pe})})\neq 0\big\}%
\rightarrow1.
\end{equation*}
It is a contradiction. Hence, $\widehat{\boldsymbol{\theta}}_n^{(2)}=%
\boldsymbol{\mathbf{0}}$.

Nextly, we consider the second result. From (\ref{eq:block}), it yields
\begin{equation*}
\begin{split}
&[E\{\nabla_{\boldsymbol{\theta}}g_t(\boldsymbol{\theta}_0)\}]^{\prime
}V_M^{-1}[E\{\nabla_{\boldsymbol{\theta}}g_t(\boldsymbol{\theta}%
_0)\}]([E\{\nabla_{\boldsymbol{\theta}}g_t(\boldsymbol{\theta}_0)\}]^{\prime
}V_M^{-1}V_nV_M^{-1}[E\{\nabla_{\boldsymbol{\theta}}g_t(\boldsymbol{\theta}%
_0)\}])^{-1} \\
&~~~~~~~~~~~~~~~~\times[E\{\nabla_{\boldsymbol{\theta}}g_t(\boldsymbol{\theta%
}_0)\}]^{\prime }V_M^{-1}[E\{\nabla_{\boldsymbol{\theta}}g_t(\boldsymbol{%
\theta}_0)\}]=\left(
\begin{array}{cc}
({\mathbf{S}}_{11}-{\mathbf{S}}_{12}{\mathbf{S}}_{22}^{-1}{\mathbf{S}}%
_{21})^{-1} & \ast \\
\ast & \ast \\
&
\end{array}
\right).
\end{split}%
\end{equation*}
Let
\begin{equation*}
[E\{\nabla_{\boldsymbol{\theta}}g_t(\boldsymbol{\theta}_0)\}]^{\prime
}V_M^{-1}[E\{\nabla_{\boldsymbol{\theta}}g_t(\boldsymbol{\theta}%
_0)\}]([E\{\nabla_{\boldsymbol{\theta}}g_t(\boldsymbol{\theta}_0)\}]^{\prime
}V_M^{-1}V_nV_M^{-1}[E\{\nabla_{\boldsymbol{\theta}}g_t(\boldsymbol{\theta}%
_0)\}])^{-1/2}=\left(
\begin{array}{cc}
{\mathbf{U}} & {\mathbf{V}} \\
{\mathbf{V}}^{\prime } & \ast \\
&
\end{array}
\right),
\end{equation*}
where ${\mathbf{U}}$ is a $s\times s$ symmetric matrix, then ${\mathbf{U}}{%
\mathbf{U}}^{\prime }+{\mathbf{V}}{\mathbf{V}}^{\prime }=({\mathbf{S}}_{11}-{%
\mathbf{S}}_{12}{\mathbf{S}}_{22}^{-1}{\mathbf{S}}_{21})^{-1}$. For any $%
\boldsymbol{\alpha}_n\in\mathbb{R}^s$ such that $\|\boldsymbol{\alpha}%
_n\|_2=1$, define
\begin{equation*}
\widetilde{\boldsymbol{\alpha}}_n=\left(
\begin{array}{c}
{\mathbf{U}}^{\prime } \\
{\mathbf{V}}^{\prime } \\
\end{array}
\right)({\mathbf{S}}_{11}-{\mathbf{S}}_{12}{\mathbf{S}}_{22}^{-1}{\mathbf{S}}%
_{21})^{1/2}\boldsymbol{\alpha}_n.
\end{equation*}
Then,
\begin{equation*}
\widetilde{\boldsymbol{\alpha}}_n^{\prime }\widetilde{\boldsymbol{\alpha}}_n=%
\boldsymbol{\alpha}_n^{\prime }({\mathbf{S}}_{11}-{\mathbf{S}}_{12}{\mathbf{S%
}}_{22}^{-1}{\mathbf{S}}_{21})^{1/2}({\mathbf{U}}{\mathbf{U}}^{\prime }+{%
\mathbf{V}}{\mathbf{V}}^{\prime })({\mathbf{S}}_{11}-{\mathbf{S}}_{12}{%
\mathbf{S}}_{22}^{-1}{\mathbf{S}}_{21})^{1/2}\boldsymbol{\alpha}_n=1.
\end{equation*}
Following the same argument for Proposition 2, we know it still holds for $%
\widehat{\boldsymbol{\theta}}_n^{(\mathrm{pe})}$. Note that
\begin{equation*}
\begin{split}
&\widetilde{\boldsymbol{\alpha}}_n^{\prime }([E\{\nabla_{\boldsymbol{\theta}%
}g_t(\boldsymbol{\theta}_0)\}]^{\prime }V_M^{-1}V_nV_M^{-1}[E\{\nabla_{%
\boldsymbol{\theta}}g_t(\boldsymbol{\theta}_0)\}])^{-1/2}[E\{\nabla_{%
\boldsymbol{\theta}}g_t(\boldsymbol{\theta}_0)\}]^{\prime
}V_M^{-1}[E\{\nabla_{\boldsymbol{\theta}}g_t(\boldsymbol{\theta}_0)\}](%
\widehat{\boldsymbol{\theta}}_n^{(\mathrm{pe})}-\boldsymbol{\theta}_0) \\
&~~~~~~~~~~~~~~~~=\boldsymbol{\alpha}_n^{\prime }({\mathbf{S}}_{11}-{\mathbf{%
S}}_{12}{\mathbf{S}}_{22}^{-1}{\mathbf{S}}_{21})^{-1/2}(\widehat{\boldsymbol{%
\theta}}_n^{(1)}-\boldsymbol{\theta}_0^{(1)}),
\end{split}%
\end{equation*}
then we establish the second result following Proposition 2. $\hfill\Box$
%\section*{Acknowledgements}
%
%We thank the Associate Editor for very constructive comments and suggestions which have improved the presentation of the paper.

%\begin{supplement}
%\sname{Supplement A}\label{suppA} \stitle{Title of the Supplement A}
%\slink[url]{http://www.e-publications.org/ims/support/dowload/imsart-ims.zip}
%\sdescription{Dum esset rex in accubitu suo, nardus mea dedit odorem
%suavitatis. Quoniam confortavit seras portarum tuarum, benedixit
%filiis tuis in te. Qui posuit fines tuos}
%\end{supplement}

\newpage

\begin{table}[tbp]
\caption{Empirical medians of the squared estimation errors ($\times 10^2$)
of the empirical likelihood (EL), the exponential tilting (ET), the
continuous updating (CU) and the optimal GMM for the high dimensional mean
model with $p=\lfloor 10n^{2/15}\rfloor$. }
\label{tb:ex1}%/home/statlab/yichaowu/TangCL/rerunEx2center/summary.R
% The rows headed with EL, ET, CUE
%and GMM are , respectively. The
%magnitudes reported in the table are actual values multiplied by
%$10000$.
\bigskip \centering {\small
\begin{tabular}{c|c||ccc|ccc|ccc}
\hline\hline
& Sample size & \multicolumn{3}{c|}{$n=500$} & \multicolumn{3}{c|}{$n=1000$}
& \multicolumn{3}{c}{$n=2000$} \\ \hline\hline
$L, M$ & \backslashbox{Method}{$\psi$} & $0.1$ & $0.3$ & $0.5$ & $0.1$ & $%
0.3 $ & $0.5$ & $0.1$ & $0.3$ & $0.5$ \\ \hline\hline
(i) & EL & 0.99 & 1.25 & 1.69 & 0.61 & 0.73 & 1.04 & 0.46 & 0.56 & 0.69 \\
& ET & 0.98 & 1.24 & 1.68 & 0.61 & 0.72 & 1.03 & 0.44 & 0.55 & 0.71 \\
& CU & 0.99 & 1.26 & 1.70 & 0.62 & 0.73 & 1.03 & 0.42 & 0.54 & 0.70 \\
& GMM & 1.20 & 1.42 & 1.95 & 0.81 & 0.98 & 1.32 & 0.57 & 0.70 & 0.95 \\
\hline
(ii) & EL & 0.97 & 1.23 & 1.67 & 0.61 & 0.73 & 1.00 & 0.40 & 0.51 & 0.69 \\
& ET & 0.97 & 1.22 & 1.68 & 0.62 & 0.74 & 1.00 & 0.42 & 0.52 & 0.69 \\
& CU & 0.98 & 1.22 & 1.67 & 0.60 & 0.72 & 1.01 & 0.44 & 0.51 & 0.68 \\
& GMM & 1.20 & 1.41 & 1.97 & 0.82 & 0.99 & 1.32 & 0.58 & 0.72 & 0.97 \\
\hline
(iii) & EL & 0.93 & 1.21 & 1.63 & 0.59 & 0.71 & 0.99 & 0.40 & 0.51 & 0.69 \\
& ET & 0.95 & 1.20 & 1.63 & 0.62 & 0.72 & 0.99 & 0.41 & 0.50 & 0.69 \\
& CU & 0.95 & 1.21 & 1.64 & 0.58 & 0.70 & 1.00 & 0.41 & 0.51 & 0.67 \\
& GMM & 1.20 & 1.40 & 1.94 & 0.81 & 0.98 & 1.31 & 0.57 & 0.71 & 0.96 \\
\hline
(iv) & EL & 0.98 & 1.25 & 1.67 & 0.61 & 0.72 & 1.00 & 0.39 & 0.50 & 0.66 \\
& ET & 0.98 & 1.24 & 1.68 & 0.62 & 0.71 & 1.01 & 0.40 & 0.51 & 0.66 \\
& CU & 0.97 & 1.25 & 1.69 & 0.61 & 0.71 & 1.02 & 0.41 & 0.52 & 0.68 \\
& GMM & 1.23 & 1.45 & 1.94 & 0.85 & 1.03 & 1.32 & 0.57 & 0.72 & 0.98 \\
\hline
(v) & EL & 0.96 & 1.21 & 1.64 & 0.57 & 0.71 & 0.98 & 0.37 & 0.50 & 0.66 \\
& ET & 0.98 & 1.20 & 1.68 & 0.57 & 0.69 & 0.99 & 0.39 & 0.51 & 0.65 \\
& CU & 0.97 & 1.22 & 1.68 & 0.58 & 0.70 & 1.00 & 0.40 & 0.51 & 0.68 \\
& GMM & 1.21 & 1.42 & 1.92 & 0.85 & 1.01 & 1.30 & 0.57 & 0.72 & 0.97 \\
\hline\hline
\end{tabular}
}
\end{table}

\begin{table}[tbp]
\caption{Empirical medians of the squared estimation errors ($\times 10^2$)
of the empirical likelihood (EL), the exponential tilting (ET), the
continuous updating (CU) and the optimal GMM for the high dimensional mean
model with $p=\lfloor 12n^{2/15}\rfloor$. }
\label{tb:ex1}%/home/statlab/yichaowu/TangCL/rerunEx2center/summary.R
\bigskip \centering {\small
\begin{tabular}{c|c||ccc|ccc|ccc}
\hline\hline
& Sample size & \multicolumn{3}{c|}{$n=500$} & \multicolumn{3}{c|}{$n=1000$}
& \multicolumn{3}{c}{$n=2000$} \\ \hline\hline
$L, M$ & \backslashbox{Method}{$\psi$} & $0.1$ & $0.3$ & $0.5$ & $0.1$ & $%
0.3 $ & $0.5$ & $0.1$ & $0.3$ & $0.5$ \\ \hline\hline
(i) & EL & 0.79 & 1.02 & 1.39 & 0.44 & 0.57 & 0.77 & 0.30 & 0.35 & 0.44 \\
& ET & 0.80 & 0.99 & 1.39 & 0.45 & 0.58 & 0.78 & 0.33 & 0.36 & 0.45 \\
& CU & 0.79 & 1.01 & 1.40 & 0.46 & 0.57 & 0.79 & 0.32 & 0.34 & 0.43 \\
& GMM & 1.15 & 1.49 & 2.00 & 0.79 & 1.01 & 1.30 & 0.55 & 0.73 & 0.89 \\
\hline
(ii) & EL & 0.79 & 0.97 & 1.36 & 0.43 & 0.54 & 0.75 & 0.28 & 0.35 & 0.45 \\
& ET & 0.78 & 0.96 & 1.38 & 0.41 & 0.53 & 0.74 & 0.30 & 0.36 & 0.44 \\
& CU & 0.77 & 0.98 & 1.39 & 0.40 & 0.52 & 0.75 & 0.26 & 0.33 & 0.42 \\
& GMM & 1.15 & 1.48 & 1.99 & 0.81 & 1.01 & 1.31 & 0.56 & 0.71 & 0.90 \\
\hline
(iii) & EL & 0.76 & 0.93 & 1.33 & 0.41 & 0.53 & 0.75 & 0.28 & 0.33 & 0.43 \\
& ET & 0.77 & 0.94 & 1.35 & 0.40 & 0.52 & 0.73 & 0.24 & 0.33 & 0.43 \\
& CU & 0.74 & 0.92 & 1.36 & 0.38 & 0.51 & 0.74 & 0.25 & 0.32 & 0.40 \\
& GMM & 1.14 & 1.47 & 1.99 & 0.80 & 1.01 & 1.29 & 0.54 & 0.70 & 0.89 \\
\hline
(iv) & EL & 0.77 & 0.98 & 1.35 & 0.41 & 0.52 & 0.71 & 0.26 & 0.32 & 0.42 \\
& ET & 0.78 & 0.98 & 1.36 & 0.41 & 0.53 & 0.72 & 0.27 & 0.33 & 0.44 \\
& CU & 0.79 & 0.97 & 1.37 & 0.40 & 0.51 & 0.71 & 0.28 & 0.36 & 0.43 \\
& GMM & 1.15 & 1.48 & 1.99 & 0.80 & 1.02 & 1.31 & 0.56 & 0.72 & 0.91 \\
\hline
(v) & EL & 0.75 & 0.95 & 1.34 & 0.40 & 0.51 & 0.69 & 0.23 & 0.32 & 0.40 \\
& ET & 0.77 & 0.96 & 1.35 & 0.38 & 0.50 & 0.71 & 0.26 & 0.30 & 0.42 \\
& CU & 0.76 & 0.97 & 1.36 & 0.39 & 0.49 & 0.70 & 0.26 & 0.32 & 0.41 \\
& GMM & 1.14 & 1.46 & 1.98 & 0.79 & 1.01 & 1.27 & 0.55 & 0.72 & 0.90 \\
\hline\hline
\end{tabular}
}
\end{table}

\begin{table}[tbp]
\caption{Empirical medians of the squared estimation errors ($\times 10^2$)
of the empirical likelihood (EL), the penalized empirical likelihood (PEL),
the exponential tilting (ET), the penalized exponential tilting (PET), the
continuous updating (CU), the penalized continuous updating (PCU) and the
optimal GMM for the high dimensional generalized linear model with $%
p=\lfloor 5n^{2/15}\rfloor$. }
\label{tb:ex2}%/home/statlab/yichaowu/TangCL/rerunEx2center/summary.R
\bigskip \centering {\small
\begin{tabular}{c|c||ccc|ccc|ccc}
\hline\hline
& Sample size & \multicolumn{3}{c|}{$n=500$} & \multicolumn{3}{c|}{$n=1000$}
& \multicolumn{3}{c}{$n=2000$} \\ \hline\hline
$L, M$ & \backslashbox{Method}{$\psi$} & $0.1$ & $0.3$ & $0.5$ & $0.1$ & $%
0.3 $ & $0.5$ & $0.1$ & $0.3$ & $0.5$ \\ \hline\hline
(i) & EL & 3.16 & 3.59 & 3.67 & 3.13 & 3.42 & 3.62 & 2.85 & 2.89 & 3.00 \\
& PEL & 0.95 & 1.05 & 1.10 & 0.84 & 0.93 & 1.05 & 0.74 & 0.78 & 0.89 \\
& ET & 3.05 & 3.19 & 3.38 & 3.03 & 3.10 & 3.20 & 2.70 & 2.83 & 2.92 \\
& PET & 0.84 & 0.95 & 1.00 & 0.82 & 0.89 & 0.95 & 0.67 & 0.80 & 0.83 \\
& CU & 3.77 & 4.01 & 4.35 & 3.73 & 3.94 & 4.16 & 3.18 & 3.30 & 3.46 \\
& PCU & 1.14 & 1.22 & 1.26 & 1.05 & 1.18 & 1.22 & 0.95 & 1.01 & 1.05 \\
& GMM & 6.72 & 7.06 & 7.12 & 6.30 & 6.43 & 6.60 & 5.83 & 5.88 & 5.92 \\
\hline
(ii) & EL & 3.11 & 3.32 & 3.48 & 3.07 & 3.10 & 3.35 & 2.63 & 2.86 & 2.93 \\
& PEL & 0.95 & 1.05 & 1.14 & 0.83 & 0.89 & 1.00 & 0.77 & 0.84 & 0.87 \\
& ET & 3.05 & 3.10 & 3.11 & 2.85 & 2.95 & 3.00 & 2.37 & 2.53 & 2.76 \\
& PET & 0.89 & 0.93 & 0.95 & 0.77 & 0.85 & 0.89 & 0.74 & 0.82 & 0.84 \\
& CU & 3.26 & 3.45 & 3.77 & 3.19 & 3.29 & 3.62 & 2.79 & 3.02 & 3.24 \\
& PCU & 1.00 & 1.10 & 1.14 & 0.95 & 1.00 & 1.07 & 0.84 & 0.95 & 0.98 \\
& GMM & 6.80 & 7.11 & 7.34 & 6.33 & 6.55 & 6.63 & 5.83 & 5.92 & 5.96 \\
\hline
(iii) & EL & 2.81 & 2.90 & 2.97 & 2.76 & 2.85 & 2.92 & 2.53 & 2.70 & 2.85 \\
& PEL & 0.77 & 0.84 & 0.93 & 0.71 & 0.77 & 0.85 & 0.63 & 0.70 & 0.77 \\
& ET & 2.65 & 2.70 & 2.76 & 2.45 & 2.65 & 2.68 & 2.21 & 2.41 & 2.51 \\
& PET & 0.84 & 0.89 & 0.93 & 0.71 & 0.82 & 0.86 & 0.66 & 0.71 & 0.77 \\
& CU & 3.16 & 3.42 & 3.56 & 3.11 & 3.32 & 3.46 & 2.70 & 2.83 & 3.18 \\
& PCU & 0.95 & 1.05 & 1.10 & 0.84 & 0.95 & 1.00 & 0.77 & 0.84 & 0.89 \\
& GMM & 6.76 & 7.06 & 7.23 & 6.28 & 6.50 & 6.60 & 5.81 & 5.86 & 5.90 \\
\hline
(iv) & EL & 2.24 & 2.35 & 2.43 & 2.19 & 2.28 & 2.39 & 1.97 & 2.12 & 2.26 \\
& PEL & 0.71 & 0.75 & 0.84 & 0.63 & 0.65 & 0.71 & 0.55 & 0.63 & 0.68 \\
& ET & 2.14 & 2.28 & 2.33 & 2.07 & 2.21 & 2.26 & 1.92 & 2.05 & 2.12 \\
& PET & 0.63 & 0.72 & 0.79 & 0.60 & 0.66 & 0.77 & 0.55 & 0.57 & 0.66 \\
& CU & 2.86 & 2.90 & 3.11 & 2.74 & 2.81 & 3.03 & 2.28 & 2.43 & 2.74 \\
& PCU & 0.90 & 0.95 & 1.00 & 0.77 & 0.84 & 0.95 & 0.63 & 0.73 & 0.82 \\
& GMM & 7.16 & 7.24 & 7.30 & 6.30 & 6.50 & 6.59 & 5.85 & 5.92 & 6.03 \\
\hline
(v) & EL & 2.19 & 2.26 & 2.32 & 2.07 & 2.14 & 2.17 & 1.61 & 1.76 & 1.84 \\
& PEL & 0.63 & 0.71 & 0.77 & 0.55 & 0.63 & 0.70 & 0.45 & 0.50 & 0.57 \\
& ET & 1.79 & 1.94 & 2.07 & 1.64 & 1.79 & 2.00 & 1.38 & 1.52 & 1.70 \\
& PET & 0.50 & 0.55 & 0.65 & 0.45 & 0.51 & 0.63 & 0.32 & 0.45 & 0.52 \\
& CU & 2.74 & 2.83 & 3.08 & 2.53 & 2.68 & 2.90 & 2.10 & 2.35 & 2.83 \\
& PCU & 0.84 & 0.87 & 0.95 & 0.71 & 0.80 & 0.84 & 0.55 & 0.68 & 0.80 \\
& GMM & 7.04 & 7.13 & 7.20 & 6.29 & 6.47 & 6.55 & 5.84 & 5.92 & 5.97 \\
\hline\hline
\end{tabular}
}
\end{table}

\begin{table}[tbp]
\caption{Empirical medians of the squared estimation errors ($\times 10^2$)
of the empirical likelihood (EL), the penalized empirical likelihood (PEL),
the exponential tilting (ET), the penalized exponential tilting (PET), the
continuous updating (CU), the penalized continuous updating (PCU) and the
optimal GMM for the high dimensional generalized linear model with $%
p=\lfloor 6n^{2/15}\rfloor$. }
\label{tb:ex2}%/home/statlab/yichaowu/TangCL/rerunEx2center/summary.R
\bigskip \centering {\small
\begin{tabular}{c|c||ccc|ccc|ccc}
\hline\hline
& Sample size $n$ & \multicolumn{3}{c|}{$n=500$} & \multicolumn{3}{c|}{$%
n=1000$} & \multicolumn{3}{c}{$n=2000$} \\ \hline\hline
$L, M$ & \backslashbox{Method}{$\psi$} & $0.1$ & $0.3$ & $0.5$ & $0.1$ & $%
0.3 $ & $0.5$ & $0.1$ & $0.3$ & $0.5$ \\ \hline\hline
(i) & EL & 2.63 & 2.75 & 2.82 & 2.48 & 2.60 & 2.76 & 2.31 & 2.42 & 2.54 \\
& PEL & 0.81 & 0.84 & 0.87 & 0.75 & 0.78 & 0.81 & 0.68 & 0.72 & 0.77 \\
& ET & 2.46 & 2.60 & 2.71 & 2.35 & 2.51 & 2.62 & 2.19 & 2.23 & 2.34 \\
& PET & 0.74 & 0.79 & 0.85 & 0.69 & 0.76 & 0.79 & 0.64 & 0.69 & 0.71 \\
& CU & 2.83 & 2.96 & 3.04 & 2.75 & 2.83 & 2.95 & 2.62 & 2.71 & 2.87 \\
& PCU & 0.88 & 0.92 & 0.95 & 0.82 & 0.88 & 0.93 & 0.76 & 0.80 & 0.86 \\
& GMM & 5.84 & 5.93 & 6.09 & 5.47 & 5.53 & 5.62 & 4.73 & 4.83 & 4.99 \\
\hline
(ii) & EL & 2.50 & 2.79 & 2.85 & 2.42 & 2.58 & 2.66 & 2.18 & 2.30 & 2.39 \\
& PEL & 0.73 & 0.83 & 0.86 & 0.66 & 0.78 & 0.80 & 0.56 & 0.66 & 0.71 \\
& ET & 2.35 & 2.44 & 2.59 & 2.28 & 2.37 & 2.45 & 2.15 & 2.23 & 2.32 \\
& PET & 0.65 & 0.73 & 0.79 & 0.62 & 0.67 & 0.75 & 0.57 & 0.61 & 0.66 \\
& CU & 2.79 & 2.82 & 2.94 & 2.61 & 2.74 & 2.81 & 2.58 & 2.67 & 2.73 \\
& PCU & 0.82 & 0.87 & 0.93 & 0.78 & 0.86 & 0.89 & 0.70 & 0.73 & 0.81 \\
& GMM & 5.89 & 5.96 & 6.13 & 5.45 & 5.56 & 5.64 & 4.70 & 4.81 & 5.05 \\
\hline
(iii) & EL & 2.37 & 2.65 & 2.69 & 2.09 & 2.22 & 2.42 & 1.96 & 2.12 & 2.31 \\
& PEL & 0.68 & 0.76 & 0.80 & 0.58 & 0.70 & 0.73 & 0.46 & 0.60 & 0.67 \\
& ET & 2.20 & 2.36 & 2.44 & 1.98 & 2.15 & 2.27 & 1.90 & 1.99 & 2.08 \\
& PET & 0.61 & 0.72 & 0.76 & 0.45 & 0.57 & 0.63 & 0.42 & 0.50 & 0.57 \\
& CU & 2.64 & 2.80 & 2.84 & 2.59 & 2.68 & 2.76 & 2.50 & 2.58 & 2.63 \\
& PCU & 0.76 & 0.84 & 0.89 & 0.74 & 0.79 & 0.82 & 0.65 & 0.73 & 0.77 \\
& GMM & 5.86 & 5.92 & 6.10 & 5.42 & 5.50 & 5.61 & 4.68 & 4.75 & 5.02 \\
\hline
(iv) & EL & 2.21 & 2.30 & 2.41 & 1.99 & 2.16 & 2.24 & 1.73 & 1.80 & 1.97 \\
& PEL & 0.58 & 0.67 & 0.72 & 0.51 & 0.58 & 0.63 & 0.40 & 0.46 & 0.53 \\
& ET & 1.92 & 2.04 & 2.18 & 1.82 & 1.89 & 1.97 & 1.46 & 1.62 & 1.83 \\
& PET & 0.54 & 0.59 & 0.61 & 0.46 & 0.52 & 0.54 & 0.34 & 0.39 & 0.48 \\
& CU & 2.58 & 2.63 & 2.74 & 2.31 & 2.46 & 2.70 & 2.23 & 2.31 & 2.47 \\
& PCU & 0.68 & 0.78 & 0.82 & 0.65 & 0.73 & 0.78 & 0.57 & 0.66 & 0.71 \\
& GMM & 6.04 & 6.10 & 6.16 & 5.60 & 5.64 & 5.78 & 4.70 & 4.79 & 5.01 \\
\hline
(v) & EL & 2.16 & 2.25 & 2.32 & 1.84 & 2.02 & 2.15 & 1.63 & 1.76 & 1.92 \\
& PEL & 0.56 & 0.60 & 0.67 & 0.47 & 0.52 & 0.61 & 0.37 & 0.42 & 0.52 \\
& ET & 1.90 & 1.96 & 2.06 & 1.71 & 1.76 & 1.84 & 1.34 & 1.43 & 1.59 \\
& PET & 0.54 & 0.58 & 0.61 & 0.43 & 0.46 & 0.51 & 0.30 & 0.32 & 0.36 \\
& CU & 2.53 & 2.58 & 2.66 & 2.24 & 2.44 & 2.61 & 2.17 & 2.21 & 2.29 \\
& PCU & 0.68 & 0.72 & 0.77 & 0.56 & 0.67 & 0.75 & 0.56 & 0.60 & 0.67 \\
& GMM & 5.93 & 6.08 & 6.15 & 5.52 & 5.59 & 5.71 & 4.64 & 4.72 & 4.88 \\
\hline\hline
\end{tabular}
}
\end{table}

%\end{spacing}

\end{document}